%% file: ms.tex
\documentclass[opre]{informs3} %

\TheoremsNumberedThrough
\ECRepeatTheorems
\EquationsNumberedThrough
\usepackage{endnotes}
\let\footnote=\endnote
\let\enotesize=\normalsize
\def\notesname{Endnotes}%
\def\makeenmark{$^{\theenmark}$}
\def\enoteformat{\rightskip0pt\leftskip0pt\parindent=1.75em
  \leavevmode\llap{\theenmark.\enskip}}

\usepackage{natbib}
 \bibpunct[, ]{(}{)}{,}{a}{}{,}%
 \def\bibfont{\small}%
 \def\bibsep{\smallskipamount}%
 \def\bibhang{24pt}%
 \def\BIBand{and}%

\input{math_commands.tex}

\usepackage{booktabs}
\usepackage{float}
\usepackage{adjustbox}
\usepackage{subcaption}
\usepackage{hyperref}
\usepackage{xcolor}
\usepackage[ruled,vlined]{algorithm2e}

\def\finalver{0} %

\if\finalver1
\newcommand{\felix}[1]{}
\newcommand{\vinodquote}[1]{}
\newcommand{\sergey}[1]{}
\newcommand{\TODO}[1]{}
\newcommand{\todoct}[1]{}
\newcommand{\todopl}[1]{}
\newcommand{\todoivg}[1]{}
\newcommand{\yujia}[1]{}
\newcommand{\nic}[1]{}
\newcommand{\oriol}[1]{}
\newcommand{\ravi}[1]{}
\newcommand{\pengming}[1]{}
\newcommand{\ira}[1]{}

\else
\newcommand{\felix}[1]{\textcolor{magenta}{\textit{\textbf{[Comment by Felix: #1]}}}}
\newcommand{\vinodquote}[1]{\textcolor{cyan}{\textit{\textbf{[Vinod quote: #1]}}}}
\newcommand{\sergey}[1]{\textcolor{green}{\textit{\textbf{[TODO Sergey: #1]}}}}
\newcommand{\TODO}[1]{\textcolor{red}{[TODO: #1]}}
\newcommand{\todoct}[1]{\textcolor{blue}{[ctjandra: #1]}}
\newcommand{\todopl}[1]{\textcolor{blue}{[pawell: #1]}}
\newcommand{\todoivg}[1]{\textcolor{blue}{[ingridvg: #1]}}
\newcommand{\yujia}[1]{\textcolor{orange}{[yujiali: #1]}}
\newcommand{\nic}[1]{\textcolor{teal}{[sonnerat: #1]}}
\newcommand{\oriol}[1]{\textcolor{pink}{[Oriol: #1]}}
\newcommand{\ravi}[1]{\textcolor{blue}{[Ravi: #1]}}
\newcommand{\pengming}[1]{\textcolor{blue}{[Pengming: #1]}}

\newcommand{\ira}[1]{\textcolor{purple}{[iraktena: #1]}}

\fi

\input{bod_macros}

\begin{document}

\RUNAUTHOR{Nair et al.}

\RUNTITLE{Solving Mixed Integer Programs Using Deep Neural Networks}

\TITLE{Solving Mixed Integer Programs Using Neural Networks}

\ARTICLEAUTHORS{%
\AUTHOR{Vinod Nair$^{*\dagger 1}$, Sergey Bartunov$^{*1}$, Felix Gimeno$^{*1}$, Ingrid von Glehn$^{*1}$, Pawel Lichocki$^{*2}$, Ivan Lobov$^{*1}$, Brendan O'Donoghue$^{*1}$, Nicolas Sonnerat$^{*1}$, Christian Tjandraatmadja$^{*2}$, Pengming Wang$^{*1}$, Ravichandra Addanki$^{1}$, Tharindi Hapuarachchi$^{1}$, Thomas Keck$^{1}$, James Keeling$^{1}$, Pushmeet Kohli$^{1}$, Ira Ktena$^{1}$, Yujia Li$^{1}$, Oriol Vinyals$^{1}$, Yori Zwols$^{1}$}
\AFF{$^{1}$DeepMind, $^{2}$Google Research \footnotetext{$^{*}$Equal contributors, listed alphabetically by last name after the lead author.} \footnotetext{$^{\dagger}$Corresponding author, email: vinair@google.com.}}
} %

\ABSTRACT{
Mixed Integer Programming (MIP) solvers rely on an array of sophisticated heuristics developed with decades of research to solve large-scale MIP instances encountered in practice. Machine learning offers to automatically construct better heuristics from data by exploiting shared structure among instances in the data. This paper applies learning to the two key sub-tasks of a MIP solver, generating a high-quality joint variable assignment, and bounding the gap in objective value between that assignment and an optimal one. Our approach constructs two corresponding neural network-based components, \emph{Neural Diving} and \emph{Neural Branching}, to use in a base MIP solver such as SCIP. Neural Diving learns a deep neural network to generate multiple partial assignments for its integer variables, and the resulting smaller MIPs for un-assigned variables are solved with SCIP to construct high quality joint assignments. Neural Branching learns a deep neural network to make variable selection decisions in branch-and-bound to bound the objective value gap with a small tree. This is done by imitating a new variant of Full Strong Branching we propose that scales to large instances using GPUs. We evaluate our approach on diverse real-world datasets, including two Google production datasets and MIPLIB, by training separate neural networks on each. Most instances in all the datasets combined have $10^3-10^6$ variables and constraints after presolve, which is significantly larger than previous learning approaches.  Comparing solvers with respect to primal-dual gap averaged over a held-out set of instances at large time limits, the learning-augmented SCIP achieves $1.5\times$, $2\times$, and $10^4\times$ better gap on three out of the five datasets with the largest MIPs, achieves a $10\%$ gap $5\times$ faster on a fourth one, and matches SCIP on the fifth. To the best of our knowledge, ours is the first learning approach to demonstrate such large improvements over SCIP on both large-scale real-world application datasets and MIPLIB.
}


\HISTORY{First version December 2020.}

\maketitle

\input{introduction_new.tex}

\input{background.tex} %
\input{rep_and_arch.tex}

\input{datasets.tex}

\input{evaluation.tex}

\input{supermip.tex}
\input{supermip_results.tex}
\input{deepbrancher.tex}

\input{deepbrancher_results.tex}

\input{joint_eval_results.tex}

\input{related_work.tex}

\input{conclusions.tex}

\input{acknowledgments.tex}

\newpage


\input{ms.bbl}
\newpage

\input{appendix.tex}

\end{document}

%% file: math_commands.tex
\usepackage{amsmath,amsfonts,bm}

\def\ceil#1{\lceil #1 \rceil}
\def\floor#1{\lfloor #1 \rfloor}
\def\1{\bm{1}}

\DeclareMathAlphabet{\mathsfit}{\encodingdefault}{\sfdefault}{m}{sl}
\SetMathAlphabet{\mathsfit}{bold}{\encodingdefault}{\sfdefault}{bx}{n}

\newcommand{\R}{\mathbb{R}}

%% file: bod_macros.tex
\newcommand{\eg}{{\it e.g.}}
\newcommand{\ie}{{\it i.e.}}

\newcommand{\BA}{\begin{array}}
\newcommand{\EA}{\end{array}}

\newcommand{\BIT}{\begin{itemize}}
\newcommand{\EIT}{\end{itemize}}

\newcommand{\reals}{{\mathbb{R}}} %

\newcommand{\Ic}{\mathcal{I}}
\newcommand{\Xc}{\mathcal{X}}

%% file: introduction_new.tex
\section{Introduction}
\label{sec:intro}

\begin{figure}[t]
    \includegraphics[trim={3cm 2.5cm 3cm 2.5cm}, clip, width=0.98\textwidth]{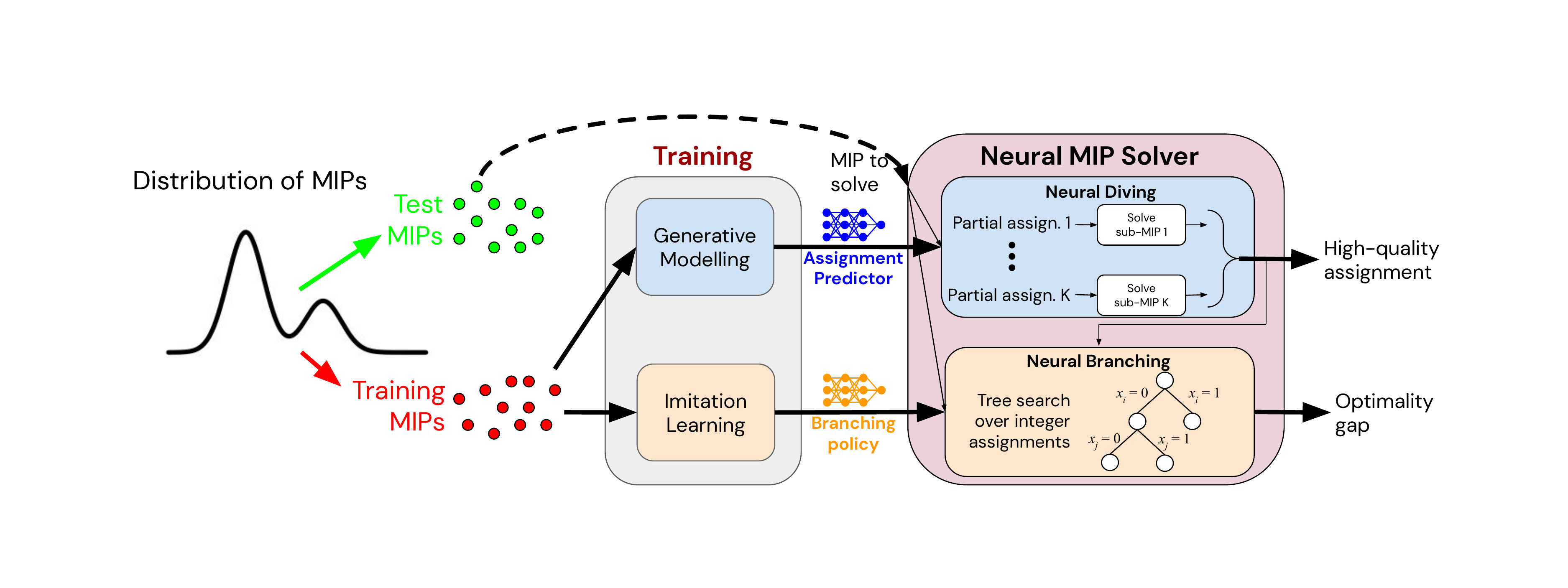}
    \caption{\small{Our approach constructs two neural network-based components to use in a MIP solver, \emph{Neural Diving} and \emph{Neural Branching}, and combine them to produce a \emph{Neural Solver} customized to a given MIP dataset. 
    }}
    \label{fig:overview}
\end{figure}

Mixed Integer Programs (MIPs) are a class of NP-hard problems where the goal is to minimize a linear objective subject to linear constraints, with some or all of the variables constrained to be integer-valued \citep{wolsey1998integer, karp1972reducibility}. They have enjoyed widespread adoption in a broad range of applications such as capacity planning, resource allocation, bin packing, etc. \citep{taha2014integer, junger200950, sierksma2015linear}. Significant research and engineering effort has gone into developing practical solvers, such as SCIP \citep{gamrath2020scip}, CPLEX \citep{cplex2019}, Gurobi \citep{gurobi2020}, and Xpress \citep{xpress}. These solvers use sophisticated heuristics to direct the search process for solving a MIP. A solver's performance on a given application depends crucially on how well its heuristics suit that application.

In this paper we show that machine learning can be used to automatically construct effective heuristics from a dataset of MIP instances. A compelling use case for this arises often in practice where an application requires solving a large set of instances of the same high-level semantic problem with different problem parameters. Examples of such ``homogenous'' datasets in this paper are 1) optimizing the choice of power plants on an electric grid to meet demand \citep{oneill17}, where grid topology remains the same while demand, renewable generation, etc. vary across instances, and 2) solving a packing problem in Google's production system where the semantics of ``items'' and ``bins'' to be packed remain mostly the same but their sizes fluctuate across instances. Even a ``heterogeneous'' dataset that combines many semantically different problems, such as MIPLIB 2017 \citep{gleixner2019miplib}, can have structure across instances that can be used to learn better heuristics, as we shall show.  Off-the-shelf MIP solvers cannot automatically construct heuristics to exploit such structure. In challenging applications users may rely on an expert to hand-design such heuristics, or forego potentially large performance improvements. Machine learning offers the possibility of large improvements without needing application-specific expertise.

We demonstrate that machine learning can construct heuristics customized to a given dataset that significantly outperform the classical ones used in a MIP solver, specifically the state-of-the-art non-commercial solver SCIP 7.0.1 \citep{gamrath2020scip}. Our approach applies learning to the two key sub-tasks of a MIP solver: a) output an assignment of values to all variables that satisfy the constraints (if such an assignment exists), and b) prove a bound for the gap in objective value between that assignment and an optimal one. They define the main components of our approach, \emph{Neural Diving} and \emph{Neural Branching} (see figure~\ref{fig:overview}).

\paragraph{\textbf{Neural Diving:}} This component finds high quality joint variable assignments. It is an instance of a \emph{primal heuristic} \citep{berthold2006thesis}, a class of search heuristics that have been identified as key to effective MIP solvers \citep{berthold2013primal}. We train a deep neural network to produce multiple partial assignments of the integer variables of the input MIP. The remaining unassigned variables define smaller `sub-MIPs', which are solved using an off-the-shelf MIP solver (e.g., SCIP) to complete the assignments. The sub-MIPs can be solved in parallel if the compute budget allows. The model is trained to give higher probability to feasible assignments that have better objective values, using training examples collected offline with an off-the-shelf solver. It learns on all available feasible assignments instead of only the optimal ones, and does not necessarily require optimal assignments, which can be expensive to collect. %

\paragraph{\textbf{Neural Branching:}} This component is mainly used to bound the gap between the objective value of the best assignment and an optimal one. MIP solvers use a form of tree search over integer variables called \emph{branch-and-bound} \citep{Land1960BranchAndBound} (see section \ref{sec:background}), which progressively tightens the bound and helps find feasible assignments. The choice of the variable to branch on at a given node is a key factor in determining search efficiency \citep{achterberg2005branching, Glankwamdee2011LookaheadBF, schubert2017thesis, yang2019multivariablebranching}. We train a deep neural network policy to imitate choices made by an expert policy. The imitation target is a well-known heuristic called Full Strong Branching (FSB), which has been empirically shown to produce small search trees \citep{achterberg2005branching}. While it is often too computationally expensive for practical MIP solving, it can still be used to generate imitation learning data offline as a slow and expensive one-time computation. Once trained, the neural network is able to approximate the expert at test time at a fraction of the computational cost. A CPU-based implementation of FSB can be too expensive on large-scale MIPs even for offline data generation. We develop a variant of FSB using the alternating directions method of multipliers (ADMM) \citep{boyd2011distributed} that scales to large-scale MIPs by performing the required computation in a batch manner on GPU.

\begin{figure}[t]
    \includegraphics[trim={0cm 0cm 0cm 0cm}, clip, width=0.98\textwidth]{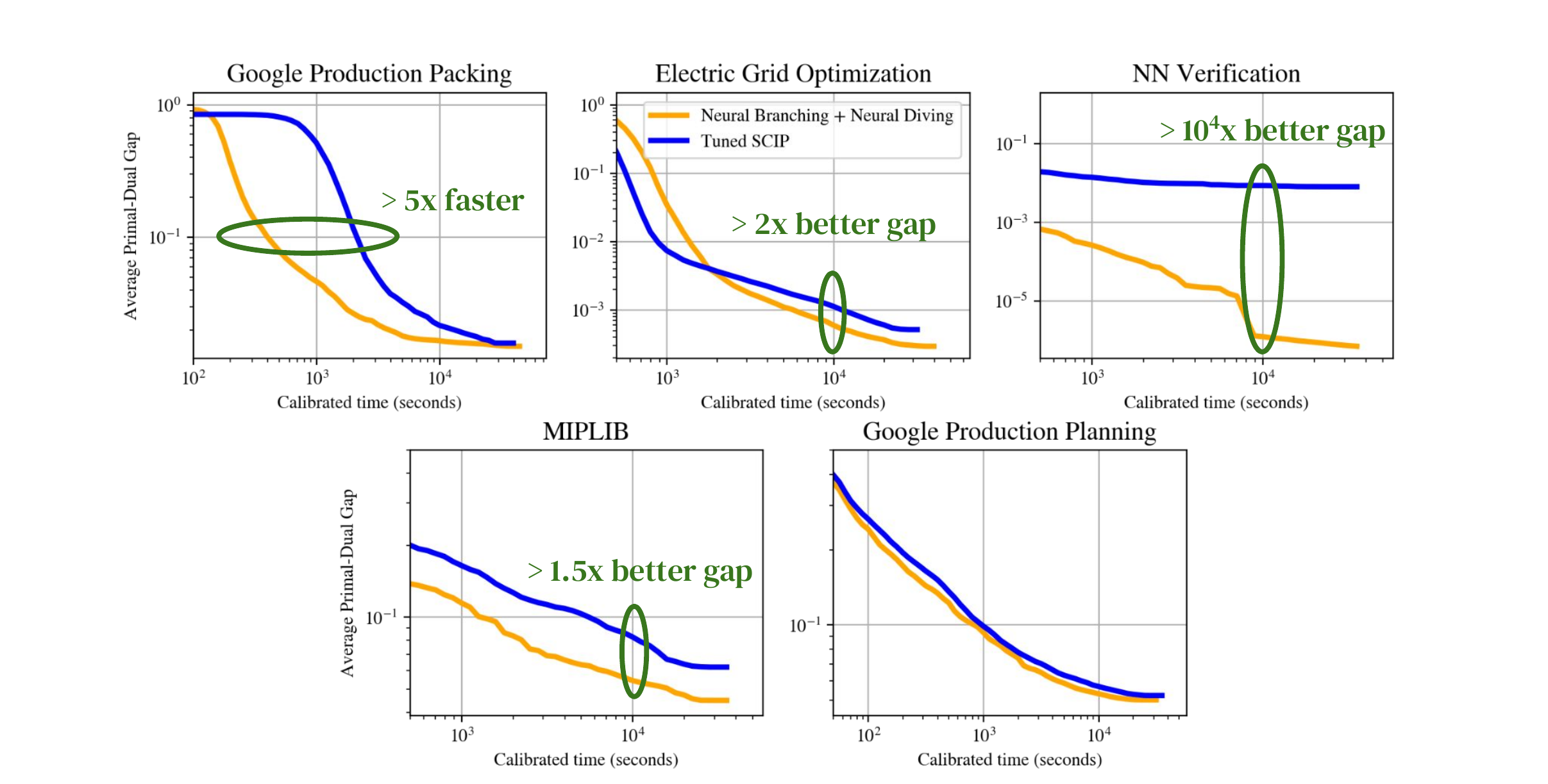}
    \caption{\small{Main result of the paper: our approach, Neural Branching + Neural Diving, matches or outperforms SCIP with respect to the primal-dual gap \citep{berthold2006thesis}, averaged on held-out instances, as a function of running time (referred to as calibrated time, see section \ref{sec:evaluation}) for the five datasets in our evaluation with the largest MIPs. Note the log scale on both axes.}}
    \label{fig:intro_result_summary}
\end{figure}

We evaluate our approach on a diverse set of datasets containing large-scale MIPs from real-world applications, including two from Google's production systems, as well as MIPLIB \citep{gleixner2019miplib}, which is a heterogeneous dataset and a standard benchmark. Most of the combined set of MIPs from all datasets have $10^3$--$10^6$ variables and constraints after presolve (see figure~\ref{fig:presolved_mip_sizes} in section~\ref{sec:datasets}), which is significantly larger than earlier works \citep{gasse2019exact,ding2019accelerating}. Once Neural Diving and Neural Branching models are trained on a given dataset, they are integrated into SCIP to form a ``Neural Solver'' specifically for that dataset. The baseline is SCIP with its emphasis parameters tuned by grid search on each dataset separately, which we refer to as Tuned SCIP. Comparing the Neural Solver and Tuned SCIP with respect to the \emph{primal-dual gap} \citep{berthold2006thesis} averaged over a held-out set of instances in figure~\ref{fig:intro_result_summary}, the Neural Solver provides either significantly better gap in the same running time, or the same gap in significantly less time, on four out of the five datasets in our evaluation with the largest MIPs, while matching Tuned SCIP's performance on the fifth. To the best of our knowledge, this is the first work to demonstrate such large improvements over SCIP on both large-scale real-world application datasets and MIPLIB using machine learning.

Tuned SCIP is the baseline we compare to since we use SCIP as the base solver for integrating learned heuristics. As a base solver SCIP provides a) extensive access to its internal state for integrating learned models, and b) permissive licensing that enables large-scale evaluation by running a large number of solver instances in parallel. We do not have access to commercial solvers with these features, which makes a fair comparison to them infeasible. We have done a partial comparison of Gurobi versus Neural Diving alone, with Gurobi as its sub-MIP solver, on two datasets (see section \ref{sec:supermip}). Comparing the \emph{primal gap} \citep{berthold2006thesis} averaged over a held-out set of instances, Neural Diving with parallel sub-MIP solving reaches 1\% average primal gap in 3$\times$ and 3.6$\times$ less time than Gurobi on the two datasets. We have also applied a modified version of Neural Diving to a subset of `open' instances in MIPLIB to find three new best known assignments, beating commercial solvers.

Several earlier works have focused on learning primal heuristics \citep{khalil2017primal, ding2019accelerating, Hendel2018alns, hottung2019neurallns, xavier2020mluc, song2020generallns, addanki2020nlns}. Unlike them, Neural Diving poses the problem of predicting variable assignments as a generative modelling problem, which provides a principled way to learn on all available feasible assignments and also to generate partial assignments at test time. Several works have also looked at learning a branching policy \citep{He2014LearningToSearch, Khalil2016LearningToBranch, Alvarez2017MLApproxStrongBranching, balcan18learningtobranch, gasse2019exact, Yang2020LearningGS, zarpellon2020parameterizing, gupta2020hybrid}. Many of these focus specifically on learning to imitate FSB as we do \citep{Khalil2016LearningToBranch, Alvarez2017MLApproxStrongBranching, gasse2019exact, gupta2020hybrid}. Unlike them, Neural Branching uses a more scalable approach to computing the target policy using GPUs, which allows it to generate more imitation data from larger instances in the same time limit than a CPU-based FSB implementation. We also go beyond earlier works that study learning individual heuristics in isolation by combining a learned primal heuristic and a learned branching policy in a solver to achieve significantly better performance on large-scale real-world application datasets and MIPLIB.

\subsection{Contributions}
\label{subsec:contributions}
\begin{enumerate}
    \item We propose Neural Diving (section \ref{sec:supermip}), a new learning-based approach to generating high-quality joint variable assignments for a MIP. On homogeneous datasets, Neural Diving achieves an average primal gap of 1\% on held out instances 3-10$\times$ faster than Tuned SCIP. On one dataset Tuned SCIP does not reach 1\% average primal gap within the time limit, while Neural Diving does. 
    \item We propose Neural Branching which learns a branching policy to use in the branch-and-bound algorithm by imitating a new scalable expert policy based on ADMM (section \ref{sec:deep_Branching}). On two of the datasets used in our evaluation for which FSB is slow due to instance sizes (e.g., with $> 10^5$ variables) or high per-iteration times, the ADMM expert generates 1.4$\times$ and 12$\times$ more training data in the same running time. The learned policy significantly outperforms SCIP's branching heuristic on four datasets with 2-20$\times$ better average \emph{dual gap} \citep{berthold2006thesis} on held out instances at large time limits, and has comparable performance on the other rest.
    \item We combine Neural Diving and Neural Branching (section \ref{sec:joint_eval}) to attain significantly better performance than SCIP with respect to the average primal-dual gap on four out of the five datasets with the largest MIPs, while matching its performance on the fifth one.
\end{enumerate}

In addition, we have open-sourced a dataset for the application of Neural Network Verification (sections \ref{sec:datasets}, \ref{subsec:appendix_datasets}), which we hope will help further research on new learning techniques for MIPs.

%% file: background.tex
\section{Integer programming background}
\label{sec:background}
In this section we present some basic integer programming concepts and techniques relevant for the paper.
A mixed integer linear program has the form
\begin{equation}
\label{eqn:mip1}
\begin{array}{ll}
\mathrm{minimize} & c^\top x\\
\text{subject to} & Ax \leq b\\
    & l \leq x \leq u \\
    & x_i \in \mathbb{Z}, \quad i \in \mathcal{I}
\end{array}
\end{equation}
over variables $x \in \reals^n$, where $A \in \mathbb{R}^{m \times n}$, $b \in \mathbb{R}^m$, $c \in \R^n$ , $l \in (\reals \cup \{-\infty\})^n$, and $u \in (\reals \cup \{\infty\})^n$ are given data, 
and $\mathcal{I} \subseteq \{1,\ldots,n\}$ refers to the index set of integer variables. We allow $l$ and $u$ to take infinite values to indicate that the associated variable has no lower or upper bound respectively, for those indices. A (complete) \emph{assignment} is any point $x \in \reals^n$. A partial assignment is when we have fixed some, but not all, of the variable values. A \emph{feasible} assignment is an assignment that satisfies all the constraints in problem \eqref{eqn:mip1}, and an \emph{optimal} assignment, or a \emph{solution}, is a feasible assignment that also minimizes the objective \citep[\S 1.1]{boyd2004convex}.

\subsection{Linear programming relaxation}
If we remove the integer constraints in problem \eqref{eqn:mip1} then it becomes a linear program (LP),
which is convex and can be solved efficiently \citep{boyd2004convex}. The optimal value of the
relaxed problem is guaranteed to be a lower bound for the original problem, since removing constraints
can only expand the feasible set. If the optimal solution to the LP satisfies the integral constraints
then it is guaranteed to be optimal for the original problem.
We shall refer to any lower bound (in our case found via linear programming) as a \emph{dual bound}.

\subsection{Branch-and-bound}
A common procedure to solve MIPs is to recursively build a search tree, with partial integer assignments at each node, and use the information gathered at each node to eventually converge on an optimal (or close to optimal) assignment \citep{Land1960BranchAndBound, lawler1966branch}. At every step we must choose a leaf node from which to `branch'. At this node we can solve the LP relaxation, where we constrain the ranges of the fixed variables at that node to their assigned values.
This gives us a valid lower bound on the true objective value of any further child nodes from this node. If this bound is larger than a known feasible assignment, then we can safely prune this part of the search tree since no optima for the original problem can exist in the subtree from that node. If we decide to expand this node, then we must choose a variable to branch on from the set of unfixed variables at that node. Once a variable is selected, we take a branching step, which adds two child nodes to the current node. One node has the domain of the selected variable constrained to be greater or equal to the ceiling of its LP relaxation value at the parent node. The other node has the selected variable's domain constrained to be less than or equal to the floor of its LP relaxation value. The tree is updated, and the procedure begins again. This algorithm is referred to as \emph{branch-and-bound}. Linear programming is a main workhorse of this procedure, both to derive the dual bounds at every node and to decide the branching variable for some of the more sophisticated branching heuristics. In theory the size of the search tree can be exponential in the input size of the problem, but in many cases the search trees can be small and it is an area of established and active research to come up with node selection and variable selection heuristics that keep the tree as small as possible.

\subsection{Primal heuristics}
A primal heuristic is a method that attempts to find a \emph{feasible}, but not necessarily optimal, variable assignment \citep{berthold2006thesis}. Any such feasible assignment provides a guaranteed upper bound on the optimal value of the MIP. Any such bound found at any point during a MIP solve is called a \emph{primal bound}.
Primal heuristics can be run independently of branch-and-bound, but they can also be run within a branch-and-bound tree and attempt to find a feasible assignment of the unfixed variables from a given node in the search tree.
Better primal heuristics that produce lower primal bounds allow the branch-and-bound procedure to prune more of the tree. Simple rounding, where the fractional variables are rounded to integer values (possibly randomly) is an example of a primal heuristic. Another case is \emph{diving} which attempts to find a feasible solution by exploring the search tree from a given node in a depth-first manner \citep{berthold2006thesis, eckstein2007pivot}.

\subsection{Primal-dual gap}
When running branch-and-bound we keep track of the \emph{global} primal bound
(the minimum objective value of any feasible assignment) and the \emph{global} dual bound (the minimum dual bound across all leaves of the branch-and-bound tree).
We can combine these to define a sub-optimality gap
\[
\text{gap} = \text{global primal bound} - \text{global dual bound}.
\]
The gap is always nonnegative by construction, and if it is zero then we have solved the problem, the feasible point that corresponds to the primal bound is optimal and the dual bound is a certificate of optimality. In practice we terminate branch-and-bound when the \emph{relative} gap (\ie, normalized in some way, see \S \ref{sec:evaluation}) is below some application-dependent quantity, and produce the best found primal solution as the approximately optimal solution.

%% file: rep_and_arch.tex
\begin{figure*}
    \centering
    \includegraphics[trim={0cm 3.9cm 0cm 2.0cm},clip,width=\textwidth]{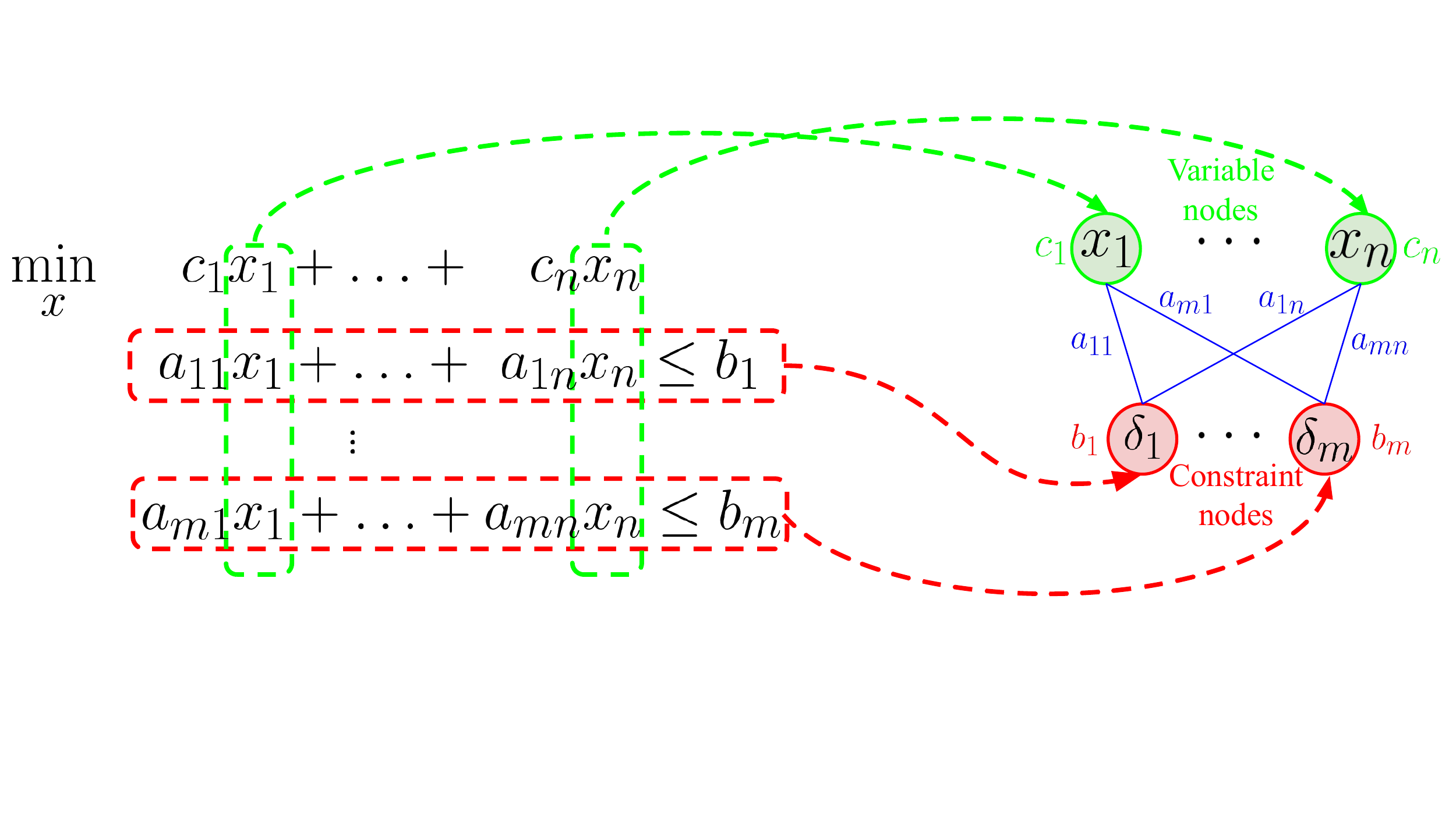}
    \caption{Bipartite graph representation of a MIP used as the input to a neural network. The set of $n$ variables $\{x_1,\ldots,x_n\}$ and the set of $m$ constraints $\{\delta_1,\ldots,\delta_m\}$ form the two sets of nodes of the bipartite graph. The coefficients are encoded as features of the nodes and edges.}
    \label{fig:mip_bipartite_graph}
\end{figure*}

\section{MIP Representation and Neural Network Architecture}
\label{sec:input_representation_and_architecture}
We describe how a MIP is represented as an input to a neural network, and the architecture we use to learn models for both Neural Diving and Neural Branching. The key deep learning architecture we use
is a form of \emph{graph neural network} (\cite{scarselli2009graph}, see survey by \cite{battaglia_survey_full_author_list}), specifically a \emph{graph convolutional network} (GCN) \citep{kipf2016semi}.

\subsection{Representing a MIP as an Input to a Neural Network}
\label{subsec:mip_bipartite_graph}
We use a bipartite graph representation of a MIP, as done in \citet{gasse2019exact}. Equation~\eqref{eqn:mip1} can be used to define a bipartite graph where one set of $n$ nodes in the graph correspond to the $n$ variables being optimized, and the other set of $m$ nodes correspond to the $m$ constraints; see figure~\ref{fig:mip_bipartite_graph}. The edge is present between a variable node and a constraint node if the corresponding variable appears in that constraint, and so the number of edges correspond to the number of non-zeros in the constraint matrix. The objective coefficients $\{c_1,\ldots,c_n\}$, the constraint bounds $\{b_1,\ldots,b_m\}$, and the non-zero coefficients of the constraint matrix are encoded as scalar ``features'' that annotate the corresponding variable nodes, constraint nodes, and edges, respectively. The variable type (continuous or integer) can also be encoded as a (categorical) feature of the variable nodes. This defines a lossless representation of the MIP that can be used as an input to a graph neural network. A similar representation is used for MIPs by \cite{ding2019accelerating}, and for Boolean Satisfiability instances in an earlier work by \cite{selsam2019NeuroSAT}. Both nodes and edges can be annotated by multi-dimensional feature vectors that encode additional information about the MIP that can be useful for learning (e.g., the solution of the LP relaxation as additional variable node features). We use the code provided by \citet{gasse2019exact} to compute the same set of features using SCIP.\footnotemark
\footnotetext{For the list of features, see \url{https://papers.nips.cc/paper/2019/file/d14c2267d848abeb81fd590f371d39bd-Supplemental.zip}}

\subsection{Neural Network Architecture}
\label{subsec:architecture}
We describe here the common aspects of the network architecture used by both Neural Diving and Neural Branching. Those aspects that differ between the two, e.g., the outputs and the loss functions, are given in their corresponding sections \ref{sec:supermip} and \ref{sec:deep_Branching}.

Given the bipartite graph representation of a MIP, we use a GCN to learn models for both Neural Diving and Neural Branching. Let the input to the GCN be a graph $G = (\mathcal{V}, \mathcal{E}, \mathcal{A})$ defined by the set of nodes $\mathcal{V}$, the set of edges $\mathcal{E}$, and the graph adjacency matrix $\mathcal{A}$. In the case of MIP bipartite graphs, $\mathcal{V}$ is the union of $n$ variable nodes and $m$ constraint nodes, of size $N := |\mathcal{V}| = n+m$. $\mathcal{A}$ is an $N\times N$ binary matrix with $\mathcal{A}_{ij}=1$ if nodes indexed by $i$ and $j$ are connected by an edge, 0 otherwise, and $\mathcal{A}_{ii}=1$ for all i. Each node has a $D$-dimensional feature vector, denoted by $u_i \in \mathbb{R}^{D}$ for the $i^\text{th}$ node. Let $U \in \mathbb{R}^{N \times D}$ be the matrix containing feature vectors of all nodes as rows, \ie, the $i^\text{th}$ row is $u_i$.
A single-layer GCN learns to compute an $H$-dimensional continuous vector representation for each node of the input graph, referred to as a \emph{node embedding}. Let $z_i \in \mathbb{R}^{H}$ be the node embedding computed by the GCN for the $i^\text{th}$ node, and $Z \in \mathbb{R}^{N \times H}$ be the matrix containing all node embeddings as rows. We define the function computing $Z$ as follows:
\begin{equation}
    Z = \mathcal{A}f_\theta(U),
\end{equation}
where $f_\theta : \reals^{D} \rightarrow \reals^{H}$ is a Multi-Layer Perceptron (MLP) \citep{Goodfellow-et-al-2016} with learnable parameters $\theta \in \Theta$. (Here we have generalized $f_\theta$ from a linear mapping followed by a fixed nonlinearity in the standard GCN by \cite{kipf2016semi} to an MLP.) We overload the notation to allow $f_\theta$ to operate on $N$ nodes simultaneously, \ie,  $f_\theta(U)$ denotes applying the MLP to each row of $U$ to compute the corresponding row of its output matrix of size $N \times H$. Multiplying by $\mathcal{A}$ combines the MLP outputs of the $i^\text{th}$ node's neighbors to compute its node embedding.
The above definition can be generalized to $L$ layers as follows:
\begin{align}
    Z^{(0)} &= U\\
    Z^{(l+1)} &= \mathcal{A}f_{\theta(l)}(Z^{(l)}), \quad l=0, \ldots, L-1,
\end{align}
where $Z^{(l)}$ and $f_{\theta(l)}()$ denote the node embeddings and the MLP, respectively, for the $l^\text{th}$ layer. The $L^\text{th}$ layer's node embeddings can be used as input to another MLP that compute the outputs for the final prediction task, as shown in sections~\ref{sec:supermip} and~\ref{sec:deep_Branching}. 
For each variable $x_d$ we further denote the corresponding node embedding from the last layer as $v_d$.

Two key properties of the bipartite graph representation of the MIP and the GCN architecture are: 1) the network output is invariant to permutations of variables and constraints, and 2) the network can be applied to MIPs of different sizes using the same set of parameters. Both of these are important because there may not be any canonical ordering for variables and constraints, and different instances within the same application can have different number of variables and constraints.

\subsection{Improvements to the Architecture}
\label{subsec:arch_improvements}
We describe changes to the above architecture, as well as from previous work, which gave performance improvements.
\begin{enumerate}
    \item We modify the MIP bipartite graph's adjacency matrix $\mathcal{A}$ to contain coefficients from the MIP's constraint matrix $A$, instead of binary values indicating the presence of edges. Specifically, for the $i^\text{th}$ variable and $j^\text{th}$ constraint, their two corresponding entries in $\mathcal{A}$ is set to $a_{ji}$ where $a_{ji}$ is the coefficient at row $j$ and column $i$ of $A$. This results in edges weighted by the entries of $A$.
    \item We extend the node embeddings for layer $l+1$ by concatenating the node embeddings from layer $l$. Specifically, we now define the embedding for layer $l+1$ to be $\widetilde{Z}^{(l+1)} = (Z^{(l+1)}, \widetilde{Z}^{(l)})$, \ie, the concatenation of the matrices row-wise, with $\widetilde{Z}^{(0)} = Z^{0}$. This is a form of skip connection \citep{skip_connections2016} and also similar to the jumping knowledge networks architecture \citep{xu2018representation}.
    \item We apply layer norm \citep{layer_normalization2016} at the output of each layer, so that $Z^{(l+1)} = \text{LayerNorm}\left(\mathcal{A}f_{\theta(l)}(Z^{(l)})\right)$.
\end{enumerate}
We have explored alternative architectures which use embeddings for both nodes and edges with separate MLPs to compute them. When using such networks with a high-dimensional edge embedding at every layer, their memory usage can become much higher than that of GCNs, which only need the adjacency matrix, and may not fit the GPU memory unless the number of layers is reduced at the cost of accuracy. GCNs are a better fit for our goal of scaling to large-scale MIPs.

%% file: datasets.tex
\section{Datasets}
\label{sec:datasets}

\begin{table}
\caption{MIP datasets used for evaluation. See additional details in section \ref{subsec:appendix_datasets}}.
\begin{adjustbox}{center}
    \begin{tabular}{|p{3.2cm}|p{12cm}|}
    \hline
        \textbf{Name} & \textbf{Description} \\
        
        \hline

        CORLAT & A small-scale dataset related to wildlife management \citep{conrad2007corlat1, gomes2008corlat2}. We include it as it is public, but our main focus is on the larger-scale ones.\\
 
        \hline
 
        NN Verification & Verifying whether a neural network is robust
        to input perturbations can be posed as a MIP \citep{cheng2017verification, tjeng2018evaluating}.
        Each input on which to verify the network gives rise to a different MIP.
        In this dataset, a convolutional neural network is verified on each
        image in the MNIST dataset, giving rise to a corresponding dataset of MIPs.\\
  
        \hline
  
        Google Production Packing & A packing optimization problem solved in a Google production system.\\
        
        \hline
        
        Google Production Planning & A planning optimization problem solved in a Google production system.\\
        
        \hline
        
        Electric Grid Optimization & Electric grid operators optimize the choice of power generators to use at different time intervals during a day to meet electricity demand by solving a MIP. This dataset
        is constructed for one of the largest grid operators in
        the US, PJM, using publicly available data about generators and demand, and the MIP formulation in \citep{knueven2018mip_uc}.\\
        
        \hline
        
        MIPLIB & Heterogeneous dataset containing `hard' instances of MIPs
        across many different application areas that is used as a long-standing standard benchmark for MIP solvers \citep{gleixner2019miplib}. We use instances from both the 2010 and 2017 versions of MIPLIB.\\ %
    \hline
    \end{tabular}
\end{adjustbox}
\label{tab:datasets}
\end{table}

\begin{figure*}[t]
    \centering
    \begin{tabular}{cc}
        \includegraphics[width=0.48\textwidth]{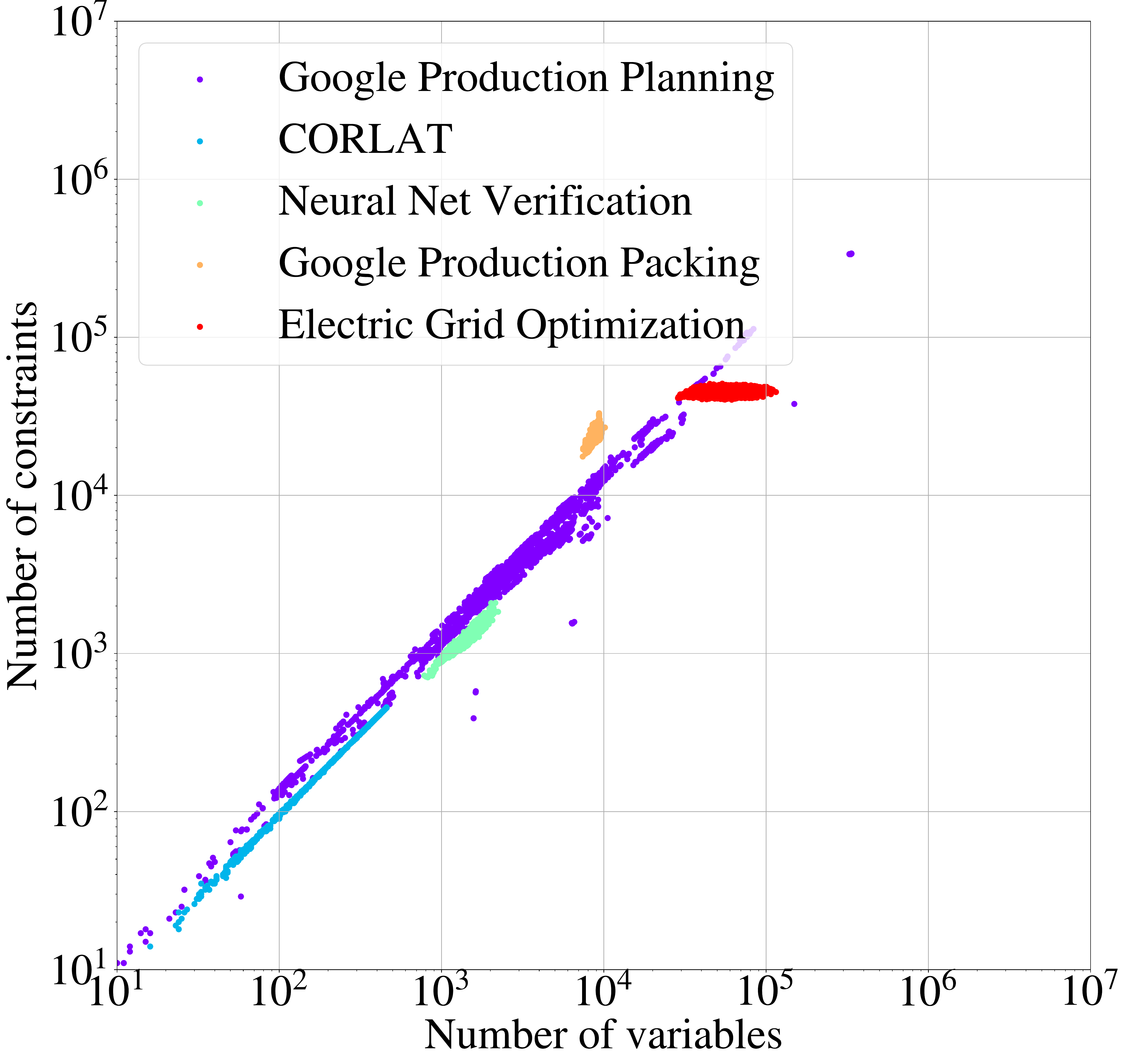}
        \hspace{0.1in}
        \includegraphics[width=0.48\textwidth]{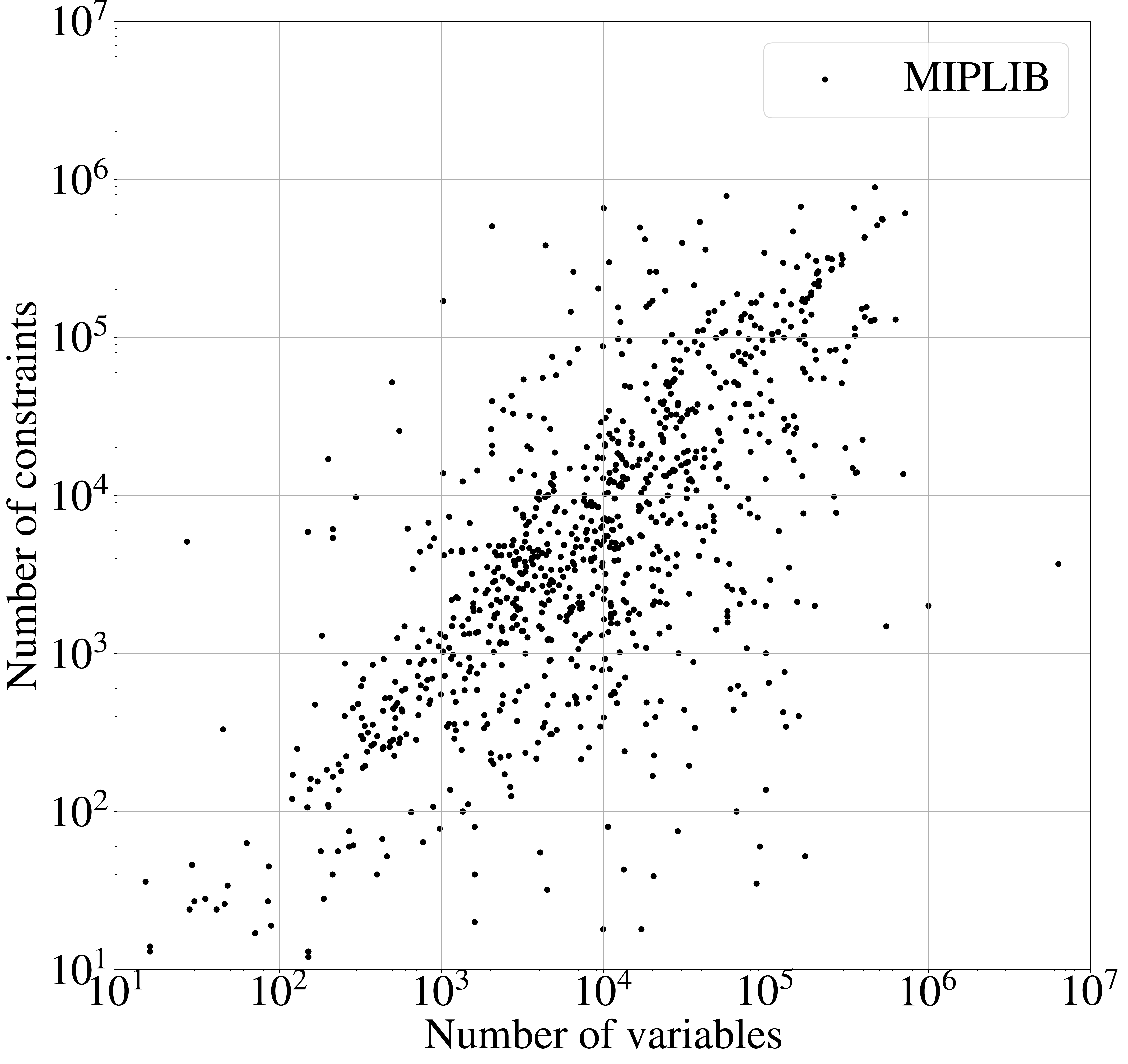}
    \end{tabular}
    \caption{
        {\small Number of variables and constraints after presolve using SCIP 7.0.1 for the datasets used in our evaluation. MIPLIB is shown separately (right) to reduce clutter.}
    }
    \label{fig:presolved_mip_sizes}
\end{figure*}

We summarize the details of our datasets in table~\ref{tab:datasets}. All datasets except MIPLIB are application-specific, i.e., they contain instances from only a single application. We include MIPLIB in our evaluation even though it contains instances from many applications. This is not the setting in which we expect learning will provide the most benefit. Nevertheless, it is a well-established benchmark for evaluating solvers, and the results on it can provide insight on the performance of learning approaches in an unfavorable setting. The MIP sizes for the various datasets after presolving are shown in figure~\ref{fig:presolved_mip_sizes}. 

For all datasets except MIPLIB, we define the training, validation, and test sets by randomly splitting the instances into disjoint subsets with 70\%, 15\%, and 15\% of the instances, respectively. For MIPLIB, we use instances from the MIPLIB 2017 Benchmark Set as the test set, since that is the set on which solvers are evaluated. The MIPLIB 2017 Collection Set and the MIPLIB 2010 set are used as the training and validation sets, respectively, after removing any overlapping instances with the MIPLIB 2017 Benchmark Set. For each dataset, the training set is used to learn models for that dataset, the validation set is used to tune the learning hyperparameters and SCIP's metaparameters, and the test set is used to report evaluation results.  Additional details are given in the Appendix, see section~\ref{subsec:appendix_datasets}.

%% file: evaluation.tex
\section{Evaluation}
\label{sec:evaluation}

We shall evaluate Neural Diving and Neural Branching individually first, and then jointly in the sequel. In all cases we evaluate on a test set of MIPs disjoint from the training set to measure generalization to unseen instances. We present results using \emph{gap plots} and \emph{survival plots}:

\paragraph{\textbf{Gap plots:}} Neural Diving is evaluated with respect to the \emph{primal gap} $\gamma_p(t)$ \citep{berthold2013primal} as a function of solving time $t$. If no feasible assignment is known then
the gap is defined to be $1$, otherwise denote by $p(t)$ the primal bound at time $t$ and denote by $p^\star$ the best known primal bound (possibly precomputed), the gap is defined as
\begin{equation}
\label{eqn:primal_gap}
    \gamma_{p}(t) =
\begin{cases}
1 & p(t) \cdot p^\star < 0\\
\frac{p(t)-p^\star}{\max{\{|p(t)|,|p^\star|, \epsilon\}}} & \rm{otherwise.} \\
\end{cases}
\end{equation}
We use $\epsilon=10^{-12}$ to avoid division by 0. To  evaluate Neural Branching we use a dual gap, which is defined analogously
\begin{equation}
\label{eqn:dual_gap}
    \gamma_{d}(t) =
\begin{cases}
1 & d(t) \cdot p^\star < 0\\
\frac{p^\star - d(t)}{\max{\{|d(t)|,|p^\star|, \epsilon\}}} & \rm{otherwise,} \\
\end{cases}
\end{equation}
where $d(t)$ is the global dual bound at time $t$. When evaluating Neural Diving and Neural Branching jointly, we use the \emph{primal-dual gap}
\begin{equation}
\label{eqn:primal_dual_gap}
    \gamma_{pd}(t) =
\begin{cases}
1 & d(t) \cdot p(t) < 0\\
\frac{p(t) - d(t)}{\max{\{|d(t)|,|p(t)|, \epsilon\}}} & \rm{otherwise.} \\
\end{cases}
\end{equation}
In each case, we plot the average gap for test MIPs as a function of running time. We pre-compute $p^\star$ for a MIP by running SCIP on it with default parameters and a time limit of 24 hours. Any suboptimality in $p^\star$ will affect all solvers being evaluated and relative comparisons using the gaps are therefore still valid. If during the solve we find a better primal bound then we replace $p^\star$ with the new value and recompute gaps at all previous times.

\paragraph{\textbf{Survival plots:}} A survival plot shows the fraction of test set MIPs solved (to the target primal-dual gap) as a function of running time. Applications may specify a target optimality gap higher than 0 as a termination condition. In such a case improving the gap below the target does not improve the solve performance. The survival plot reflects the effect of the target gap, while the gap plots do not.

\paragraph{\textbf{Calibrated time:}} The total evaluation workload across all datasets and comparisons requires more than 160,000 MIP solves and nearly a million CPU and GPU hours. To meet the compute requirements, we use a shared, heterogeneous compute cluster. Accurate running time measurement on such a cluster is difficult because the tasks may be scheduled on machines with different hardware, and interference from other unrelated tasks on the same machine increases the variance of solve times. To improve accuracy, for each solve task, we periodically solve a small \emph{calibration MIP} on a different thread from the solve task on the same machine. We use an estimate of the number of calibration MIP solves during the solve task on the same machine to measure time, which is significantly less sensitive to hardware heterogeneity and interference. This quantity is then converted into a time value using the calibration MIP's solve time on a reference machine. Section \ref{appendix:calibrated_time_results} gives the details. Results for four instances from MIPLIB show a $1.5\times$ to $30\times$ reduction in the coefficient of variation of time measurements compared to measuring wall clock time.

\paragraph{\textbf{Tuned SCIP:}} The main baseline we compare against is SCIP 7.0.1 with its parameters tuned for each test dataset. SCIP has emphasis ``meta-parameters'' for presolving, primal heuristics, and cuts, each of which has four possible settings (default, off, aggressive, fast). For each dataset, we use exhaustive search over the $4^3 = 64$ combinations to find the setting that produces the best average primal-dual gap on a subset of 200 validation MIPs in a 3 hour time limit. We call this baseline \emph{Tuned SCIP}.

\paragraph{\textbf{Performance variability with respect to random seed:}} To account for the performance variability of MIP solvers \citep{lodi2014variability} with respect to changes to a MIP that leave the problem unchanged, such as permuting rows and columns of the constraint matrix, we vary the random seed parameter used by SCIP. Specifically, we set the SCIP parameter \emph{randomization/permutevars} to \emph{True}, and assign the parameters \emph{randomization/permutationseed} and \emph{randomization/randomseedshift} both to a give seed value. The seed is set to the values $\{1,2,3,4,5\}$ for each instance. The evaluation results reported for Neural Diving, Neural Branching, and their combination are computed by aggregating over all the instance-seed pairs produced using the test set instances.

%% file: supermip.tex
\section{Neural Diving}
\label{sec:supermip}

In this section we describe our approach to learning a diving-style primal heuristic that produces high quality assignments to MIPs from a given instance distribution. The idea is to train a \emph{generative model} over assignments to a MIP's integer variables from which partial assignments can be sampled. We use SCIP to obtain high-quality assignments (not necessarily optimal) as the target labels for the training set of MIPs. Once trained on this data, the model predicts values for integer variables on unseen instances from the same problem distribution. The uncertainty represented in the model predictions is used to define partial assignments to the original MIP that fix a large fraction of the integer variables. These substantially smaller sub-MIPs can then be solved quickly using SCIP, yielding high quality feasible assignments.

\emph{Diving} refers to a set of primal heuristics that explore the branch-and-bound tree in a depth-first manner by sequentially fixing integer variables until a leaf node is reached or the assignment is deemed infeasible \citep{berthold2006thesis, eckstein2007pivot}. There are a few major differences between regular diving and what we describe here. First, diving can be started from any node, but in this work we focus on diving only from the root node, though in principle it could be performed from other nodes. Secondly, diving typically descends all the way to a leaf node, but in this case we descend only partially and then use a MIP solver to solve the remaining sub-MIP. For this reason one could reasonably call our method a hybrid of diving and neighborhood search \citep{mladenovic1997variable, shaw1998using, hansen2010variable}, but for brevity we refer to it simply as Neural Diving. Finally, diving usually proceeds in an iterative manner where decisions are made sequentially and the linear program is re-solved after each decision. The vanilla version of Neural Diving produces assignments that define the sub-MIP entirely in parallel, and only re-solves the linear program once. In section~\ref{appendix:autoregressive} we describe an extension where the decisions are made sequentially. Our approach is similar to Relaxation Enforced Neighborhood Search (RENS) \citep{berthold2007rens} in that we fix a subset of variables and solve the resulting sub-MIP, but in our case the variables are assigned values predicted by a learned model rather than based on the linear program solution, and we use multiple partial assignments instead of only one.

\subsection{Solution Prediction as Conditional Generative Modelling}

Consider an \emph{integer program} (\ie, all variables are integers) with parameters $M=(A,b,c)$ (see equation \ref{eqn:mip1}) and a nonempty feasible set over a set of integer variables $x$. Assuming minimization, we define an \emph{energy function} over $x$ using the objective function:
\begin{equation}
    E(x; M) = \begin{cases}
        c^Tx & \text{if $x$ is feasible,}\\
        \infty & \text{otherwise},\\
    \end{cases}
\end{equation}
which then defines the conditional distribution 
\begin{equation}
    p(x|M) = \frac{\exp(-E(x;M))}{Z(M)}
\end{equation}
where $Z(M)$ is the \emph{partition function} that normalizes the distribution to sum to 1:
\begin{equation}
    Z(M) = \sum_{x'}\exp(-E(x';M)).
\end{equation}
Feasible assignments with better (\ie, lower) objective values have higher probability. Infeasible assignments have zero probability. Note that multiplying $c$ by a constant $\beta$ changes the distribution, even though the feasible set and the optimal assignment(s) remain the same. ($\beta$ can be interpreted as the inverse temperature parameter of the distribution.) The distribution can be made invariant to such a re-scaling of the objective by normalizing the energy as $E(x;M)/|E(x^*;M)|$ where $x^*$ is an optimal assignment. We have not used such normalization in this work.

When a MIP contains both integer variables $x_I$ and continuous variables $x_C$, a given assignment for $x_I$ defines a linear program (LP) on $x_C$. The energy function $E(x;M)$ can then be computed by assigning $x_C$ to the optimal LP solution ${x^*}_C$ (if feasible) and then setting $E(x;M)$ to the objective value of the resulting complete assignment. %

\subsubsection{Learning}
\label{subsubsec:supermip_learning}
The learning task is to approximate $p(x|M)$ using a generative model $p_{\theta}(x|M)$ parameterized by $\theta$. The training dataset is $\mathcal{D}_{\text{train}} = \{ (X_i, M_i) \}_{i=1}^N$, where $\{M_i \sim p(M)\}_{i=1}^N$ are $N$ independent and identically distributed (IID) samples from the application-specific MIP distribution $p(M)$ and $X_i = \{ x^{i,j} \}_{j=1}^{N_i}$ is a set of unique $N_i$ assignments for the instance $M_i$. $X_i$ is obtained by running SCIP on $M_i$ and collecting feasible assignments it finds during the solve. While potentially costly, this data collection step needs to be done only once (per application), and outside the training loop. We learn on \emph{all} assignments (explained below), not just the best ones, and do not require that any of the assignments be optimal.

The model parameters $\theta$ are learned by minimizing the following weighted loss function with respect to $\theta$:
\begin{equation}\label{eq:NNS_maxlikelihood}
    L(\theta) = -\sum_{i=1}^N \sum_{j=1}^{N_i} w_{ij} \log p_{\theta}(x^{i,j} | M_i),
\end{equation}
where the weights $w_{ij}$ are used to reduce any bias in sampling the training examples $X_i = \{ x^{i,j} \}_{j=1}^{N_i}$ for the instance $M_i$. Consider the case where $X_i$ contains all possible feasible assignments of the integer variables for $M_i$ without any duplicates (\ie, $x^{i,j} \neq x^{i,k},~\forall j \neq k$). Then we can use the weights
\begin{equation}
\label{eq:importance_weights}
    w_{ij} = \frac{\exp(-c_i^T x^{i,j})}{\sum_{k=1}^{N_i} \exp(-c_i^T x^{i,k})}.
\end{equation}
in equation \ref{eq:NNS_maxlikelihood} to learn the target distribution $p(x|M_i)$ using $X_i$. Enumerating all possible assignments to construct $X_i$ is not practical for MIPs of realistic size. Instead we apply an off-the-shelf solver, SCIP in our case, to $M_i$. We define $X_i$ to include \emph{all} feasible assignments found throughout the solve, with duplicates removed. Using the weights given by equation \ref{eq:importance_weights} no longer corresponds to learning the exact target distribution $p(x|M_i)$ as the partition function is approximate. Although approximate, learning still succeeds in producing models with strong empirical performance on the datasets used in our evaluation. Intuitively we expect the approximation to produce good results because it assigns higher weights to better assignments.

\begin{figure*}
    \centering
    \includegraphics[width=\textwidth,trim={3.5cm 3.5cm 3.5cm 3cm},clip]{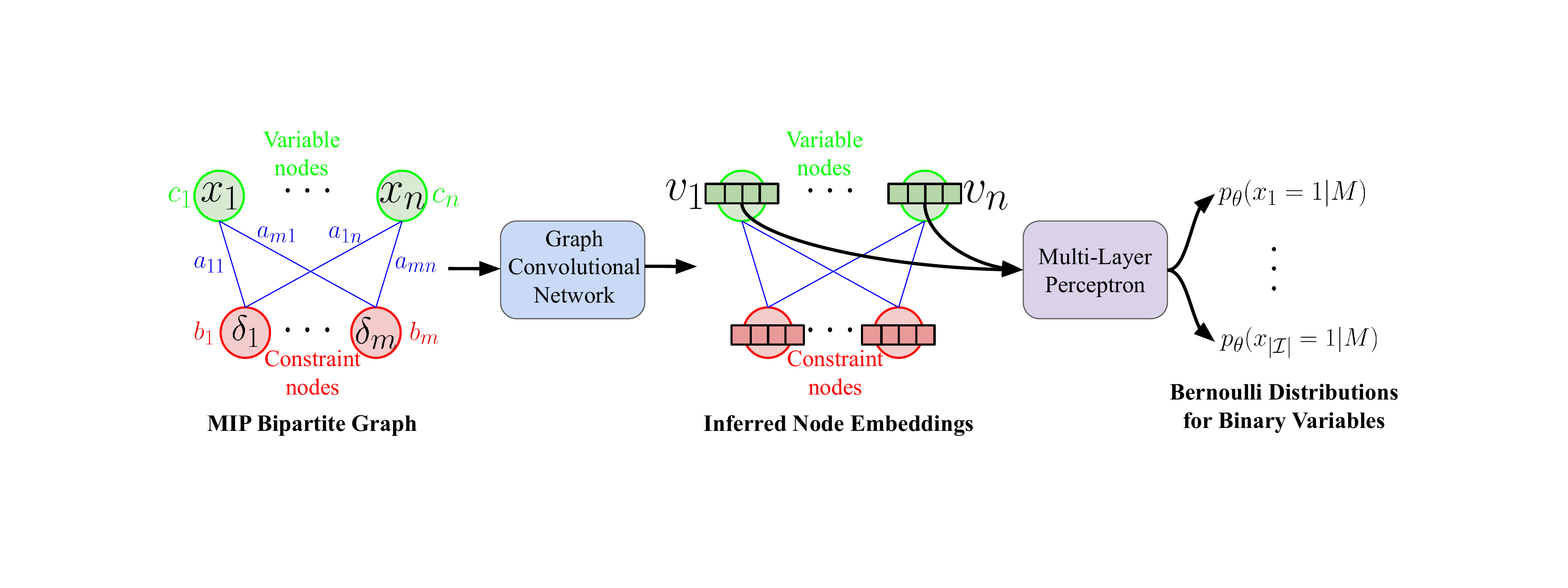}
    \caption{Conditionally independent model (Section \ref{subsubsec:ci_model}) infers node embeddings $\{v_1,\dots,v_n\}$ for the variable nodes of the bipartite graph representation of an input MIP $M$ using a graph convolutional network and applies a multi-layer perceptron to the node embeddings $\{v_d\}_{d \in \mathcal{I}}$ of the set of binary variables $\mathcal{I}$ in $M$ to compute their Bernoulli distributions. (See Section \ref{sec:NNS_integers} for the case of general integers.) These distributions are independent of each other conditioned on the input MIP.}
    \label{fig:conditionally_independent_predictor}
\end{figure*}

\subsubsection{Conditionally-independent Model}
\label{subsubsec:ci_model}
There are several choices for parameterizing the generative model (see, e.g., \cite{bishop2006prml, bengio00ar, kingma14vae}). Let $\mathcal{I}$ be the set of dimensions of $x$ corresponding to the integer variables. Let $x_d$ denote the $d^{\text{th}}$ dimension of $x$. We use a model of the form:
\begin{equation}
    p_{\theta}(x | M) = \prod_{d \in \mathcal{I}} p_{\theta}(x_d | M),
\end{equation}
which predicts a distribution for $x_d$ that is independent of the other dimensions conditioned on $M$.

For simplicity, we first assume that each $x_d$ is binary with possible values $\{0,1\}$ (the general integer case is considered in section~\ref{sec:NNS_integers}). We use the Bernoulli distribution for each such variable. The success probability $\mu_d$ for the Bernoulli distribution $p_{\theta}(x_d | M)$ is computed as 
\begin{align}
    t_d &= \text{MLP}(v_d; \theta),\\
    \mu_{d} &= p_{\theta}(x_d = 1 | M) = \frac{1}{1+\exp(-t_d)},
\end{align}
where $v_d$ is the embedding computed by a graph convolutional network for the MIP bipartite graph node corresponding to $x_d$, as described in section \ref{subsec:architecture}. Note that the same MLP is used for all variables, see figure~\ref{fig:conditionally_independent_predictor}.

While the conditionally-independent model cannot accurately model a multimodal distribution over assignments, empirically, it still shows strong performance in our experiments. We also investigated more sophisticated models, such as autoregressive models (\cite{bengio00ar}, see Section~\ref{appendix:autoregressive}), which provide modest improvements but at a much higher inference cost.

\subsubsection{Handling General Integers}
\label{sec:NNS_integers}

The model can be generalized to non-binary integer variables with some modifications. The main challenges are: a) the cardinality of an integer variable can vary significantly across instances, and b) it can potentially be very large (e.g., $10^7$). This makes simple approaches such as assuming a fixed maximum cardinality highly inefficient and difficult to learn with. 

We address these challenges by reframing the prediction task for general integer variables as a sequence of binary prediction tasks, based on the binary representation of the target integer value. For an integer variable $z$ that can be assigned values from a finite set with cardinality $\mathrm{card}(z)$, any target value can be represented as a sequence of $\lceil\log_2(\mathrm{card}(z))\rceil$ bits. We train our model to predict these bits in sequence, from most significant to least significant bit. Since the maximum cardinality of variables in a test instance becomes known only during inference, and is unknown during training, we introduce a hyperparameter $n_b$ that controls the maximum number of bits predicted for each variable along this bit sequence. We can then train our model to predict the $n_b$ most significant bits of the value for $z$, given the upper and lower bounds for $z$. During inference then,
if $\lceil\log_2(\mathrm{card}(z))\rceil$ exceeds $n_b$, we can still predict the $n_b$ most significant bits of the value for $z$, and use these predictions to tighten the bounds of $z$ by a factor of at least $2^{n_b-1}$. Otherwise, if $\lceil\log_2(\mathrm{card}(z))\rceil\leq n_b$, then we can predict $z$ exactly to a single value. Note that this approach of bit-wise prediction of values can also be seen as predicting branching decisions on variables: Predicting one bit means picking the left or right branch of a binary tree, and the possible range of the integer variable is successively tightened after each prediction.

\subsection{Combining Model Predictions with a Classical Solver}

Once a model is learned, we can then sample variable assignments $x \sim p_\theta(x | M)$ for a new problem $M$. While we could try to use these samples directly, they need not be feasible or provide the best objective value. Instead, it is more effective to only fix a (large) subset of the variables to their sampled values, and delegate the solution search for the remaining open variables to a classical solver, in our case, SCIP.

We use the SelectiveNet approach \citep{selectivenet} to train an additional binary classifier that decides which variables to predict a value for and which to refrain from predicting, and optimizes for ``coverage'' among variables, defined as the ratio of the number of variables predicted vs not predicted. Specifically, we introduce an additional output $y_d\in\{0,1\}$ for each variable $x_d$ in the input $M$ that determines whether $x_d$ should be assigned ($y_d = 1$) or not. For the conditionally-independent model, we can then train our model by minimizing the following loss function:

\begin{equation}
    l_{\rm selective}(\theta, x, M) = \frac{- \sum_{d \in \mathcal{I}} \log p_\theta(x_d|M) \cdot y_d }{\sum_{d \in \mathcal{I}} y_d} + \lambda \Psi(C - \frac{1}{|\mathcal{I}|}\sum_{d \in \mathcal{I}} y_d),
\end{equation}

\begin{equation}
    L_{\rm selective}(\theta) = \sum_{i,j} w_{ij} \cdot l_{\rm selective}(\theta, x^{i,j}, M_i).
\end{equation}
Here, $C$ is a coverage threshold representing the desired relative frequency of assigned variables, $\Psi$ is a quadratic penalty term, and $\lambda$ is a hyperparameter controlling the relative importance of achieving a coverage close to the set threshold.
We train multiple models simultaneously with different coverage thresholds (typically, values from 0.1 to 0.95, and are tuned on the validation set).

By assigning, or tightening the bounds of a large fraction of the variables we significantly reduce the problem size, and warm-start SCIP to find high quality solutions in much shorter time. The idea here is that by sampling a solution from $p_\theta(x | M)$, we move to a promising region of the solution space; and by removing again some of the assignments we expand the neighbourhood in which we then search for the optimal solution. 

This approach also offers practical computational advantages: both the sampling of predictions, and the solution search afterwards are fully parallelizable. We can repeatedly and independently extract many different samples from our model, and each partial assignment of the sample can be independently solved. For a single MIP instance, we can generate many such partial assignments, each of which we can independently solve with SCIP. The feasible solution with the best objective value across all sub-MIPs is then reported as the final solution to the original input MIP. Our approach can hence easily leverage the power of distributed computing even when solving for a single instance, which is not possible with default SCIP (without significant effort \citep{shinano2011parascip}). When comparing results, we will show results for both the parallel setting (where we make full use of the parallelism advantage of our approach) and the sequential setting (where we run our approach in a sequential manner, controlling for the total amount of computational resources).

\begin{figure}[ht!]
    \begin{tabular}{cc}
        \includegraphics[width=0.50\textwidth]{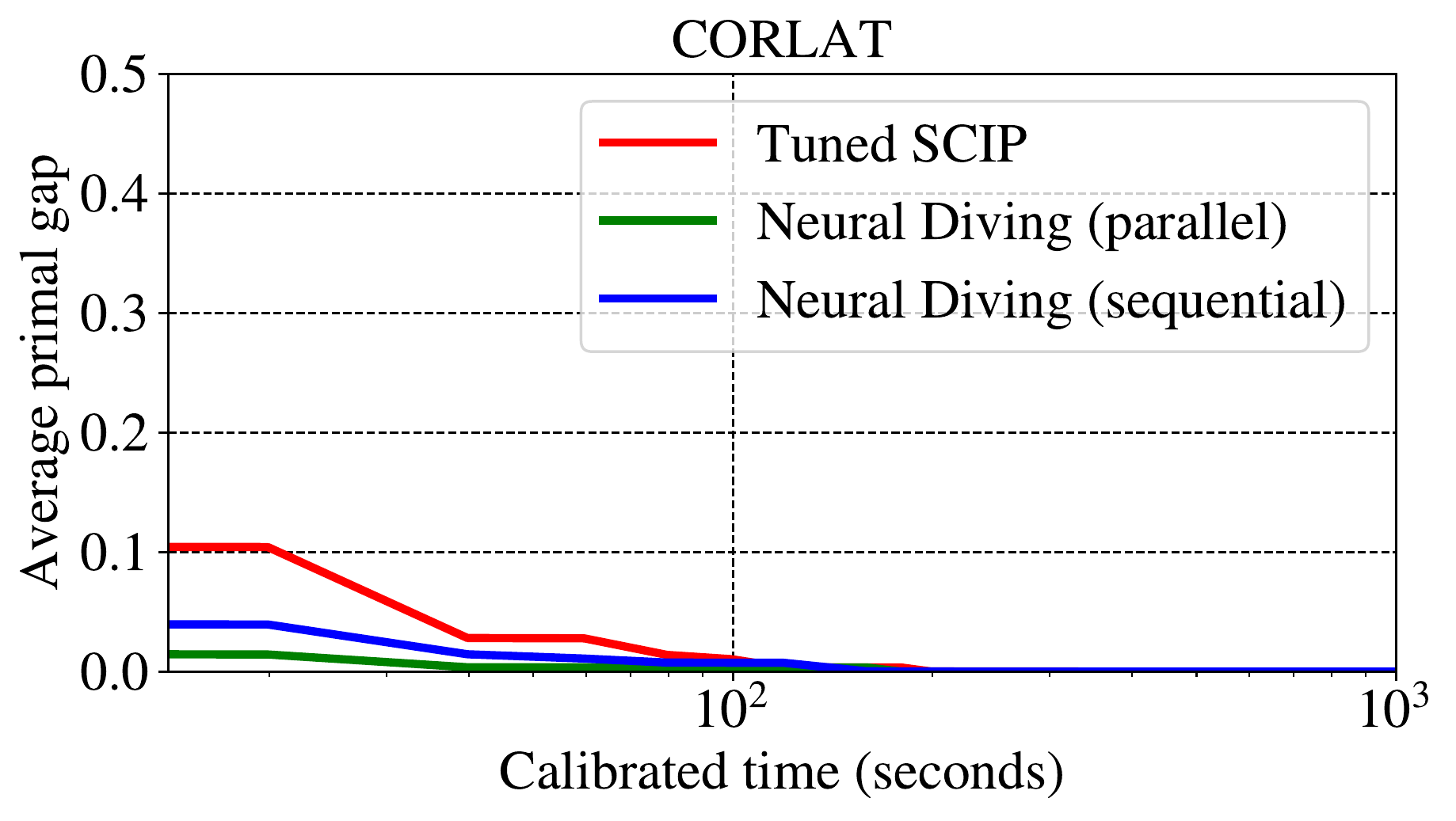} &
        \includegraphics[width=0.50\textwidth]{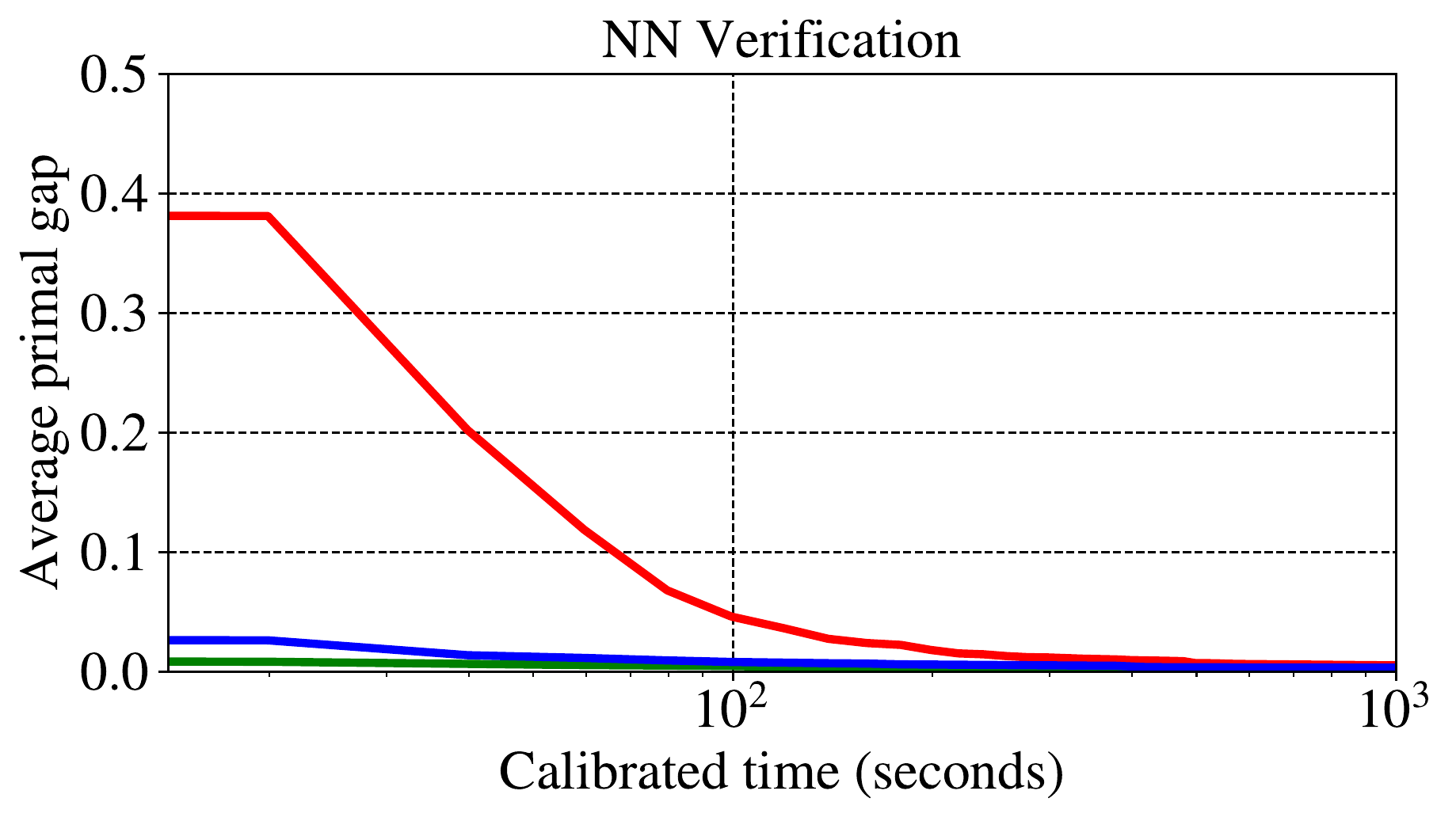} \\
        \includegraphics[width=0.50\textwidth]{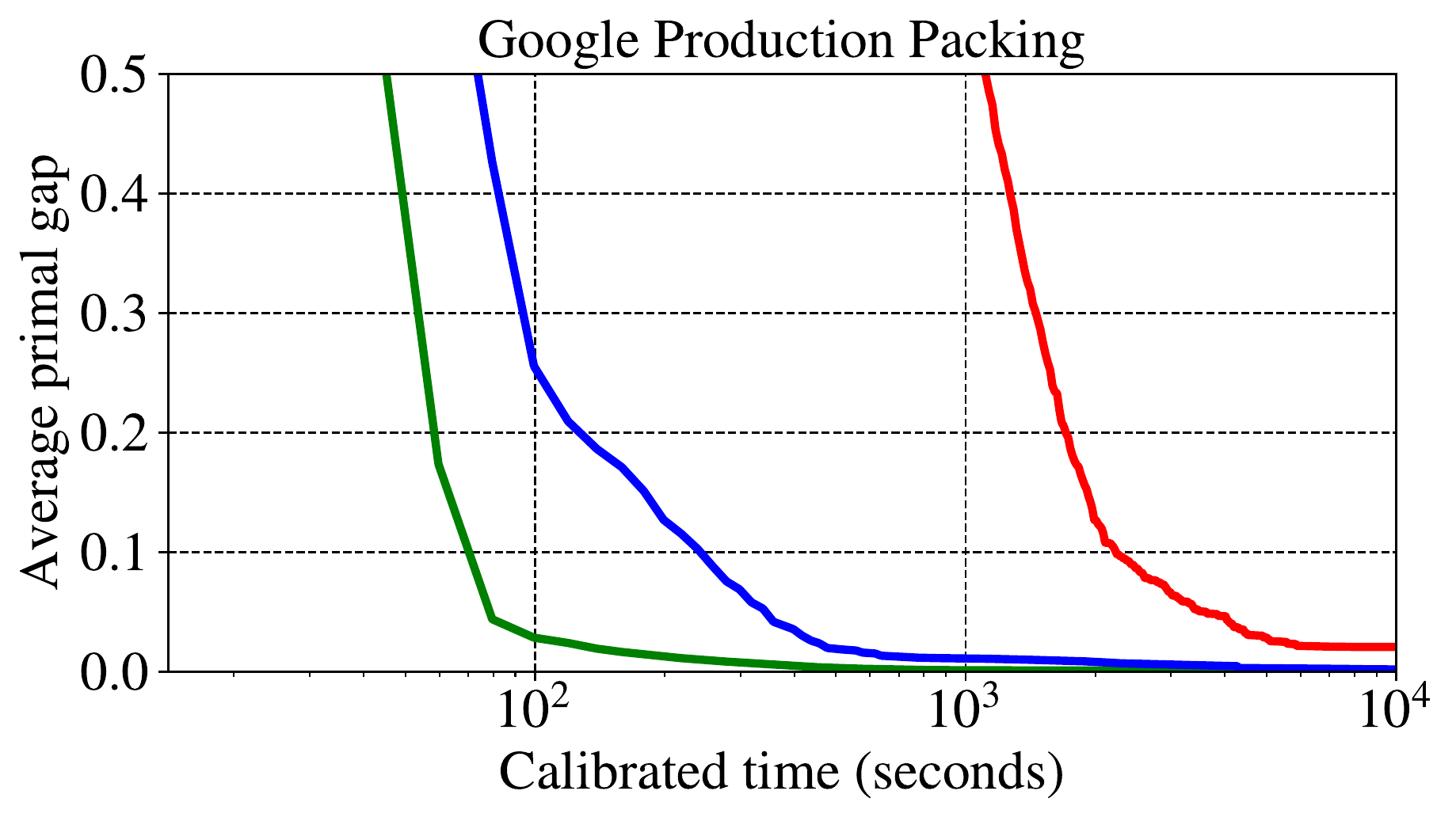} &
        \includegraphics[width=0.50\textwidth]{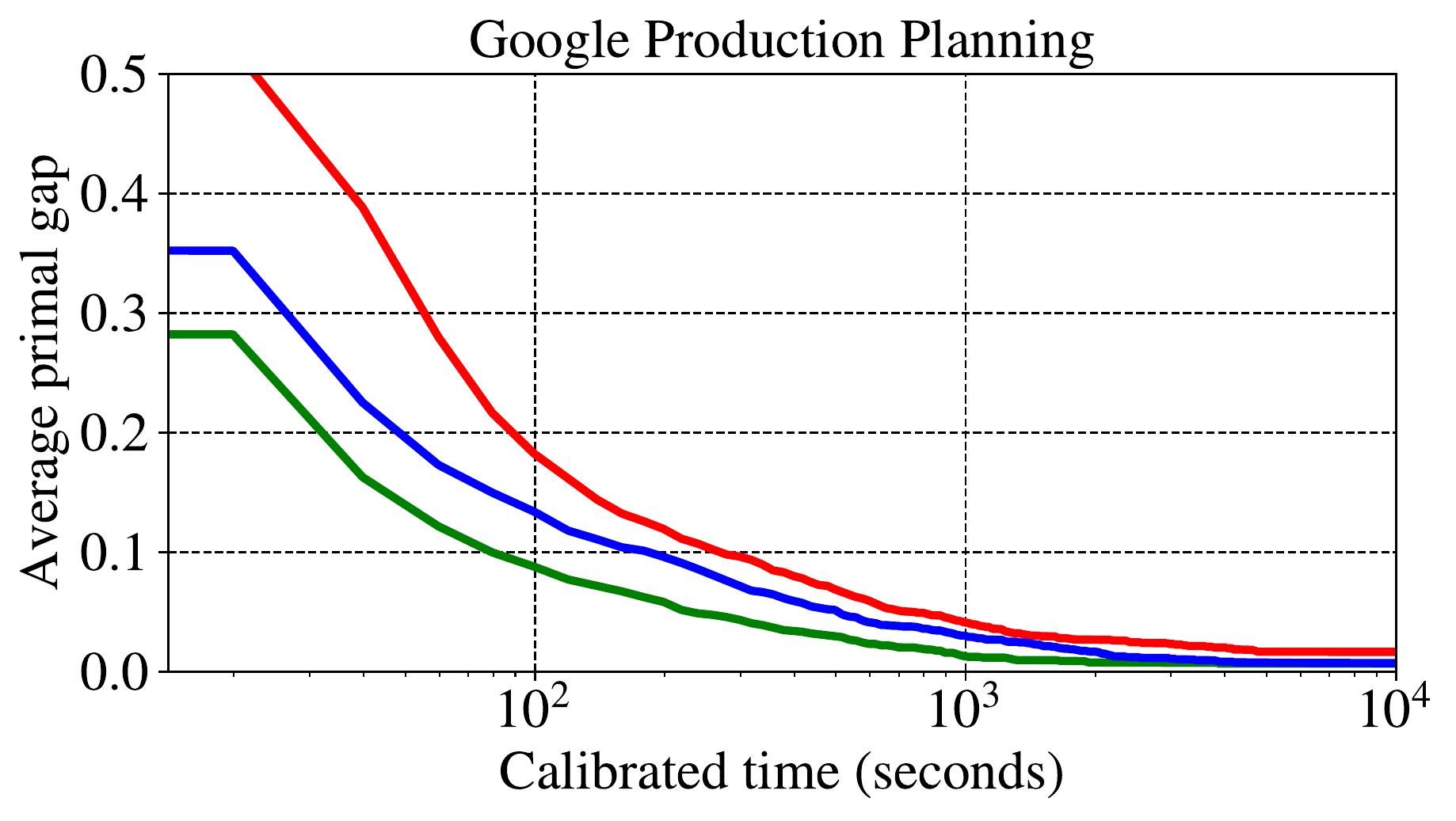} \\
        \includegraphics[width=0.50\textwidth]{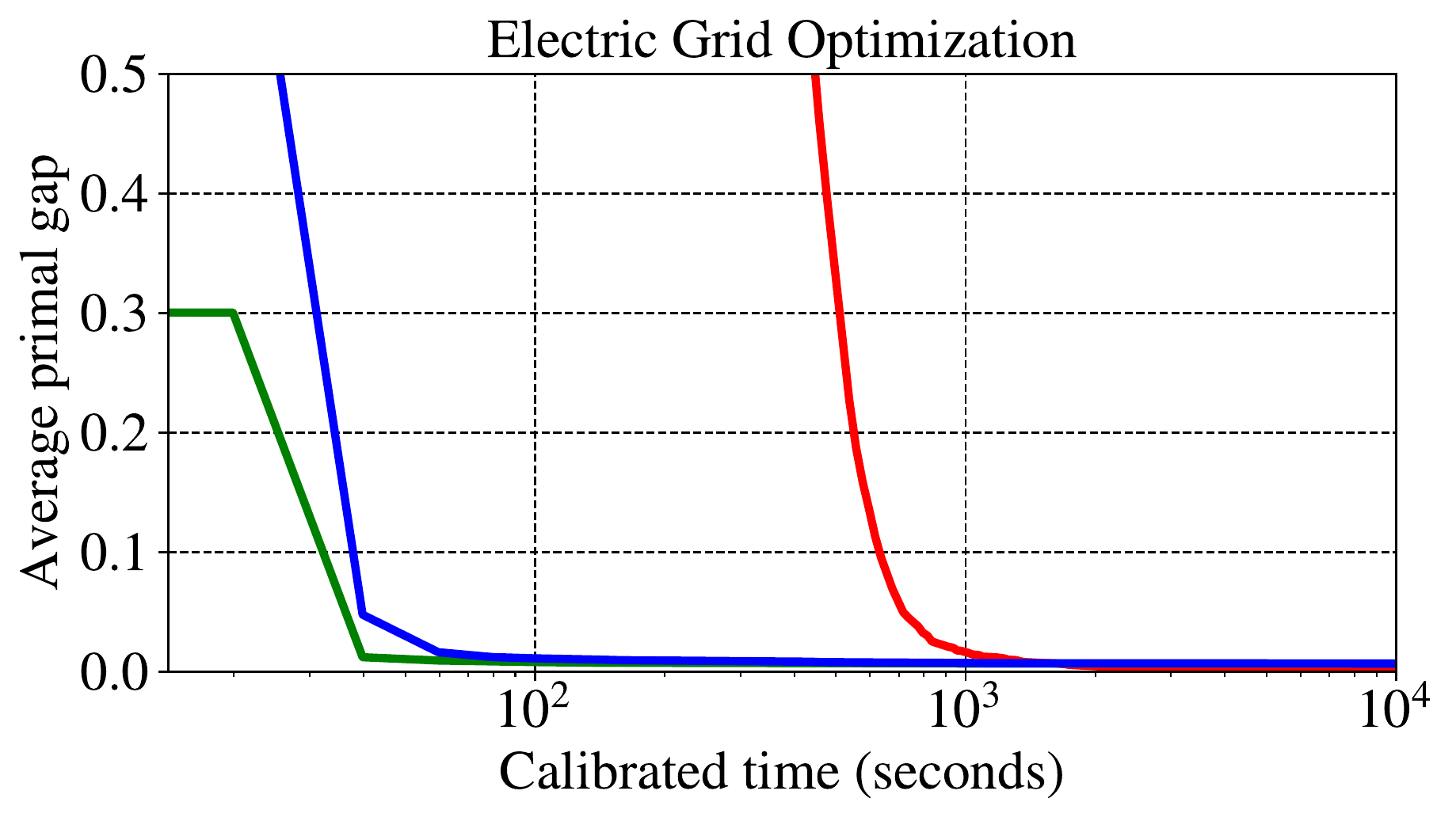} &
        \includegraphics[width=0.50\textwidth]{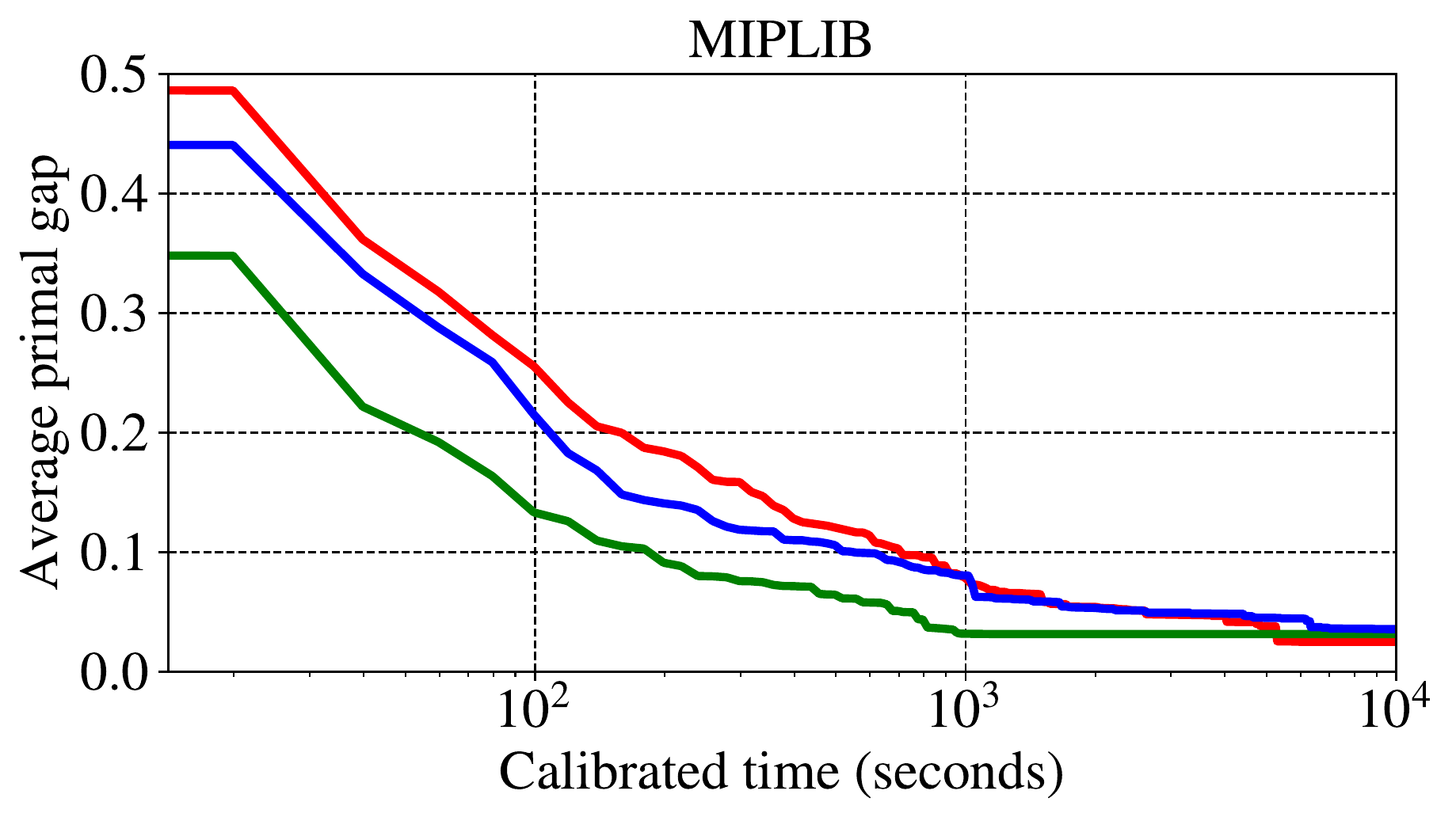} \\
    \end{tabular}
    \caption{Average primal gap between the primal bound and the best known objective value achieved by various algorithms as a function of running time on the benchmark datasets.}
    \label{fig:primal_bound_vs_calibrated_time}
\end{figure}

\begin{figure}[ht!]
    \begin{tabular}{cc}
        \includegraphics[width=0.50\textwidth]{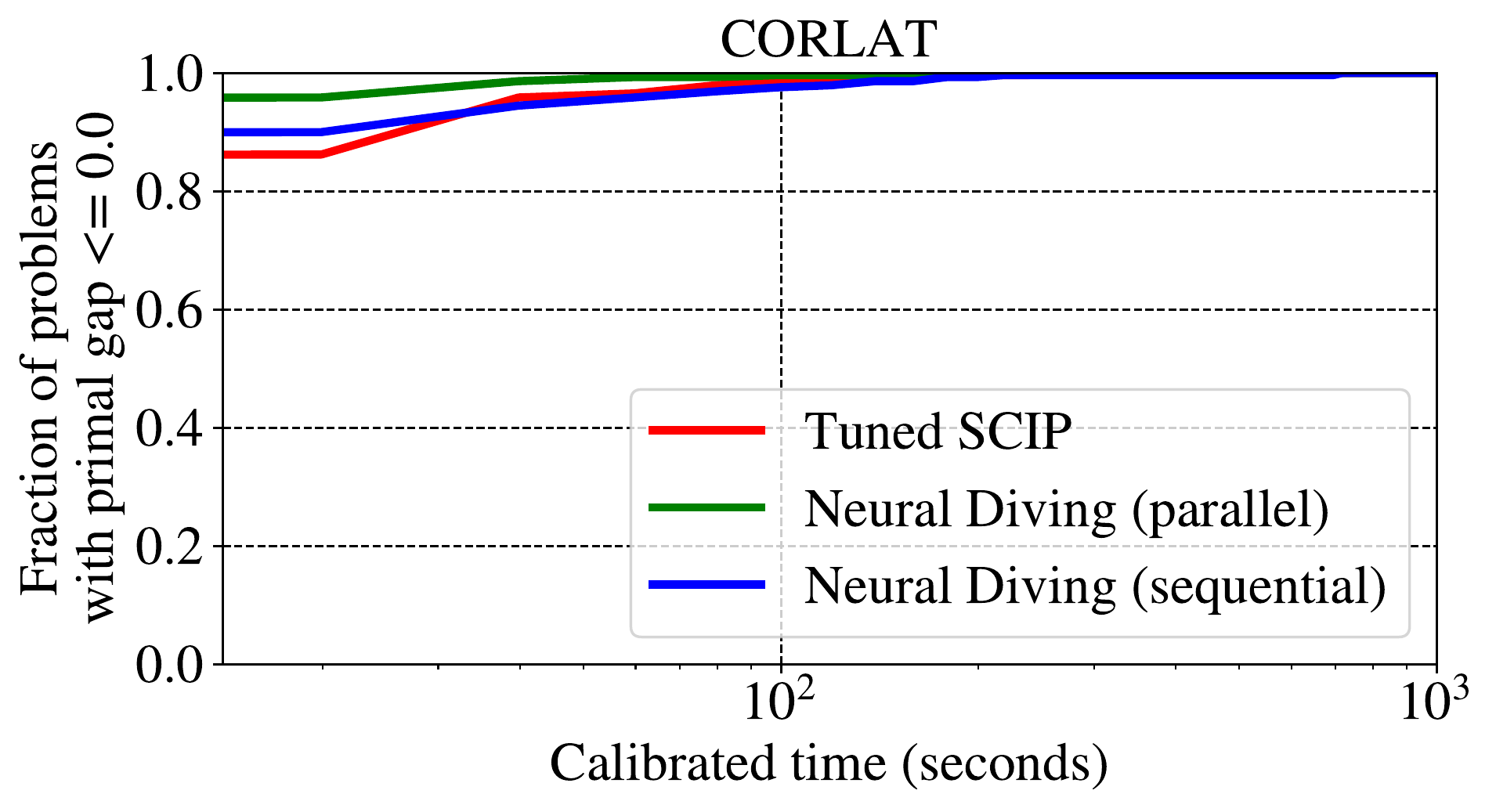} &
        \includegraphics[width=0.50\textwidth]{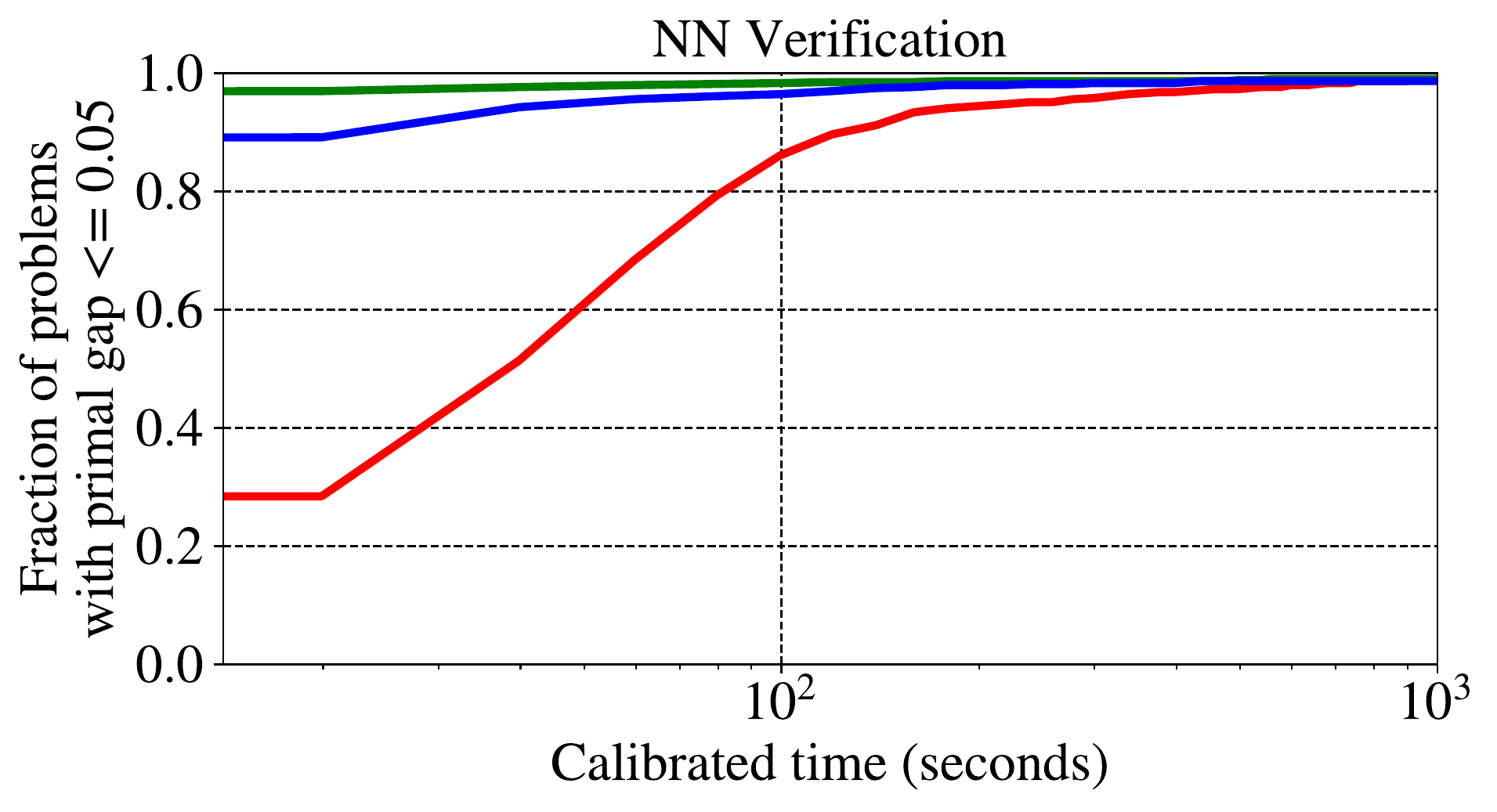} \\
        \includegraphics[width=0.50\textwidth]{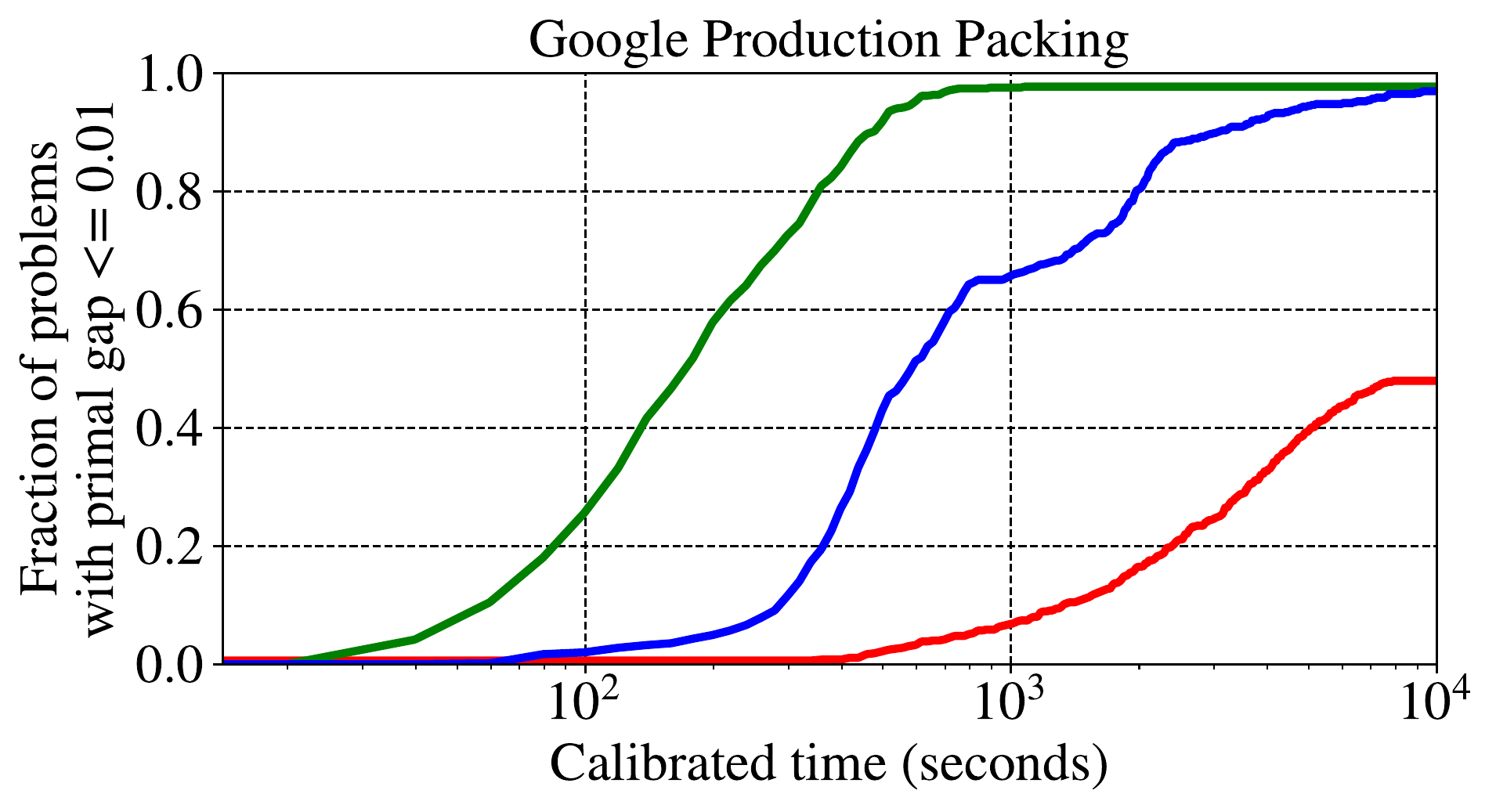} &
        \includegraphics[width=0.50\textwidth]{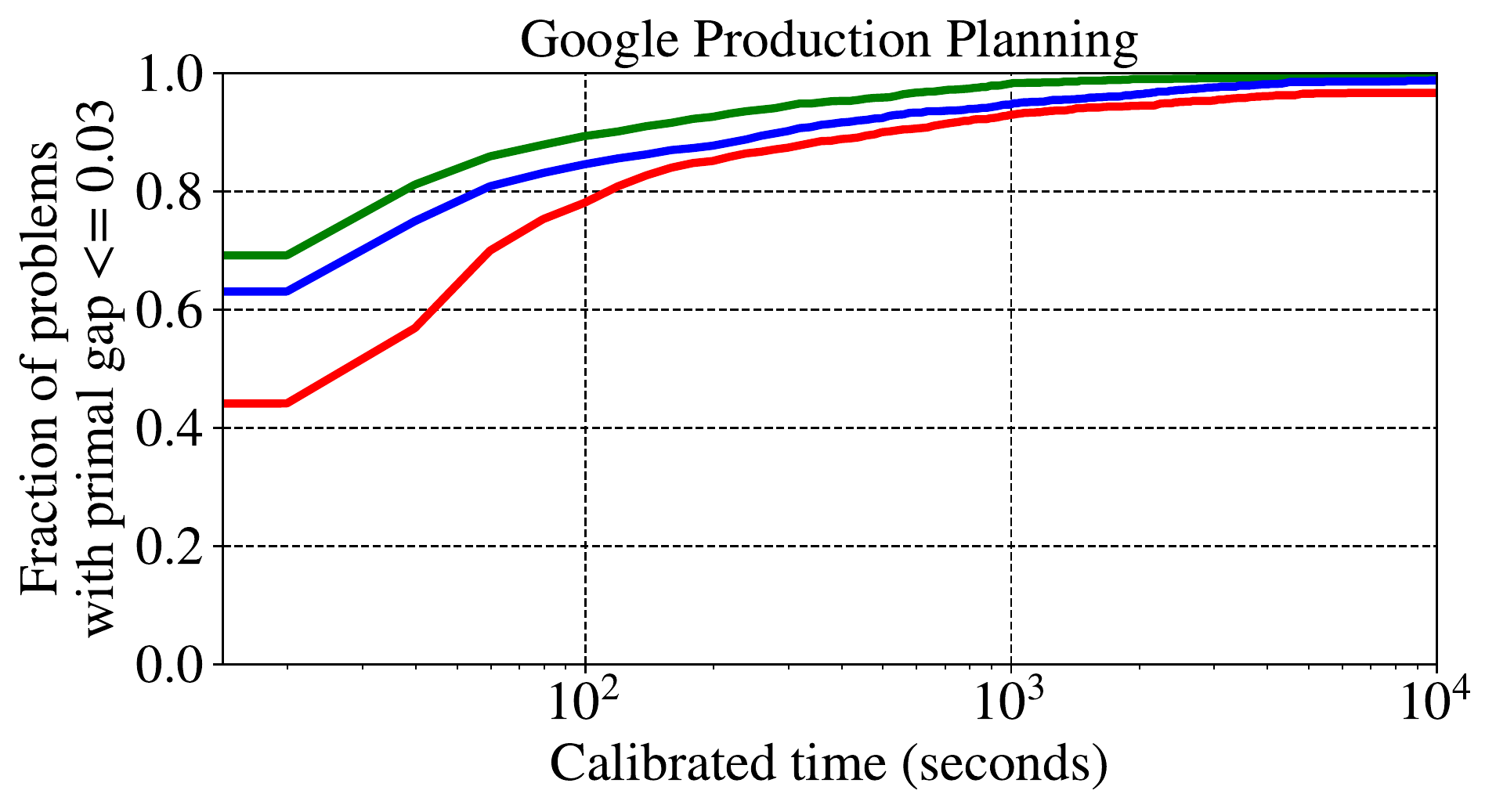} \\
        \includegraphics[width=0.50\textwidth]{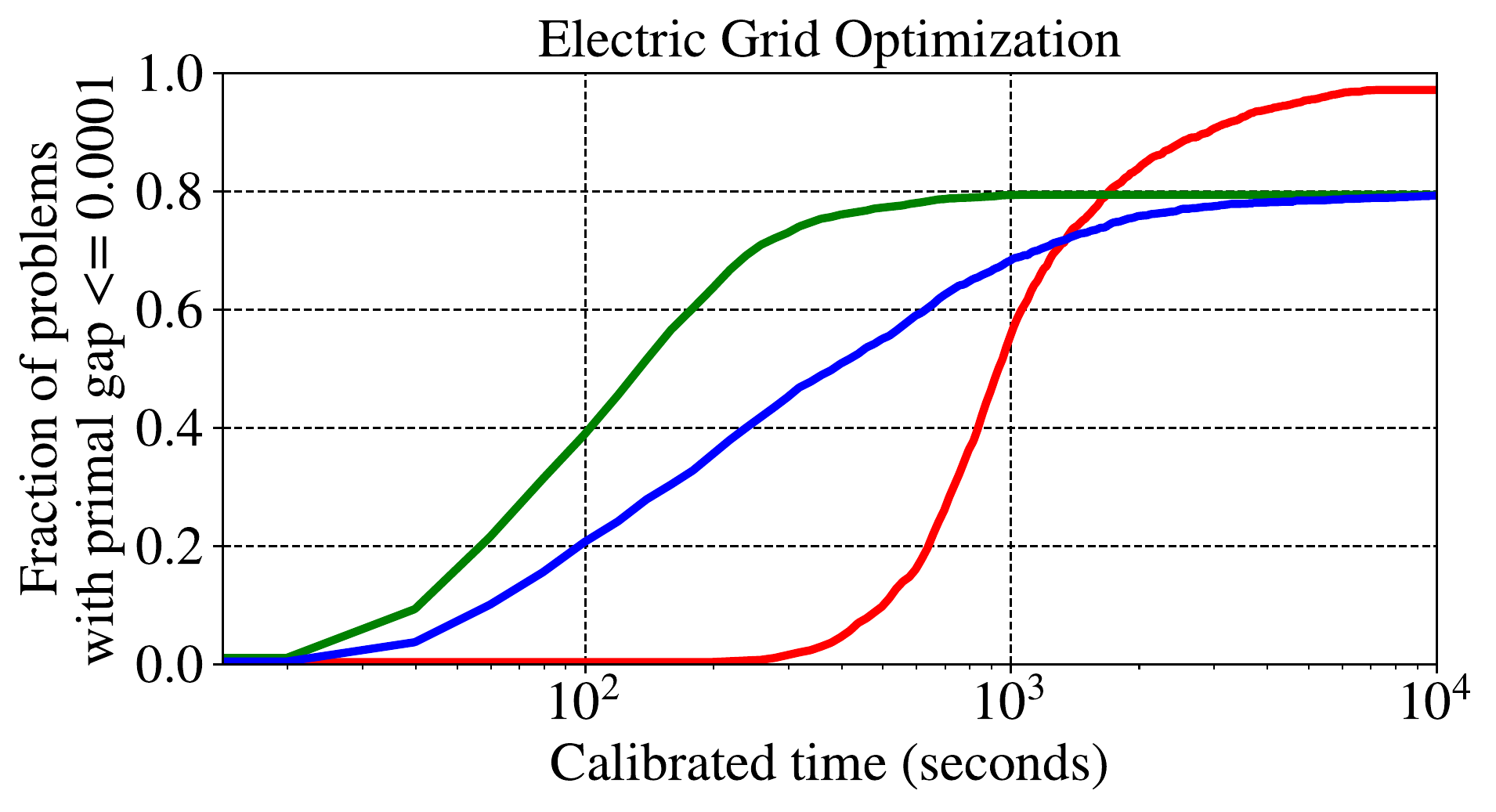} &
        \includegraphics[width=0.50\textwidth]{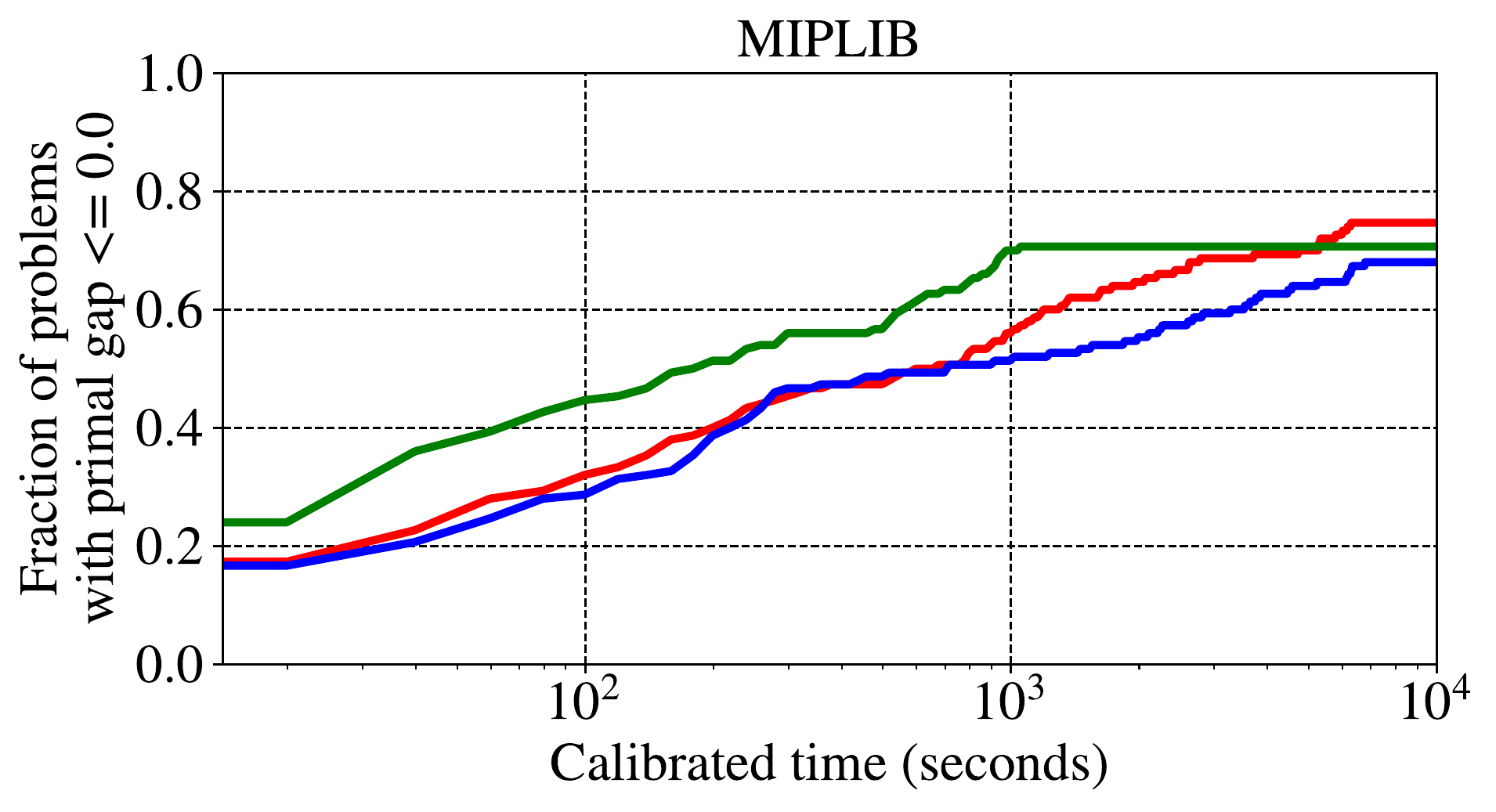} \\
    \end{tabular}
    \caption{Fraction of instances with primal gap less than or equal to the target gap achieved by various algorithms as a function of running time on the benchmark datasets.}
    \label{fig:primal_survival_vs_calibrated_time}
\end{figure}

\begin{figure}[ht!]
    \begin{tabular}{cc}
        \includegraphics[width=0.50\textwidth]{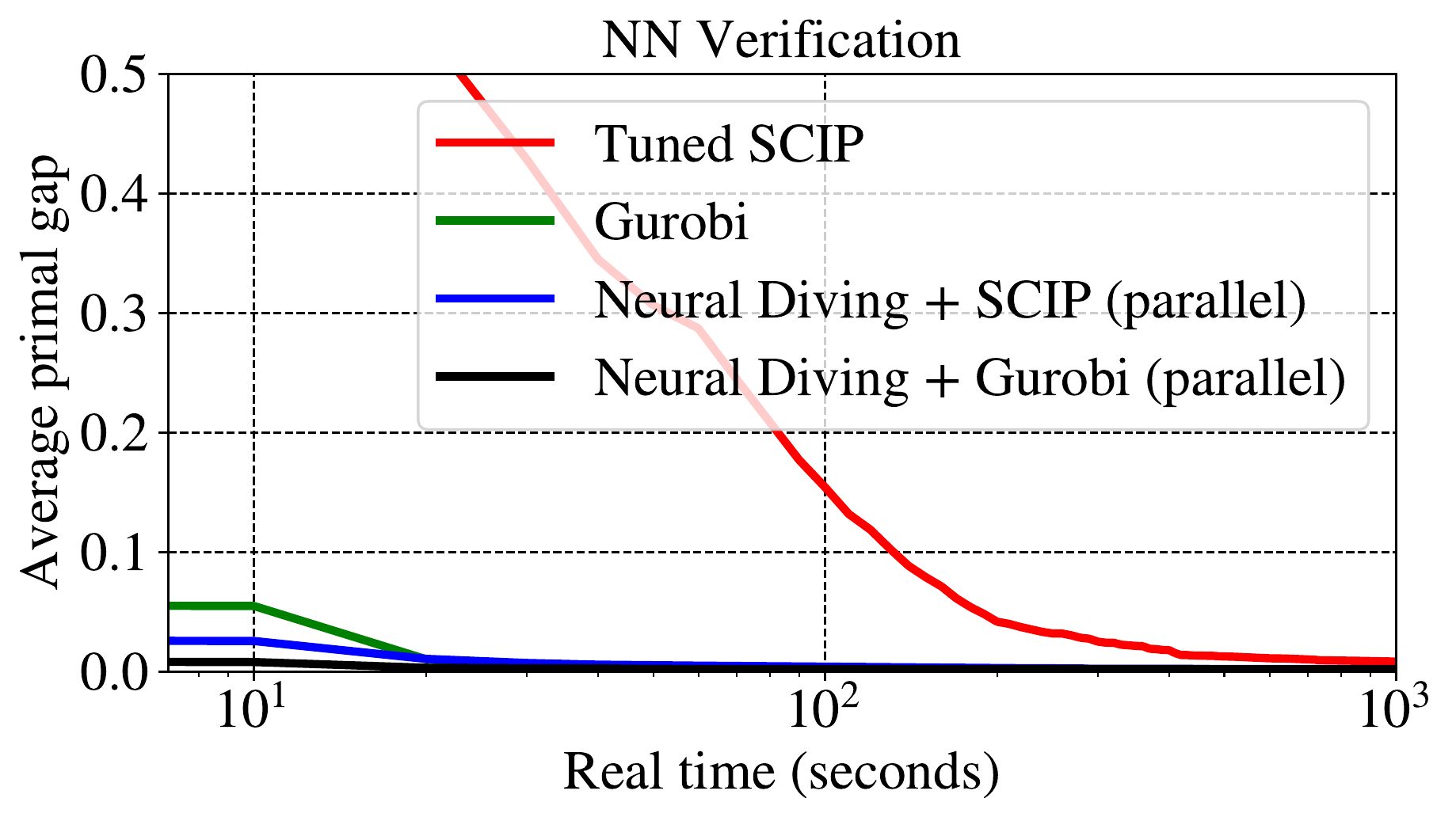} &
        \includegraphics[width=0.50\textwidth]{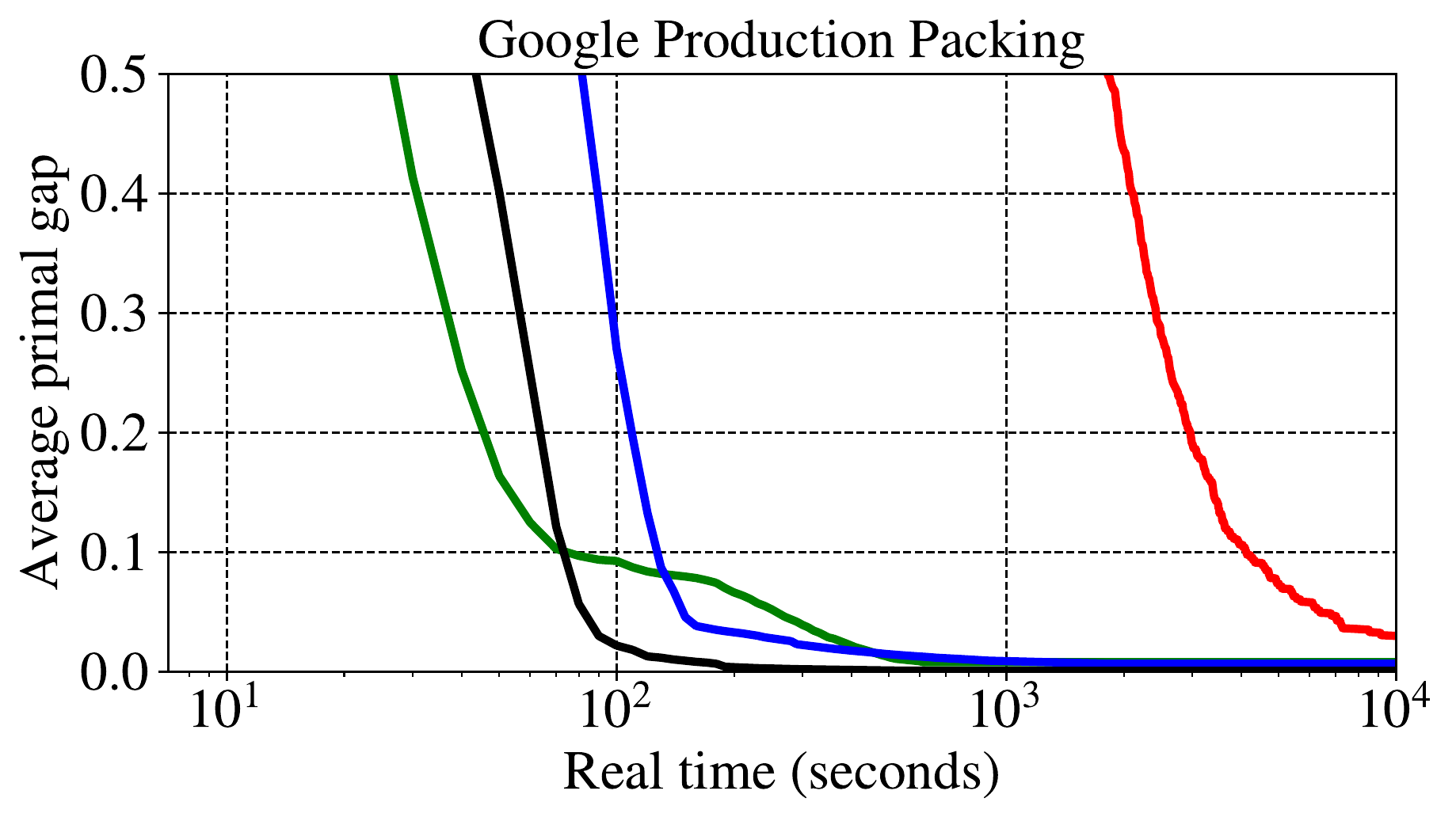}
    \end{tabular}
    \caption{Average primal gap between the primal bound and the best known objective value achieved by various algorithms as a function of running time on the benchmark datasets.}
    \label{fig:primal_bound_vs_calibrated_time_gurobi}
\end{figure}

\begin{figure}[ht!]
    \begin{tabular}{cc}
        \includegraphics[width=0.50\textwidth]{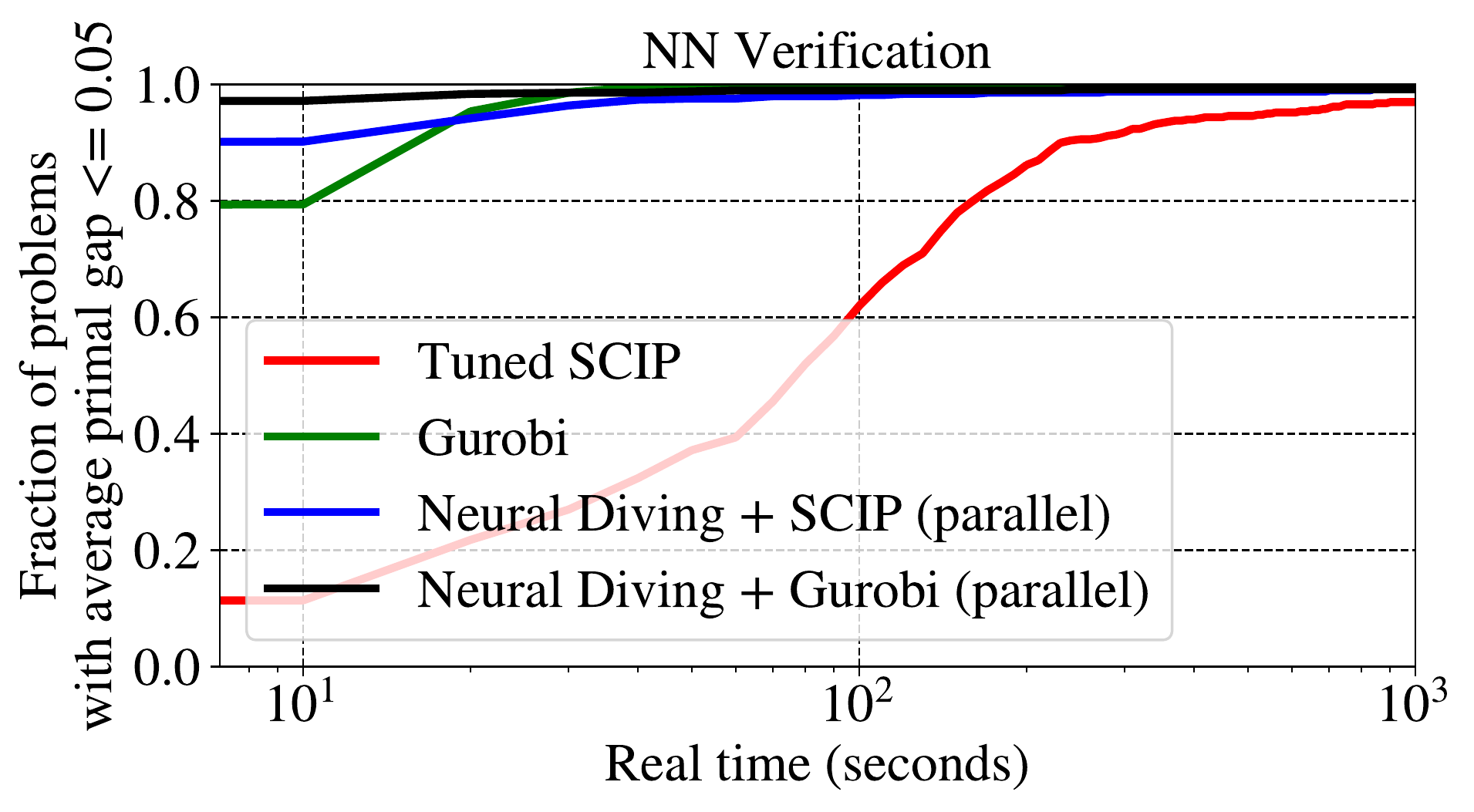} &
        \includegraphics[width=0.50\textwidth]{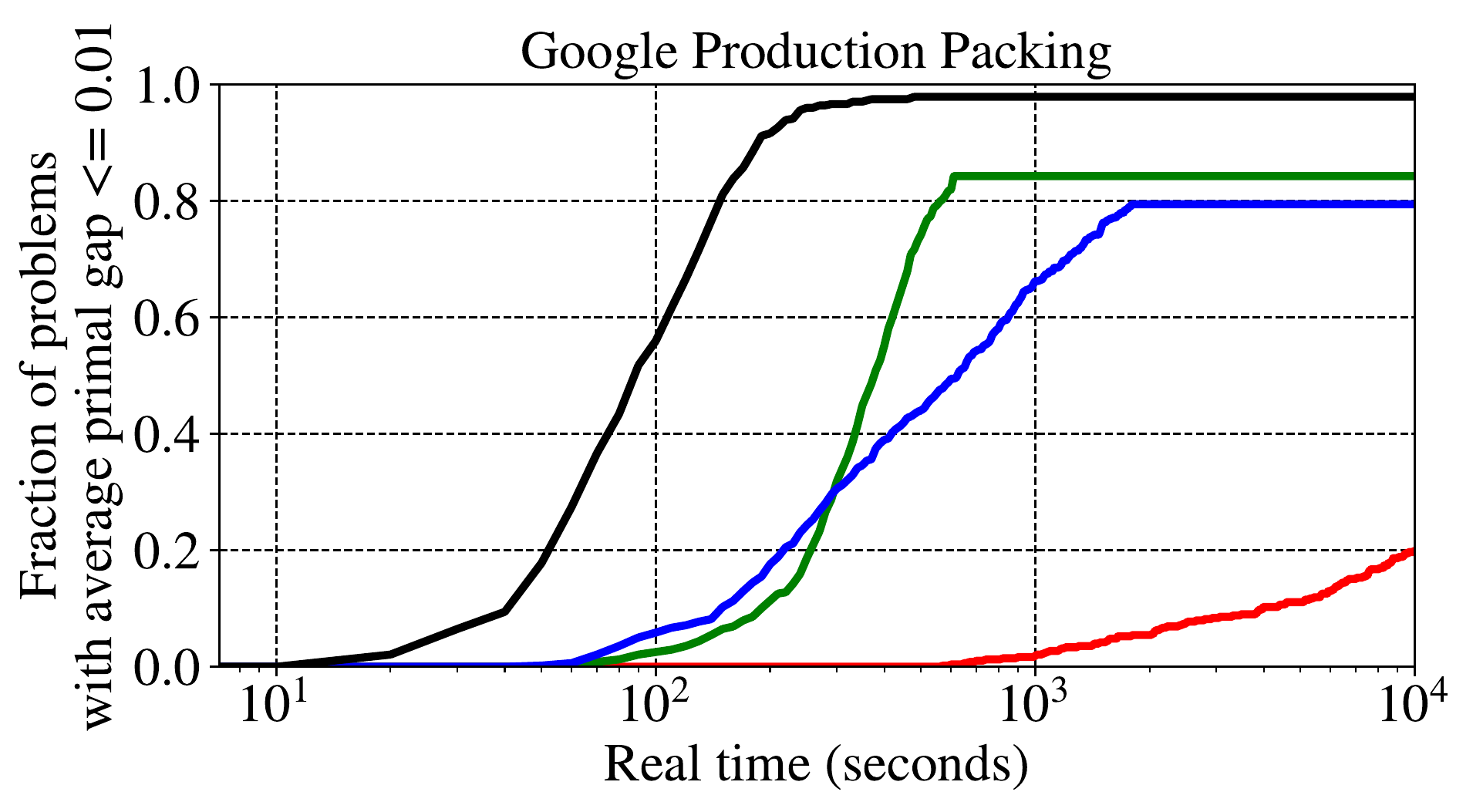}
    \end{tabular}
    \caption{Fraction of instances with primal gap less than or equal to the target gap achieved by various algorithms as a function of running time on the benchmark datasets.}
    \label{fig:primal_survival_vs_calibrated_time_gurobi}
\end{figure}

\subsubsection{Generating partial assignments}
\label{subsec:partial_assigments}
As described above, to produce partial assignments to an input MIP, we use the output of our generative models to decide both which values each integer variable should take, as well as which variables to assign: The acceptance head decides which variables are fixed or tightened, while the prediction head provides the variable values. Algorithm~\ref{alg:generate-submip} describes this procedure.

\begin{algorithm}[ht]
  \textbf{Input:} Learned distributions $p_\theta(x|M)$, $p_\theta(y|M)$, MIP instance $M$\\
  \textbf{Output:} Variable assignment and bound tightenings \\
  \smallskip

  $\mathrm{assignment} := \{\}$\\
  $\mathrm{tightenings} :=\{\}$\\

  \For{$x_i \in \mathrm{Variables}(M)$}{
    \If{$x_i$ is binary variable}{
        Sample $p_x$ from $p_\theta(x_i|M)$\\
        Sample $y$ from $p_\theta(y_i|M)$\\
      \If{$y = 1$}{
        Add ($x_i := \mathrm{round}(p_x)$) to $\mathrm{assignment}$\\
      }
    }
    \If{$x_i$ is non-binary integer variable}{
      $lb := $ lower bound of $x_i$\\
      $ub := $ upper bound of $x_i$\\
      $b_0, ..., b_k := $ binary representation of $(ub - lb)$, $b_0$ being most significant bit \\
      \For{$j\in \{0, ..., k\}$}{
        Sample $p_x$ from $p_\theta(x_{i, j}|M)$\\
        Sample $y$ from $p_\theta(y_{i, j}|M)$\\
        \eIf{$y = 1$}{
            \eIf{$\mathrm{round}(p_x) = 1$}{
                $lb := lb + \lceil(ub-lb)/2\rceil$ \\
            }{
                $ub := lb + \lfloor(ub-lb)/2\rfloor$ \\
            }
        }{
            Add $(lb \leq x_i \leq ub)$ to $\mathrm{tightenings}$\\
            \textbf{break} \\
        }
      }
    }
  }
    \textbf{return} $\mathrm{assignment}, \mathrm{tightenings}$ \\

\caption{Generating variable assignments and tightenings}\label{alg:generate-submip}
\end{algorithm}

We can repeat the above procedure many times, by obtaining multiple samples from a model, and also by using multiple models trained with different choices for the coverage threshold for the acceptance decision, to produce a diverse set of sub-MIPs. In experiments with the conditionally independent model we use one sample per model, and generate different sub-MIPs by using models trained with different coverage thresholds for $y$.

\subsubsection{Sequential and parallel sub-MIPs solving}

The sub-MIPs that we generate in the previous step can be solved fully independently from each other, allowing us to search for solutions for all sub-MIPs in parallel. In this parallel setting, we generate up to 100 sub-MIPs by applying combinations of different models and SCIP random seeds, and distribute each sub-MIP onto their own dedicated machine. Since each sub-MIP in this setting is solved in parallel, we report results in terms of (calibrated) wall clock time, meaning that we take results from the best performing sub-MIP at any given time.

While the parallel setting showcases the advantage that our approach has when distributed compute resources are available, we also set up a comparison that controls for the amount of computational resources used. That is, in addition to the parallel setting, we also compare our approach to SCIP when both are run sequentially on a single machine. In this setting, we are still generating multiple sub-MIPs as before, but process them sequentially in a random order, instead of in parallel.

%% file: supermip_results.tex
\subsection{Results}
\label{subsec:supermip_results}
We trained a GNN for every dataset on the training set and tuned hyperparameters (number of layers, width of layers, learning rate, exponential decay steps and assignment coverage thresholds) on the validation set. We report results on average primal gap vs our baseline SCIP in Figure \ref{fig:primal_bound_vs_calibrated_time}. Our parallel and sequential runs overall produce better primal bounds in shorter time on all datasets, compared to SCIP. We believe that the strength of our approach is in quickly finding good solutions, but it sometimes fails to find an optimal or near-optimal solution. This can can be seen e.g. on the survival plots in figure \ref{fig:primal_survival_vs_calibrated_time} where Neural Diving approach wins at shorter time limits, but loses to SCIP at the end on the Electric Grid Optimization and MIPLIB datasets. Note that SCIP's internal definition of gap is different from the one used in this paper, so in all experiments we set the relative gap setting to 0 in order to avoid stopping prematurely.

We include additional results on two datasets (Google Production Packing and NN Verification, see figures \ref{fig:primal_bound_vs_calibrated_time_gurobi} and \ref{fig:primal_survival_vs_calibrated_time_gurobi}) where we combine Neural Diving with the Gurobi solver \citep{gurobi2020}: we assign variables in the same way, but use Gurobi instead of SCIP to solve the remaining problem. The reported results are in non-calibrated time as we do not have access to the Gurobi internals and cannot report the calibrated time during solving.

\subsubsection{Finding incumbents to open MIPLIB instances}

In addition to evaluating Neural Diving on the datasets from Table~\ref{tab:datasets}, we also applied a modification of our approach on solving \emph{open} instances from the MIPLIB 2017 Collection Set. These instances are classified as open, because no optimal solution is known: Either there is no feasible solution known; or there exists a known feasible solution (incumbent), but it does not match any known dual bound. In this sense, these instances are the hardest problems that the MIPLIB 2017 Collection Set offers.

Apart from their inherent hardness, these instances also impose an additional challenge to our learning-based approach: The set of instances is almost entirely heterogeneous, in that there are only a handful of related instances available per problem, if at all. Our approach of learning a domain-specific primal heuristic hence needed to be adapted. To tackle this challenge, we built a training set for each attempted open instance which should represent the problem distribution of the target instance. More specifically, for a given target instance, we generated up to 500 similar instances through randomly applying a mix of the following manipulations: (1) Randomly dropping constraints of the target instance, and (2) fixing a random subset of variables to the previous incumbent. To generate each such instance, we applied 10 iterations of each step (1) and (2) on the original target instance.

After creating such a training set for each target instance, we apply our approach as described in earlier sections: We trained a separate GNN for each target instance, using the generated training set, with label solutions found on the training set using SCIP. After training, we used the corresponding GNN to produce up to 10 partial assignments to the target instance (as described in Section~\ref{subsec:partial_assigments}), and used Gurobi to solve the sub-MIPs defined by the partial assignments. For three open instances, this approach yielded new incumbents that were not known before. Table~\ref{tab:open_miplib_results} shows the improved objective values compared the previous best known solution.

\begin{table}
\caption{Incumbents found for open MIPLIB 2017 instances by Neural Diving. Lower objective value is better.}
\begin{adjustbox}{center}
{\small
    \begin{tabular}{cccc}
    \toprule
        Instance & Neural Diving & Previous best \\
         & objective value & objective value\\
        \midrule
        \emph{milo-v12-6-r1-75-1} & \textbf{1153756.398} & 1153880\\
        \emph{neos-1420790} & \textbf{3121.29} & 3121.42\\
        \emph{xmas10-2} & \textbf{-497} & -495\\
    \bottomrule
    \end{tabular}
    }
\end{adjustbox}
\label{tab:open_miplib_results}
\end{table}

%% file: deepbrancher.tex
\section{Neural Branching}
\label{sec:deep_Branching}
A branch-and-bound procedure has two decisions to make every iteration: which leaf node to expand, and which 
variable to branch on. In this work we focus on the latter. The quality of the variable selection decision can have a large impact on the number of steps taken by branch-and-bound to solve a MIP \citep{achterberg2005branching, Glankwamdee2011LookaheadBF, schubert2017thesis, yang2019multivariablebranching}. We use a deep neural network to learn a variable selection policy by imitating the actions of a node-efficient, but computationally expensive, expert. By distilling the policy of the expensive expert into the neural network we seek to maintain approximately the same decision quality but substantially reduce the time taken to make the decision. The decision at a given tree node is entirely local to that node, so a learned policy only needs to have a representation of the node as input, rather than the entire tree, which makes it more scalable.

\subsection{Expert policy}
We would like an expert policy that solves MIPs by building small branch-and-bound trees, since such an expert is making good branching decisions. Among the many branching policies proposed in the optimization literature (see, e.g., \cite{achterberg2005branching}) none provably achieve the smallest trees, but empirically \emph{full strong branching} (FSB) \citep{achterberg2005branching} tends to use fewer steps than competing approaches. It performs one-step lookahead search by simulating one branching step for \emph{all} candidate variables and picking the one that provides the largest improvement in the dual bound, as determined by the solution to the linear program at the new node. (A generalization to $k$-step lookahead search can be found in \cite{Glankwamdee2011LookaheadBF,schubert2017thesis}.) This requires solving $2 \times n_{\text{cands}}$ linear programs, where $n_{\text{cands}}$ is the number of possible branching candidates at a given step of branch-and-bound. Concretely, denote by $x^\star$ the solution to the LP relaxation of (\ref{eqn:mip1}) at the current node we are branching from with variable constraint vectors $l$ and $u$, \ie, $l \leq x^\star \leq u$. Then for each variable candidate $i$, FSB needs to solve the following two LPs:
\begin{equation}
\label{eq:up_and_down_lps}
\begin{array}{lr}
\begin{array}{ll}
\mbox{minimize} & c^\top x^\mathrm{up}\\
\mbox{subject to} & Ax^\mathrm{up} \leq b\\
& l^{(i)} \leq x^\mathrm{up} \leq u
\end{array}
\begin{array}{ll}
\mbox{minimize} & c^\top x^\mathrm{down}\\
\mbox{subject to} & Ax^\mathrm{down} \leq b\\
& l \leq x^\mathrm{down} \leq  u^{(i)}
\end{array}
\end{array}
\end{equation}
over variables $x^\mathrm{up} \in \R^n$ and $x^\mathrm{down} \in \R^n$, where $l_i^{(i)} = \ceil{x^\star_i}$ and $u_i^{(i)} = \floor{x^\star_i}$, and other entries of $l^{(i)}$ and $u^{(i)}$ are unchanged from their values at the current node, $l$ and $u$. It combines the optimal values of these LPs into a score for that candidate, and it uses those scores to decide which variable to branch on.

For practical MIP solving FSB's per-step computational cost is often so high that it can be used only for a few branch-and-bound steps before the running time becomes prohibitively high. SCIP's default variable selection policy, \emph{reliable pseudocost branching} (RPB) \citep{achterberg2005branching}, uses FSB for a small number of steps at the start of branch-and-bound and switches to a heuristic with lower per-step cost for the subsequent steps. But an expert policy used as the target for imitation learning can afford a higher per-step cost. The expert is only used to generate training data offline as a one-time expense by running it on a training set of MIPs and recording its inputs and outputs at each step. Even then, on large-scale instances (e.g., with $> 10^5$ variables), offline data generation can be prohibitively slow, which prevents learning approaches from scaling to such instances. Therefore we speed up the expert policy by exploiting GPUs, as explained next.

\subsubsection{ADMM batch LP solving}
 SCIP by default uses the LP solver SOPLEX \citep{gamrath2020scip} based on the simplex algorithm \citep{dantzig1998linear} for solving the linear programs arising from FSB. 
 Although in principle the simplex algorithm can be parallelized, it is not easy to \emph{batch} the computation
 together for several LPs that are closely related (and therefore leverage a GPU), and consequently SOPLEX simply solves the LPs sequentially on CPU.
 In order to scale to a large number of candidate variables at a reasonable cost we wrote our own batch LP solver using the alternating directions method of multipliers (ADMM) \citep{boyd2011distributed}, based on the Splitting Conic Solver (SCS) algorithm for cone programs \citep{ocpb:16, o2020operator}. This solver can simultaneously handle many LPs in a batch that we process on GPU. This can (approximately) solve all the LPs in full strong branching substantially faster than running SOPLEX sequentially. Moreover, since we are using this solver merely to produce training data, we can tune the algorithm parameters with the primary goal of producing smaller search trees, even at the cost of some data generation time. 
It is the output produced by this ADMM-based policy that we will train a graph neural network to imitate. The graph neural network is able to produce actions that are close enough to the policy of the ADMM expert to produce small search trees, but at a fraction of the computational cost. We provide full details in the Appendix, section~\ref{subsec:appendix_admm_batch_lp_solver}.

\begin{figure}[t]
    \centering
    \begin{tabular}{cc}
        \includegraphics[width=0.48\textwidth]{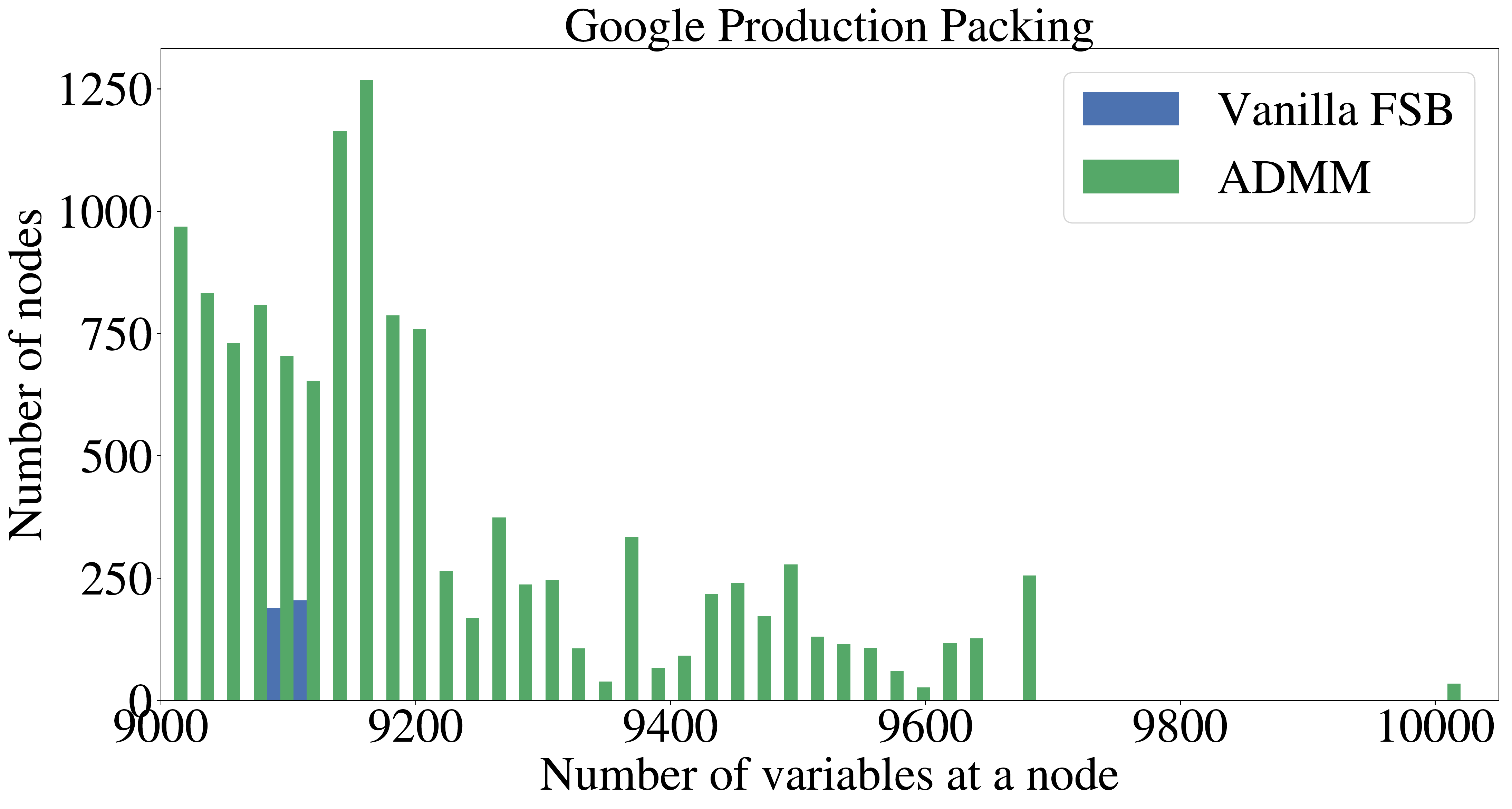}
        \includegraphics[width=0.48\textwidth]{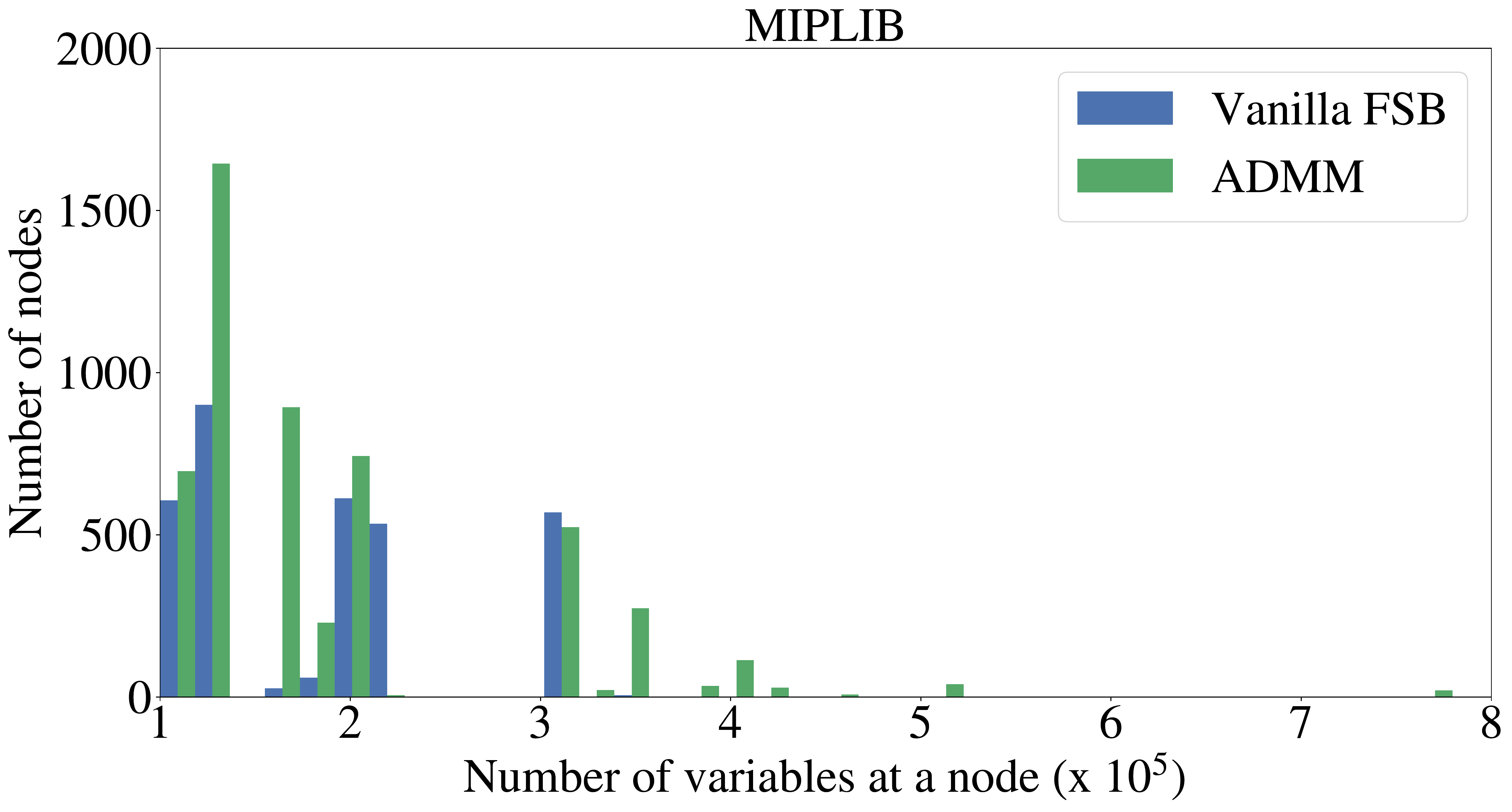}
    \end{tabular}
    \caption{Histogram of the number of imitation learning examples, where each example is a node in a branch-and-bound search tree, generated by the ADMM expert and the Vanilla Full Strong Branching expert on the training sets for Google Production Packing and MIPLIB in the same time limit.}
    \label{fig:admm_vs_vfsb_miplib_data_generation}
\end{figure}

To demonstrate the improved scalability of the expert data generation, we compare the ADMM expert to the expert used by \cite{gasse2019exact} called Vanilla Full Strong Branching (VFSB) on Google Production Packing and MIPLIB. VFSB is a modified version of SCIP's FSB implementation that is still performed sequentially but disables its `side effects', \ie, changes to the solver state not directly due to a branching decision itself, but due to information gained during the lookahead search over branching candidates.
Figure \ref{fig:admm_vs_vfsb_miplib_data_generation} shows a histogram of the number of variables in each branch-and-bound node generated by the two experts within $60$ hours,
on nodes with the largest number of variables for both datasets ($>9000$ variables for Google Production Packing and $>10^5$ variables for MIPLIB). The ADMM expert generates $31.7\times$ and $1.6\times$ more data on nodes with the largest number of variables for Google Production Packing and MIPLIB, respectively. Across the full range of the number of variables, the ADMM expert generates $12.1\times$ and $1.4\times$ more data on Google Production Packing and MIPLIB, respectively. The largest node that the VFSB expert manages to generate data from on MIPLIB has $3.4\times10^5$ variables, while the ADMM expert is able to generate data from nodes with almost $10^6$ variables.

\subsection{Imitation Learning}
Imitation learning is the name given to the broad family of algorithms that seek to learn the policy of an 
expert, given examples of expert behavior. In most cases, and in this paper, this is formulated as a supervised learning problem. We consider three variants: 1) cloning an expert policy \citep{pomerleau1989alvinn, bain1995framework}
2) distillation with random moves, and 3) DAgger \citep{ross2011reduction}. Distillation is the simplest of the three and it simply attempts to learn to predict the outputs of the expert at each node. Distillation is not robust to shifts away from the expert's state distribution. This is particularly problematic in the case of branching as mistakes near the root can cause the subsequent nodes visited in the search tree to be very different than those seen during training, since the expert never visited those nodes. Mixing random moves with expert moves is a simple heuristic that can alleviate this problem to some degree. When running the expert we take a random action with probability $10\%$ at each node. We can run the expert many times on each training MIP generating slightly different data each time, thereby increasing the total amount of training data available. In all our experiments we did this five times for each MIP. DAgger is another, more involved, procedure where we train an agent using distillation, run it in a branch-and-bound procedure to make decisions, but also compute the expert outputs at each step to use as targets for learning. This produces one step of DAgger data upon which we can re-train the agent. In theory this procedure can be repeated using the new agent trained on the previous agent DAgger data, but in this manuscript we only consider a single iteration of DAgger.

In our experiments, we treat the choice of the three imitation learning variants as a hyperparameter to tune for each dataset. We generate data and train policies using all three variants for each dataset, and we select the policy that achieves the lowest average dual gap on the validation set instances in a 3 hour time limit. The selected policy is then evaluated on the test set to report results.

\subsection{Implementation details}
The ADMM expert solves a pair of LPs for each branching candidate (the `up' and `down' branches). Let us
denote by $\mathrm{OPT}_i^\mathrm{up}$ and $\mathrm{OPT}_i^\mathrm{down}$ the (estimated) optimal
values of the LP for the up and down branches respectively for candidate $i$, and denote by $\mathrm{OPT}$ the objective
value of the LP at the node we are branching from. We combine these numbers into a single score
for each variable defined as
\begin{equation}\label{eqn:bnbvalue}
s_i = (\mathrm{OPT}_i^\mathrm{up} - \mathrm{OPT} + \epsilon)(\mathrm{OPT}_i^\mathrm{down} - \mathrm{OPT}  + \epsilon)
\end{equation}
where $\epsilon = 10^{-4}$ (this is essentially the product rule from \citet{achterberg2009scip}). Given a set $\mathcal{C}$ of candidate variables for branching, we convert their scores into a categorical distribution over the candidates using
\[
p^\mathrm{expert}_i = \frac{s_i}{\sum_{c \in \mathcal{C}} s_c}.
\]
On the training set of MIPs we ran the ADMM expert, with random moves or initialized using DAgger as appropriate,
and logged the expert scores taken at each node (we did not log the action if it was selected randomly or 
selected by the agent in the DAgger setting - only the ADMM expert action was ever logged). This provides our training dataset for the imitation learning agent. 

The neural network policy for variable selection has the same graph convolutional network architecture described in section~\ref{subsec:architecture}. It conditions on a leaf node of the branch-and-bound tree selected by SCIP to branch on. As the node itself is a MIP, it can be represented using the bipartite graph representation described in section~\ref{subsec:mip_bipartite_graph}. Its output is a categorical distribution over the set of candidates $\mathcal{C}$. Let $M$ be the MIP bipartite graph representation of the leaf node, $v_c$ be the embedding computed by a graph convolutional network for the bipartite graph node for the branching candidate variable $x_c$, as described in section~\ref{subsec:architecture}, and $\phi$ be the learnable parameters of the policy. The probability $p_{\phi}(x_c|M)$ of selecting $x_c$ for branching at the node is given by:
\begin{align}
    t_c &= \text{MLP}(v_c; \phi),\\
    p_{\phi}(x_c|M) &= \frac{\exp(-t_c)}{\sum_{c'\in \mathcal{C}}\exp(-t_{c'})}.
\end{align}

When training, the neural network receives a batch of node features with the associated expert scores and produces a set of probabilities $p^\mathrm{graphnet}$ for each node in the batch by passing the activations of the last layer through a softmax. We used batch size 8 sampled uniformly at random from the training dataset. We tried a variety of losses and found that negative cross-entropy loss performed best on average, which is given by
\[
L(\phi) = \sum_{c \in \mathcal{C}} p_c^\mathrm{expert} \log p_{\phi}(x_c|M).
\]
The loss is minimized using ADAM \citep{kingma2014adam} with a learning rate of $10^{-4}$. Figure~\ref{fig:neural_Branching_target_policy_comparison} in the appendix shows that imitation learning succeeds to accurately approximate the expert policy, both on the training and test sets.

%% file: deepbrancher_results.tex
\subsection{Results}
\label{subsec:deepbrancher_results}

We evaluate the learned branching policy on the task of optimizing the dual bound. As explained in section \ref{sec:evaluation}, for comparison we use the dual gap $\gamma_d(t)$ with respect to a pre-computed best known primal bound $p^\star$ as a function of time. The gap is averaged over all MIPs in the test set of a dataset.

\begin{figure}[H]
    \begin{tabular}{cc}
        \includegraphics[width=0.48\textwidth]{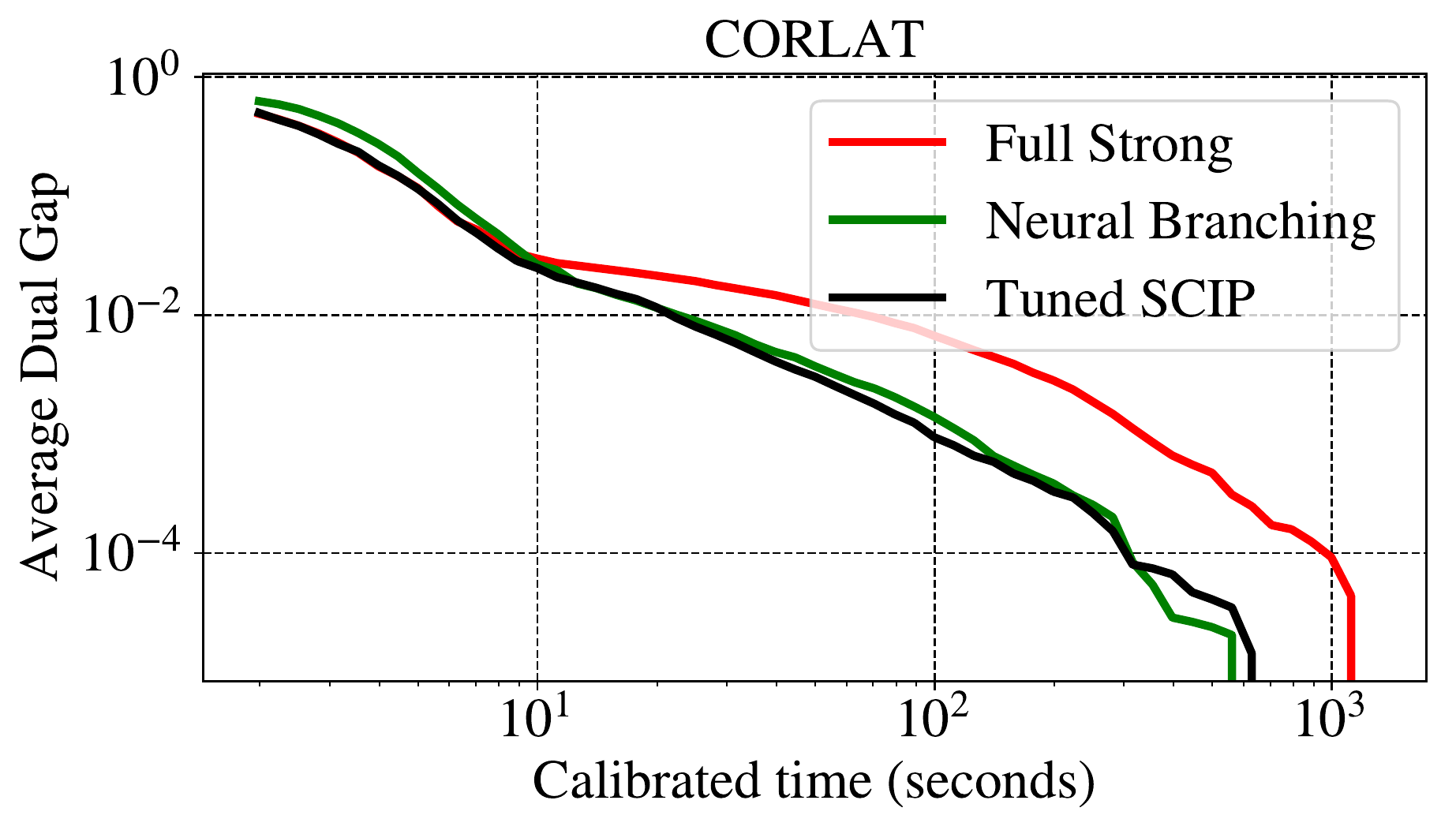} &
        \includegraphics[width=0.48\textwidth]{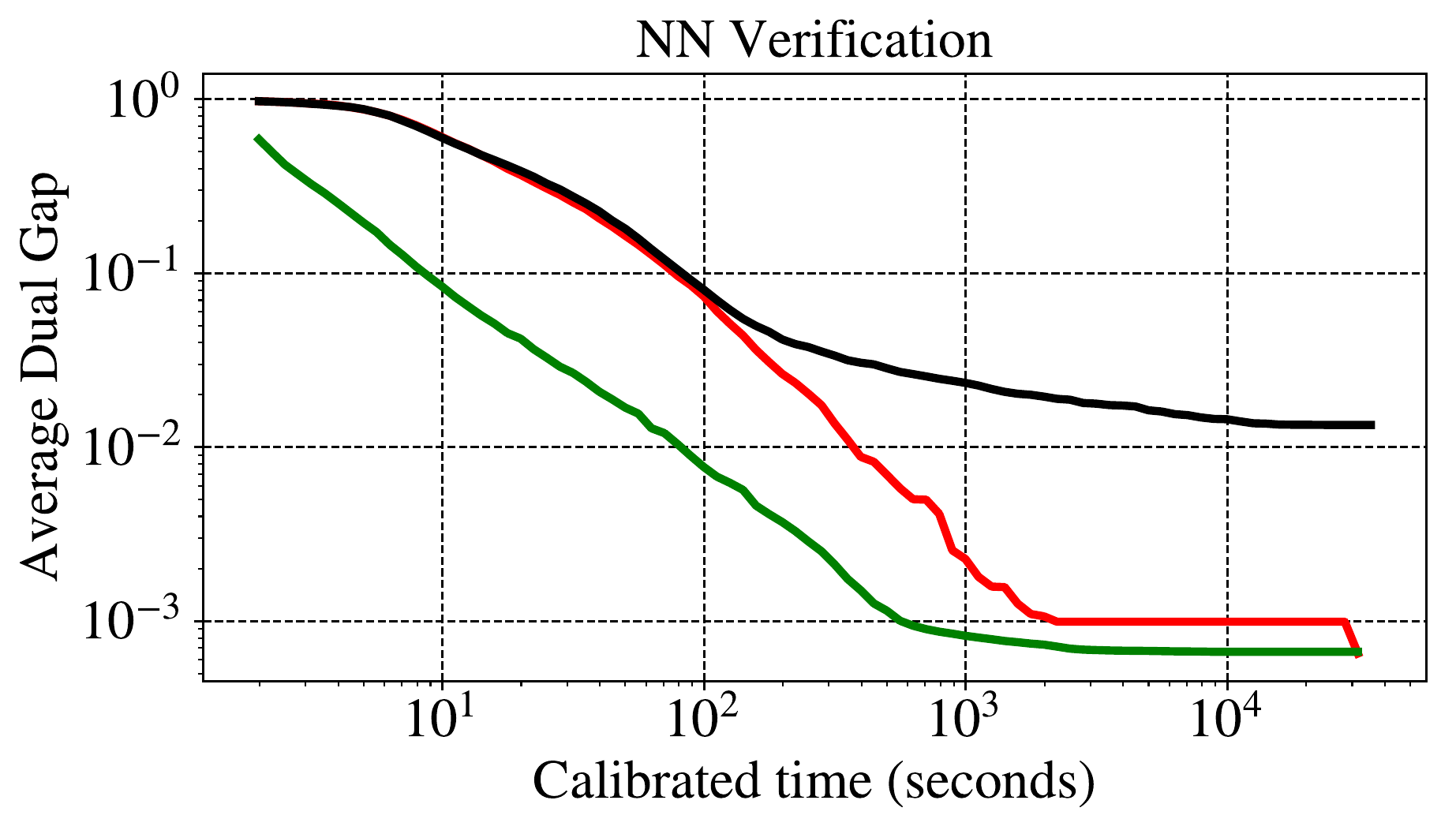} \\
        \includegraphics[width=0.48\textwidth]{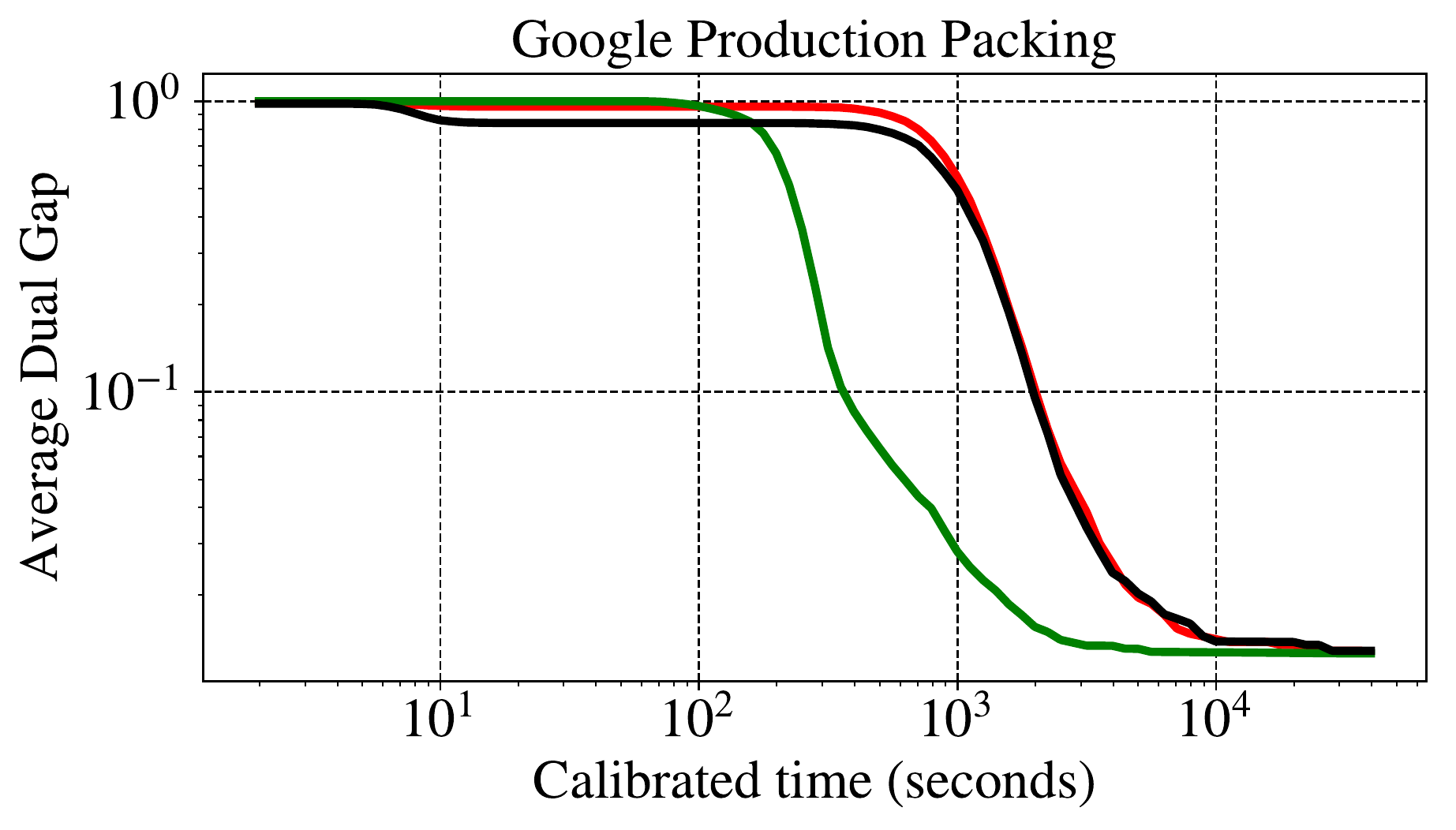} &
        \includegraphics[width=0.48\textwidth]{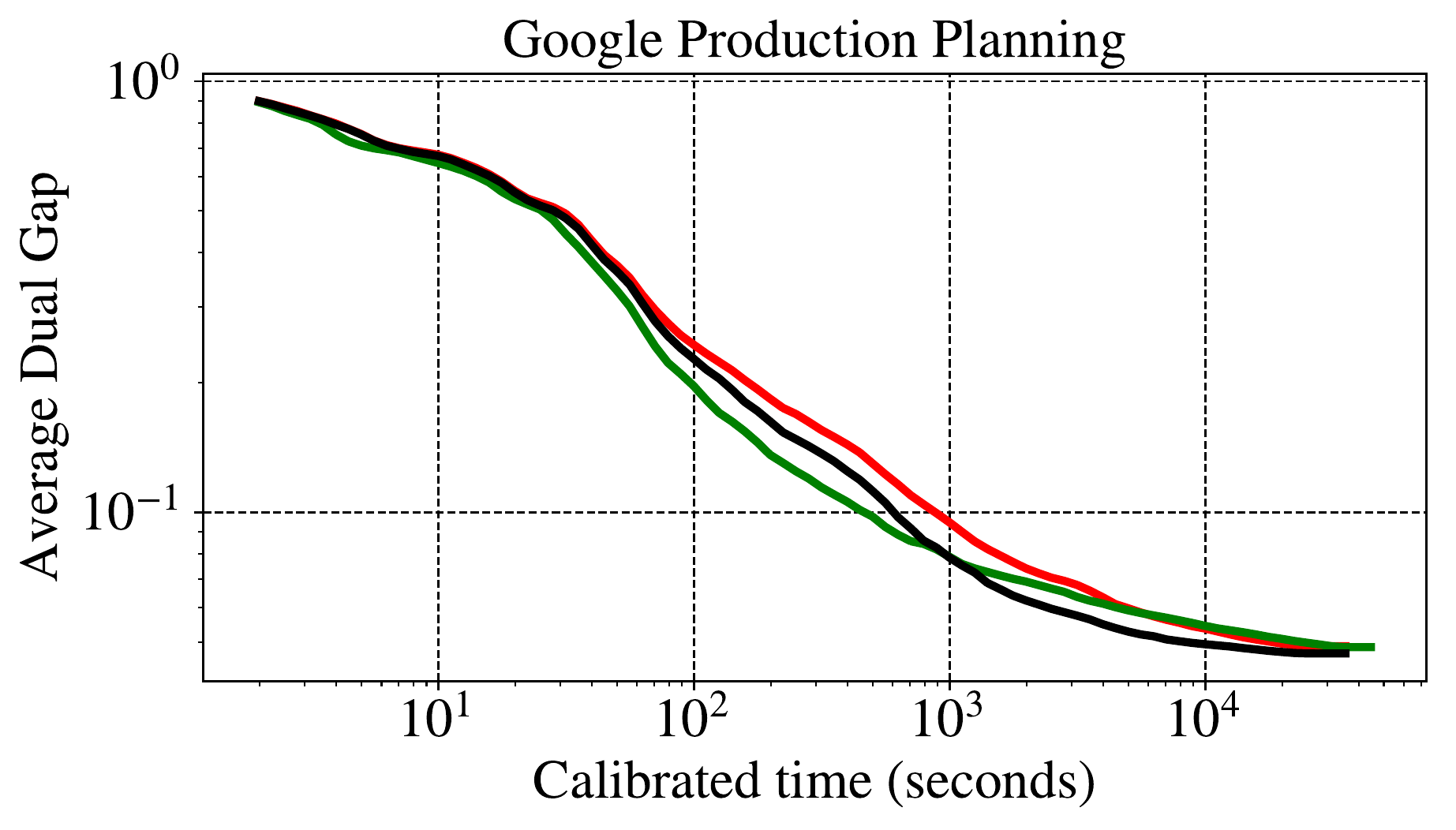} \\
        \includegraphics[width=0.48\textwidth]{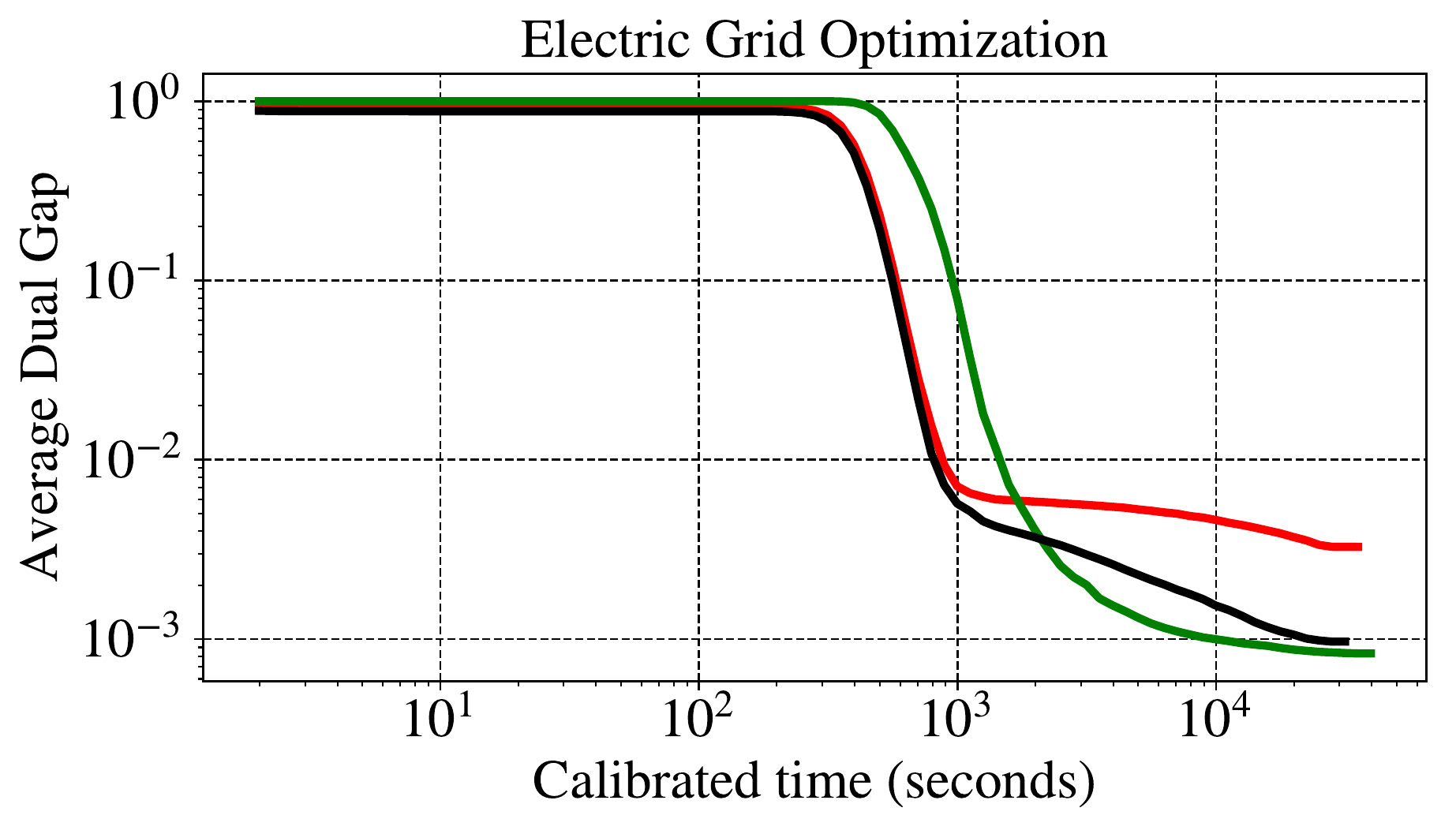} &
        \includegraphics[width=0.48\textwidth]{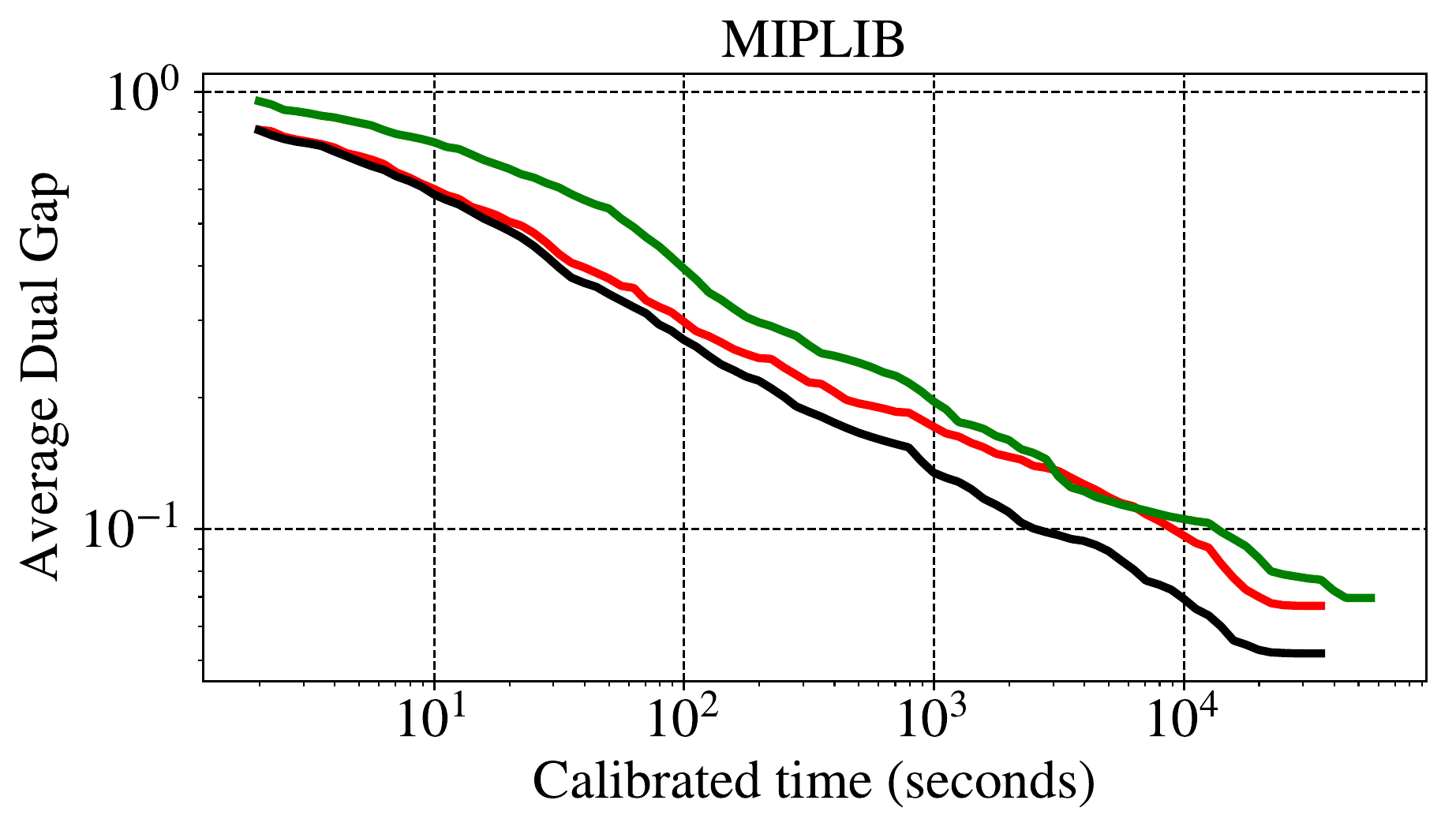} \\

    \end{tabular}
    \caption{Average dual gap with respect to a pre-computed best known primal bound as a function of running time on the test set of benchmark datasets.}
    \label{fig:dual_bound_vs_calibrated_time}
\end{figure}

Figure~\ref{fig:dual_bound_vs_calibrated_time} shows the average dual gap curves for Neural Branching and Tuned SCIP.
We also compare to SCIP's Full Strong branching policy (run with all SCIP parameters set to default, except the branching policy) as a node-efficient baseline. Neural Branching achieves significantly smaller average dual gap for the same running time on four out of six datasets. (Note the log scale of the y-axis.) MIPLIB and Google Production Planning are the only datasets where Neural Branching is not the best, but it is still comparable to Tuned SCIP. The improvements given by Neural Branching can be substantial -- e.g., its average dual gap is about 7$\times$ better on CORLAT, 20$\times$ better on Neural Network Verification, and 2$\times$ better on Electric Grid Optimization, compared to Tuned SCIP at large time limits. On Google Production Packing, Neural Branching is better at intermediate times (e.g., about 1.5$\times$ better at $10^3$ seconds), with Tuned SCIP catching up eventually.

\begin{figure}[H]
    \begin{tabular}{cc}
        \includegraphics[width=0.48\textwidth]{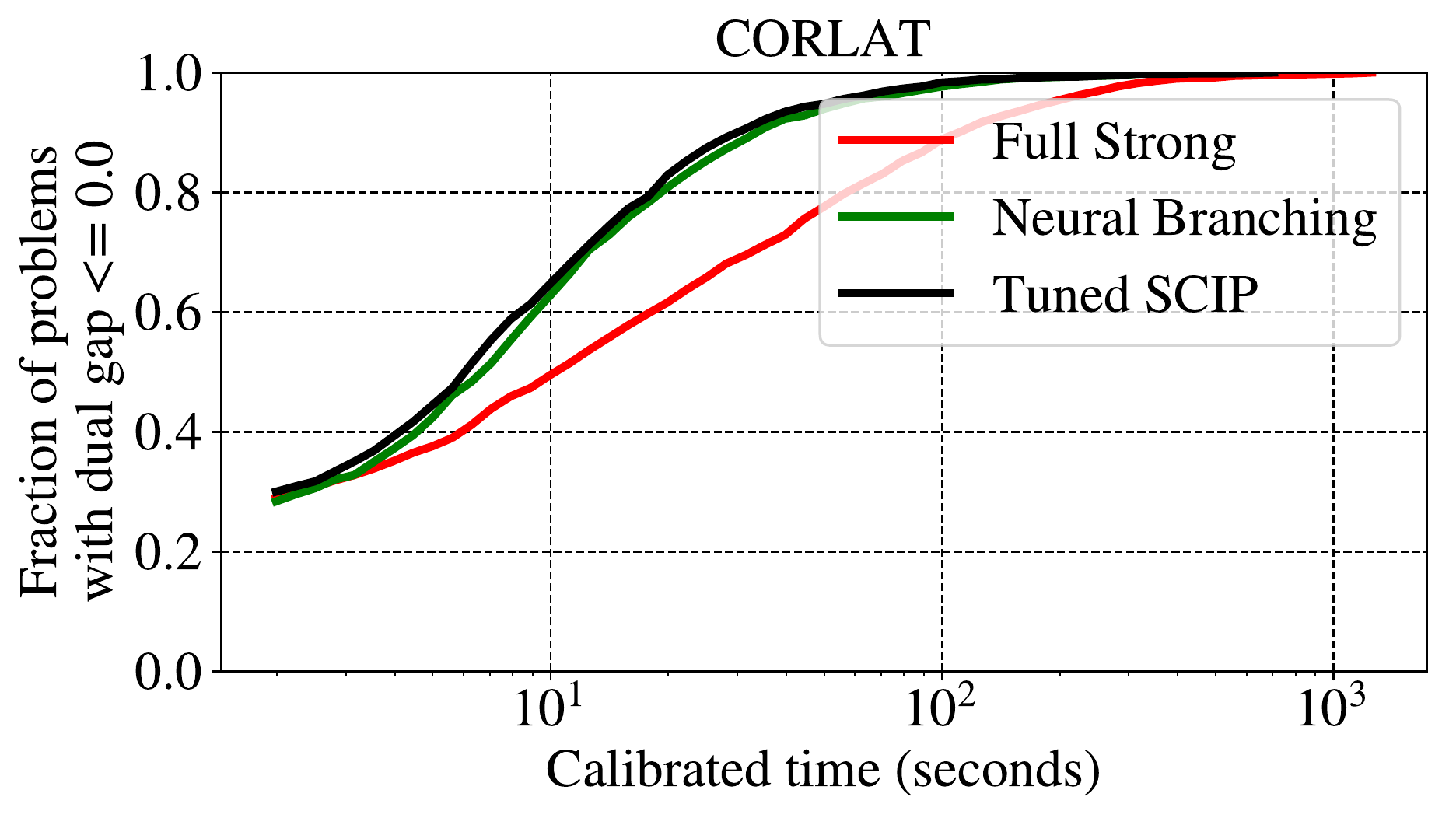} &
        \includegraphics[width=0.48\textwidth]{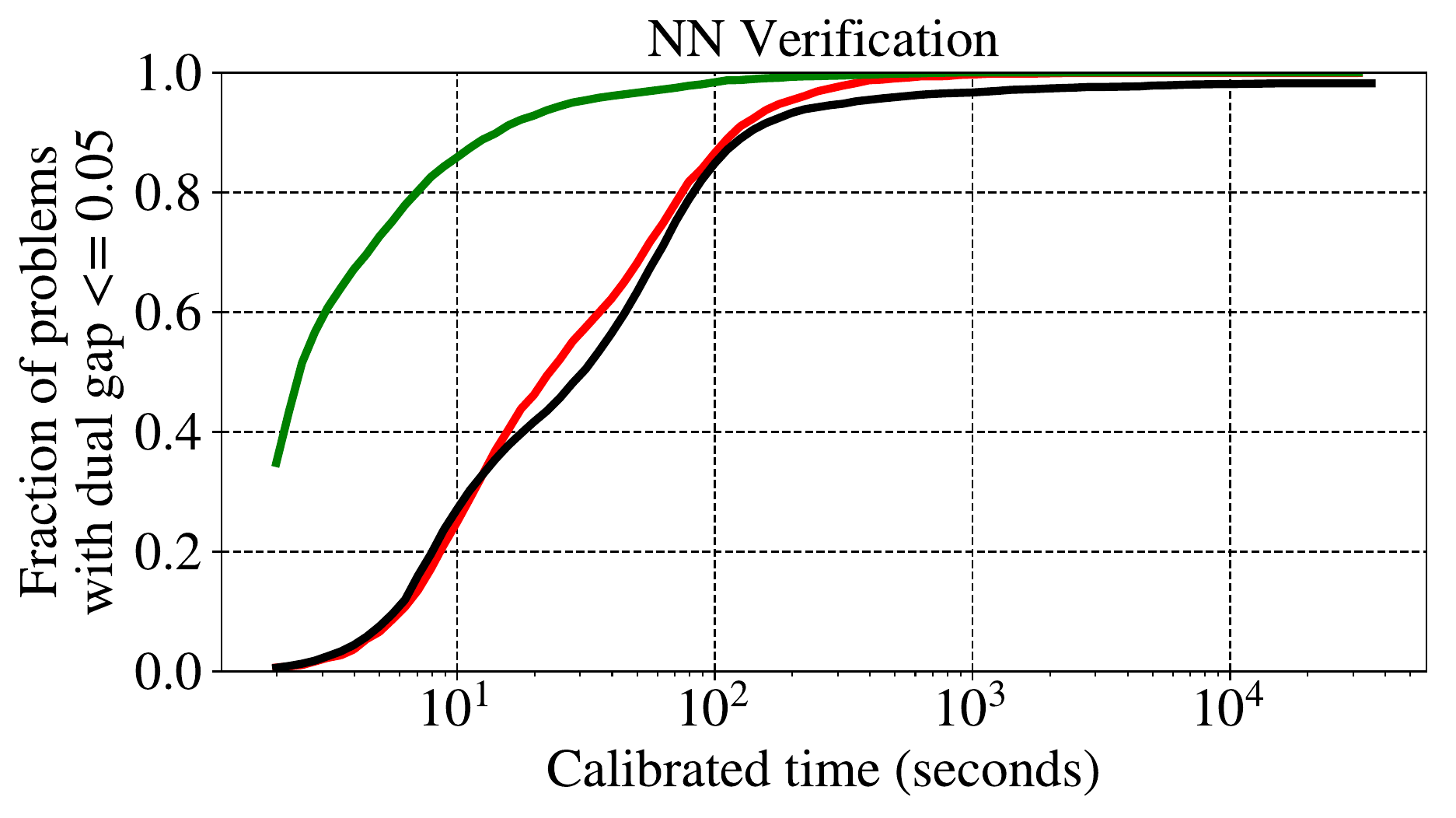} \\
        \includegraphics[width=0.48\textwidth]{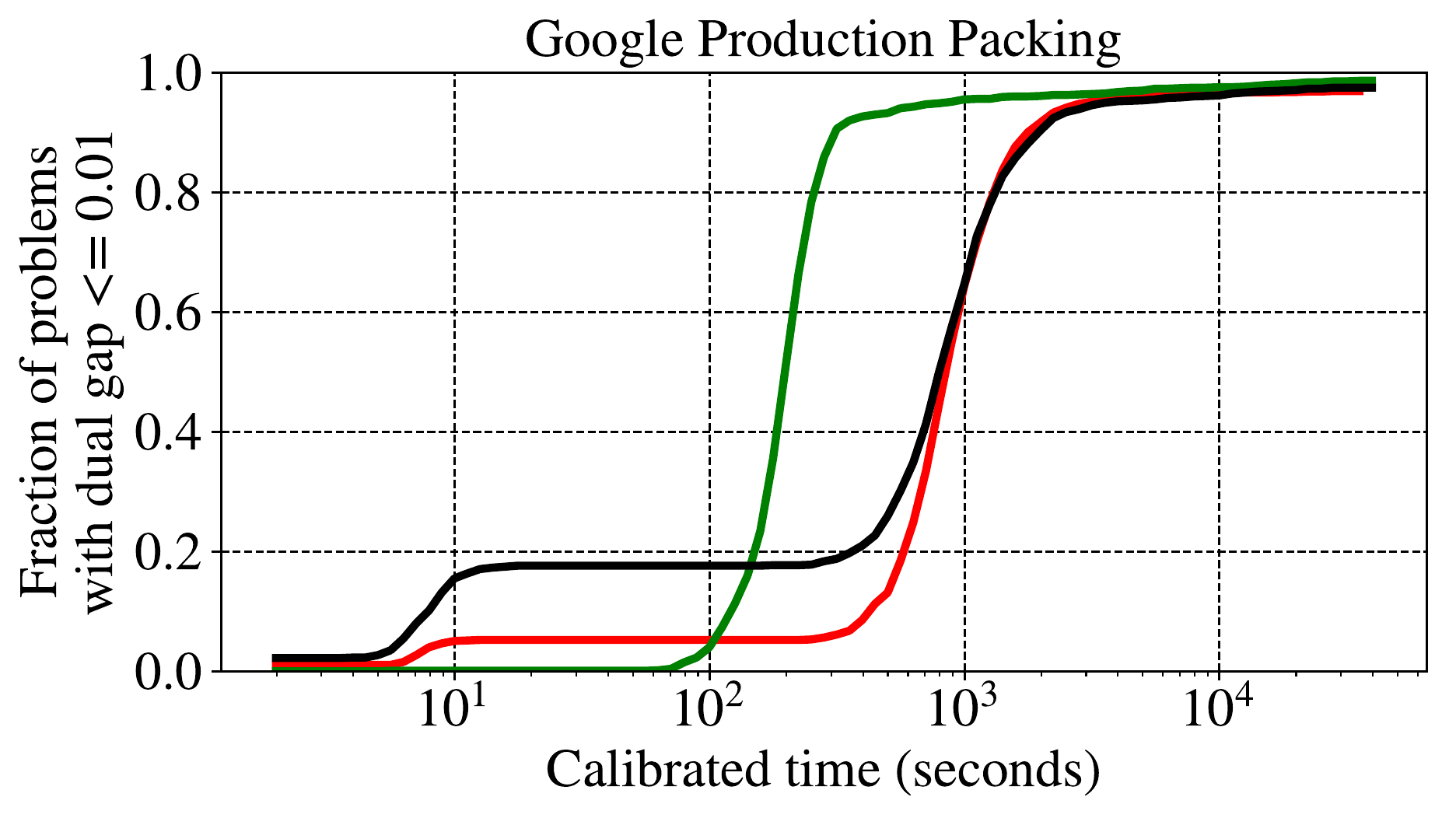} &
        \includegraphics[width=0.48\textwidth]{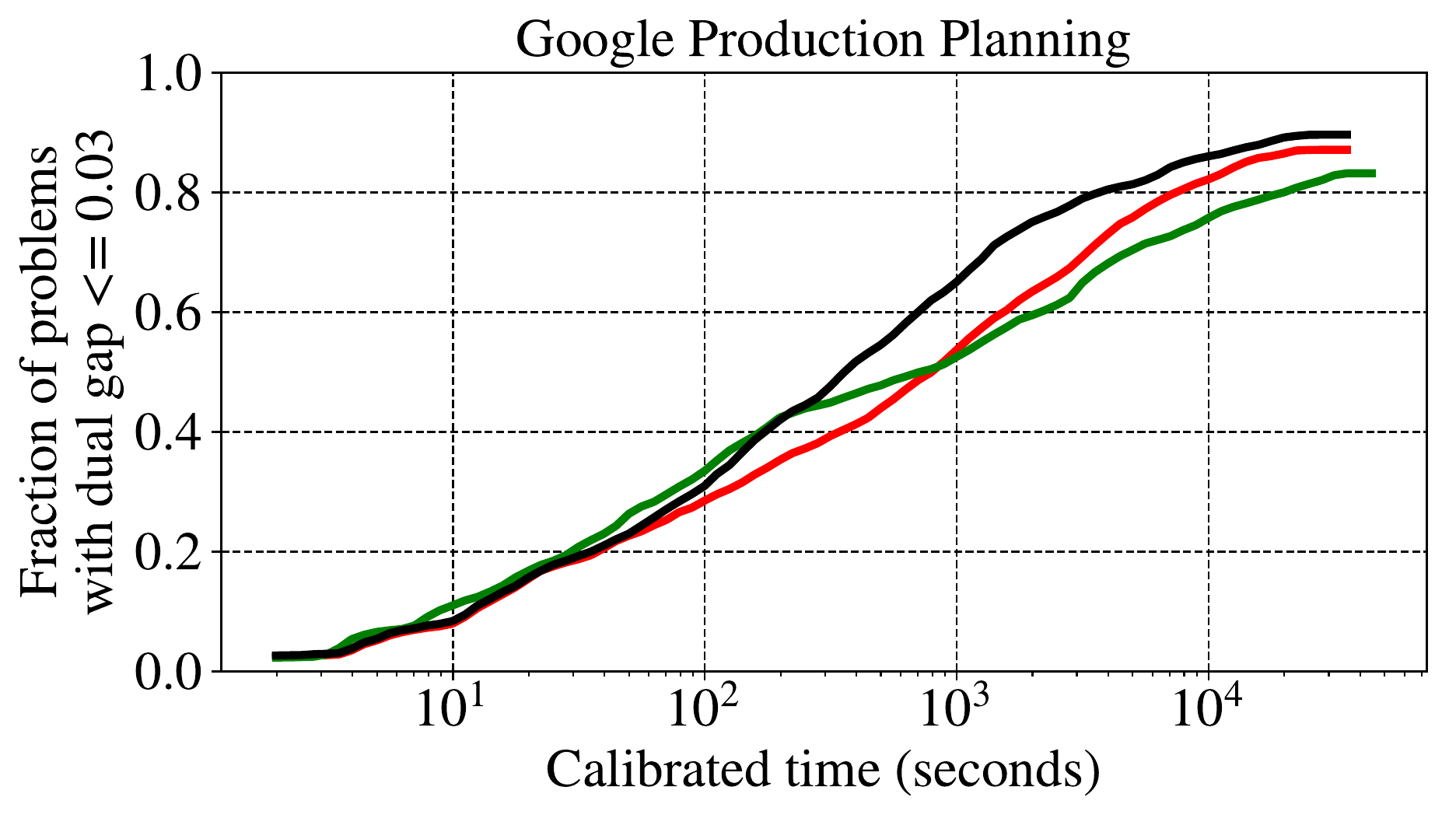} \\
        \includegraphics[width=0.48\textwidth]{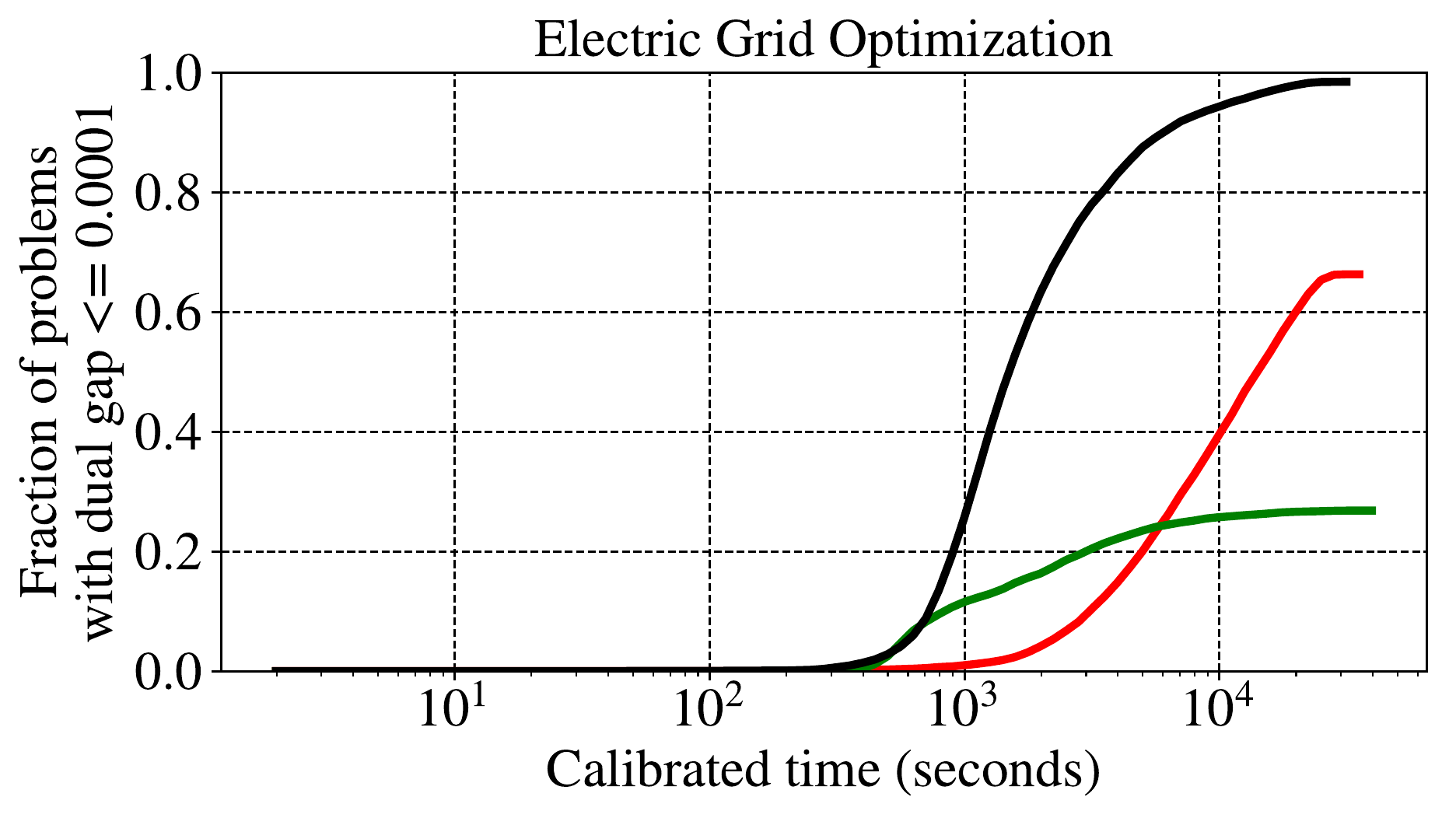} &
        \includegraphics[width=0.48\textwidth]{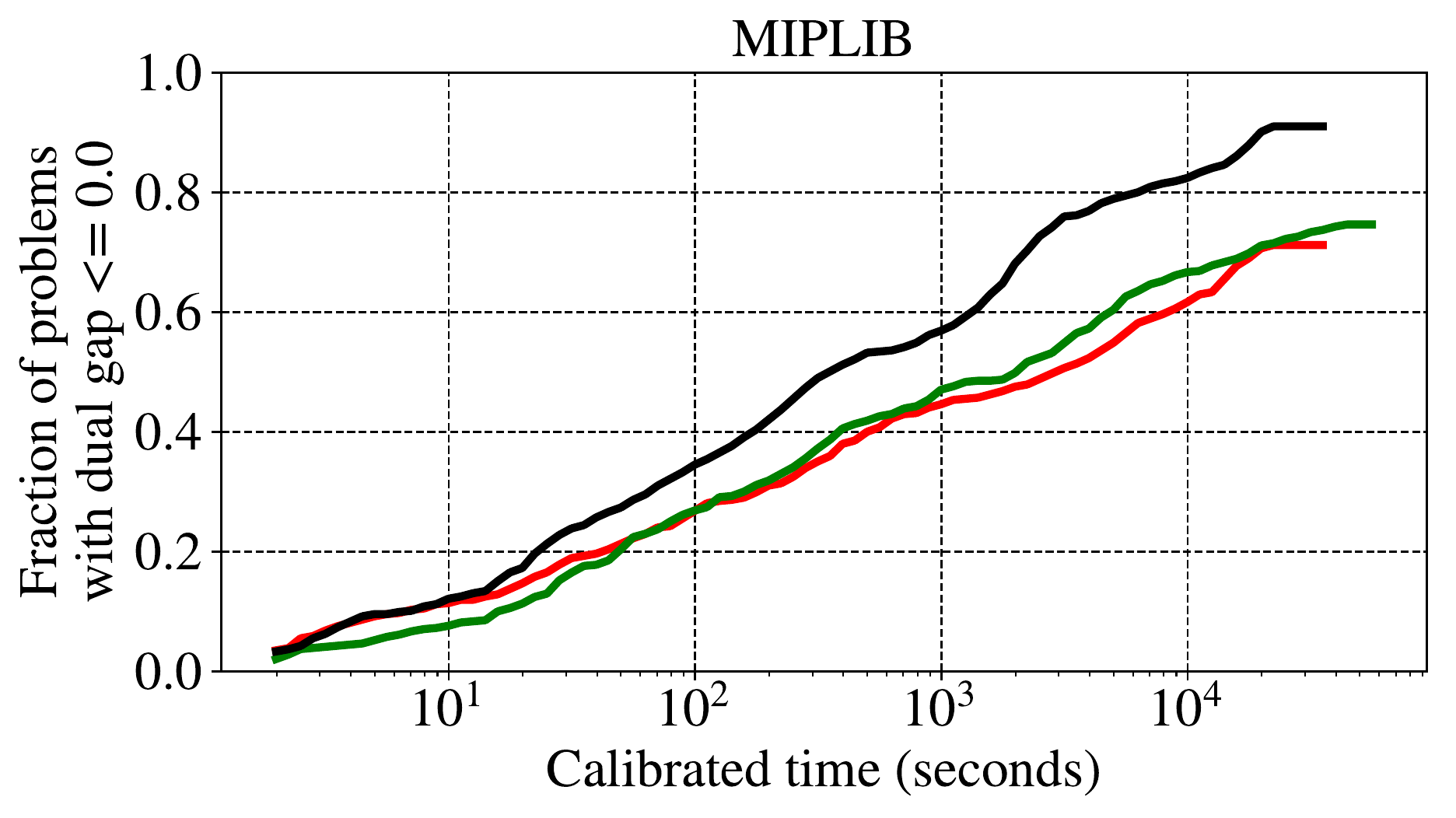} \\
    \end{tabular}
    \caption{Survival plots computed by applying a dataset's target optimality gap to the dual gap of each test set MIP instance.}
    \label{fig:survival_dual_time}
\end{figure}

Figure~\ref{fig:survival_dual_time} shows survival plots computed by applying a dataset's target optimality gap to the dual gap of each test set MIP instance. They confirm the conclusions from figure~\ref{fig:dual_bound_vs_calibrated_time} -- Neural Branching solves a higher fraction of test instances consistently across all time limits on all datasets except Google Production Planning and MIPLIB.

%% file: joint_eval_results.tex
\section{Joint Evaluation}
\label{sec:joint_eval}
We now combine Neural Branching and Neural Diving into a single solver. As we shall see, this results in significant speedups over Tuned SCIP. We use the learned heuristics by integrating them into SCIP via the interface provided by PySCIPOpt \citep{pyscipopt}. We consider all four possible ways of combining Neural Branching and Neural Diving with SCIP: 1) \emph{Tuned SCIP} alone, 2) \emph{Neural Branching + Neural Diving (Sequential)} uses both neural heuristics, 3) \emph{Neural Branching} uses only the learned branching policy, and 4) \emph{Tuned SCIP + Neural Diving (Sequential)} uses only the sequential version of Neural Diving. When Neural Diving (Sequential) is used as a primal heuristic, we disable all built-in primal heuristics of SCIP.

To make a comparison based on Neural Diving fairer in terms of compute resources, we only consider its sequential version here. For a single MIP solve, Neural Diving (Sequential) and Neural Branching each use a CPU core and a GPU. So \emph{Neural Branching + Neural Diving (Sequential)} uses two cores and two GPUs. We also give Tuned SCIP two cores, which it uses by running two independent SCIP solves on the input MIP with different random seeds. The primal, dual, and primal-dual gaps for Tuned SCIP at a given time are then computed by selecting the best primal and dual bounds across the two runs at that time. Since SCIP does not use GPUs, we do not control for that resource when comparing Tuned SCIP and Neural Solvers.

For each expert we solve each MIP in the test dataset five times, with five different seeds for SCIP. We present the primal-dual gap $\gamma_{pd}(t)$ (see section \ref{sec:evaluation}) averaged over each dataset's test MIPs as a function of calibrated time. We also present survival plots that show the fraction of test MIPs solved within the required dataset-specific optimality gap as a function of calibrated time.

\subsection{Results}
\begin{figure}[t]
    \begin{tabular}{cc}

        \includegraphics[width=0.48\textwidth]{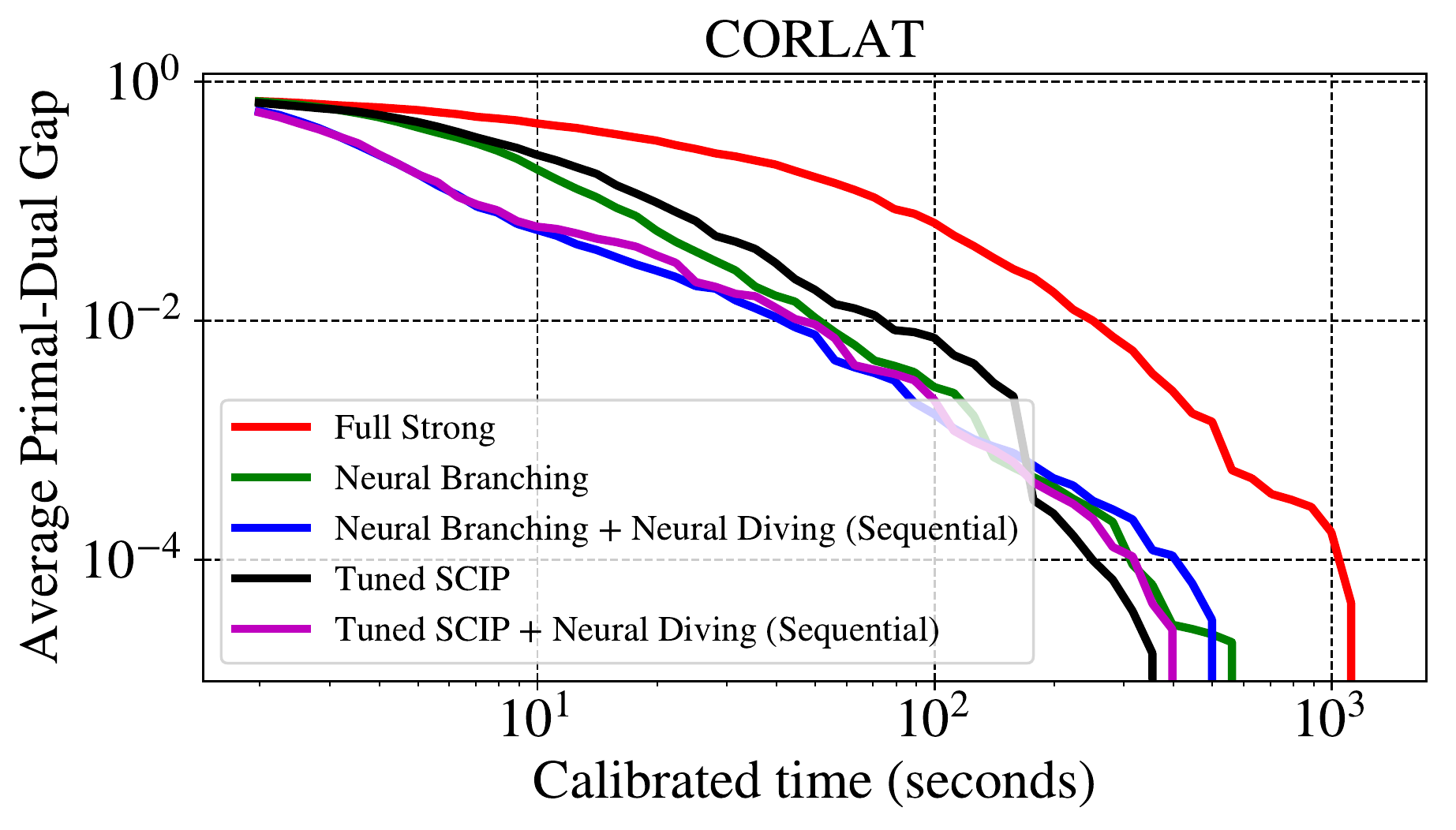} &
        \includegraphics[width=0.48\textwidth]{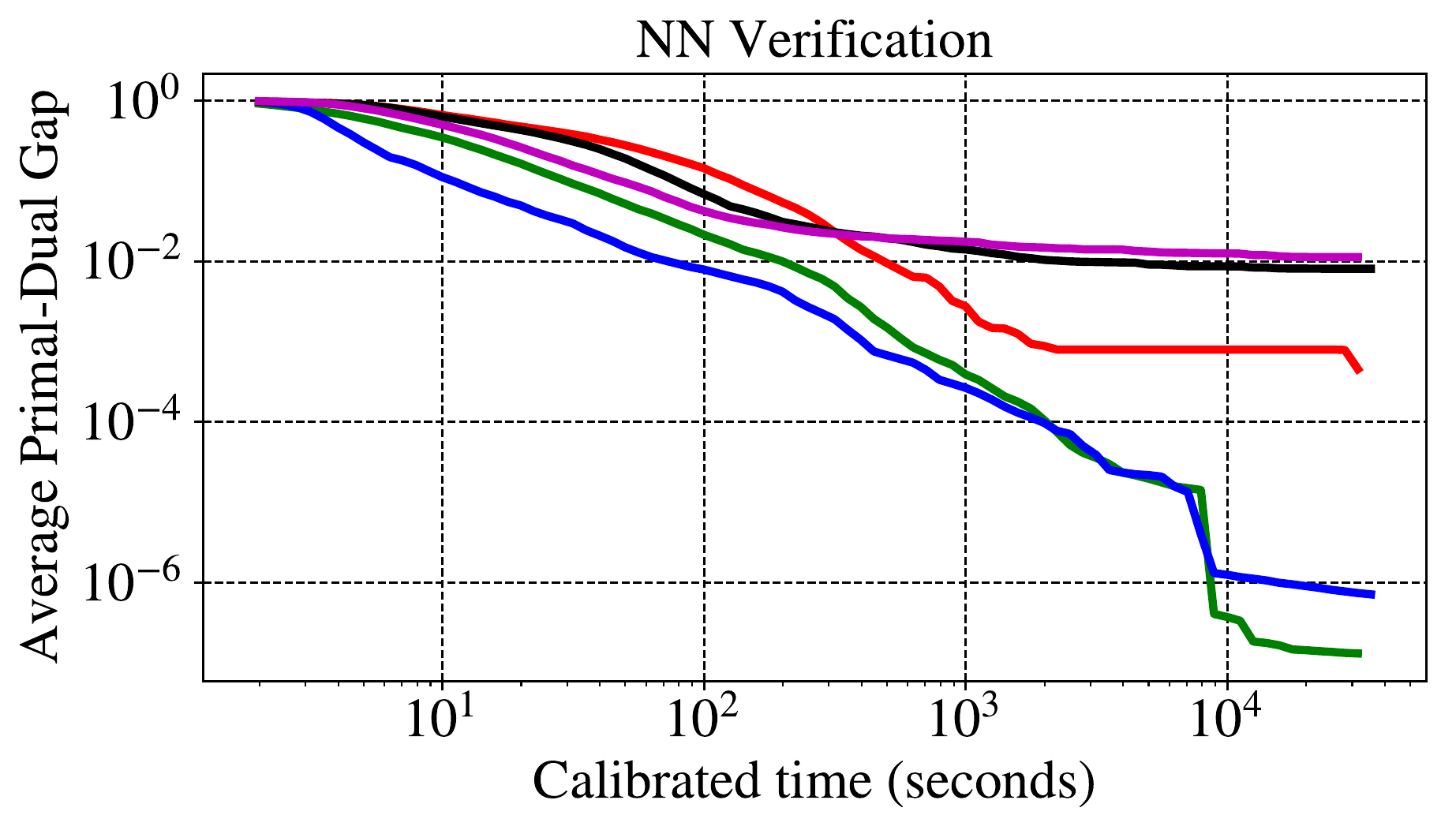} \\
        \includegraphics[width=0.48\textwidth]{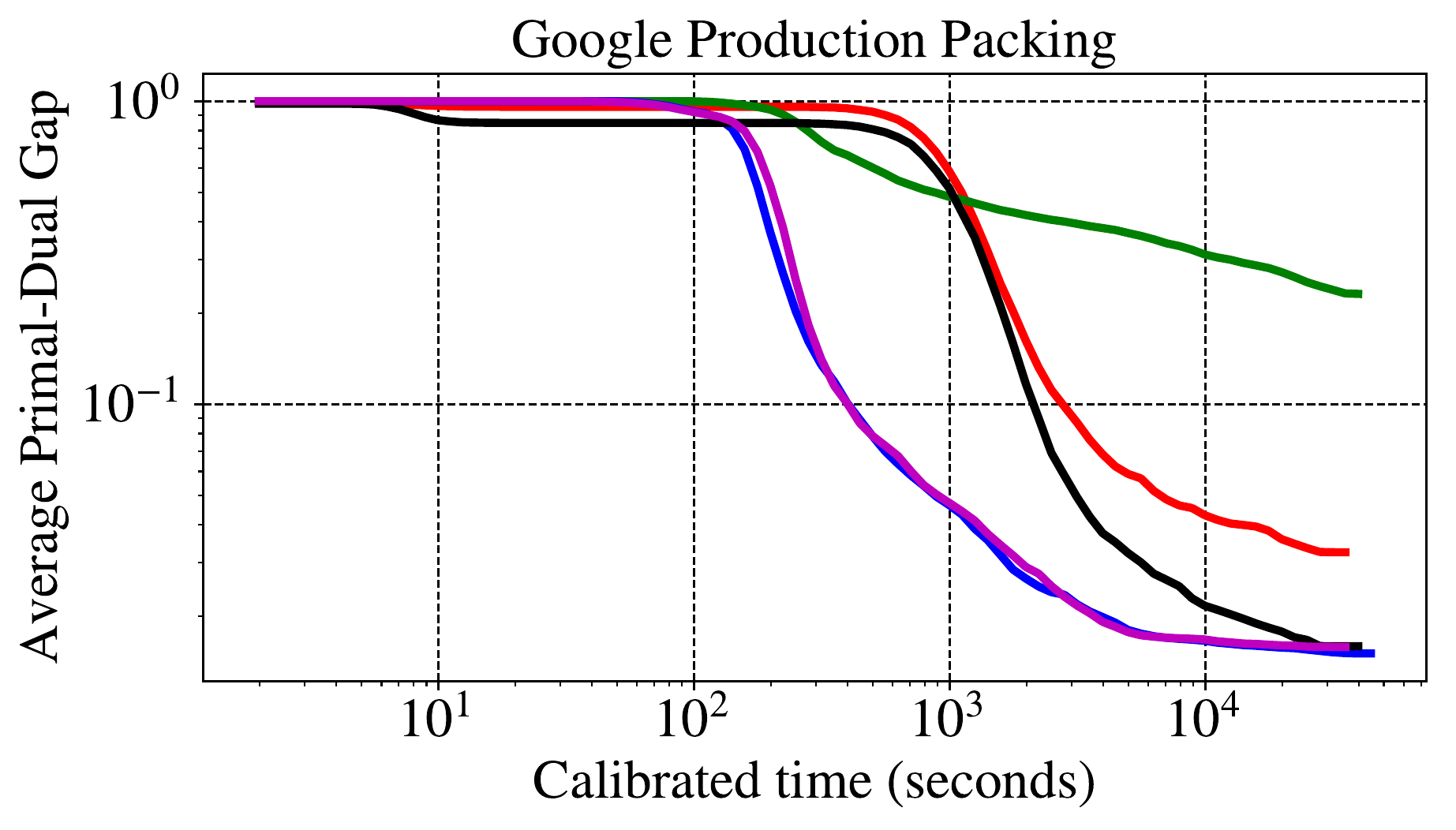} &
        \includegraphics[width=0.48\textwidth]{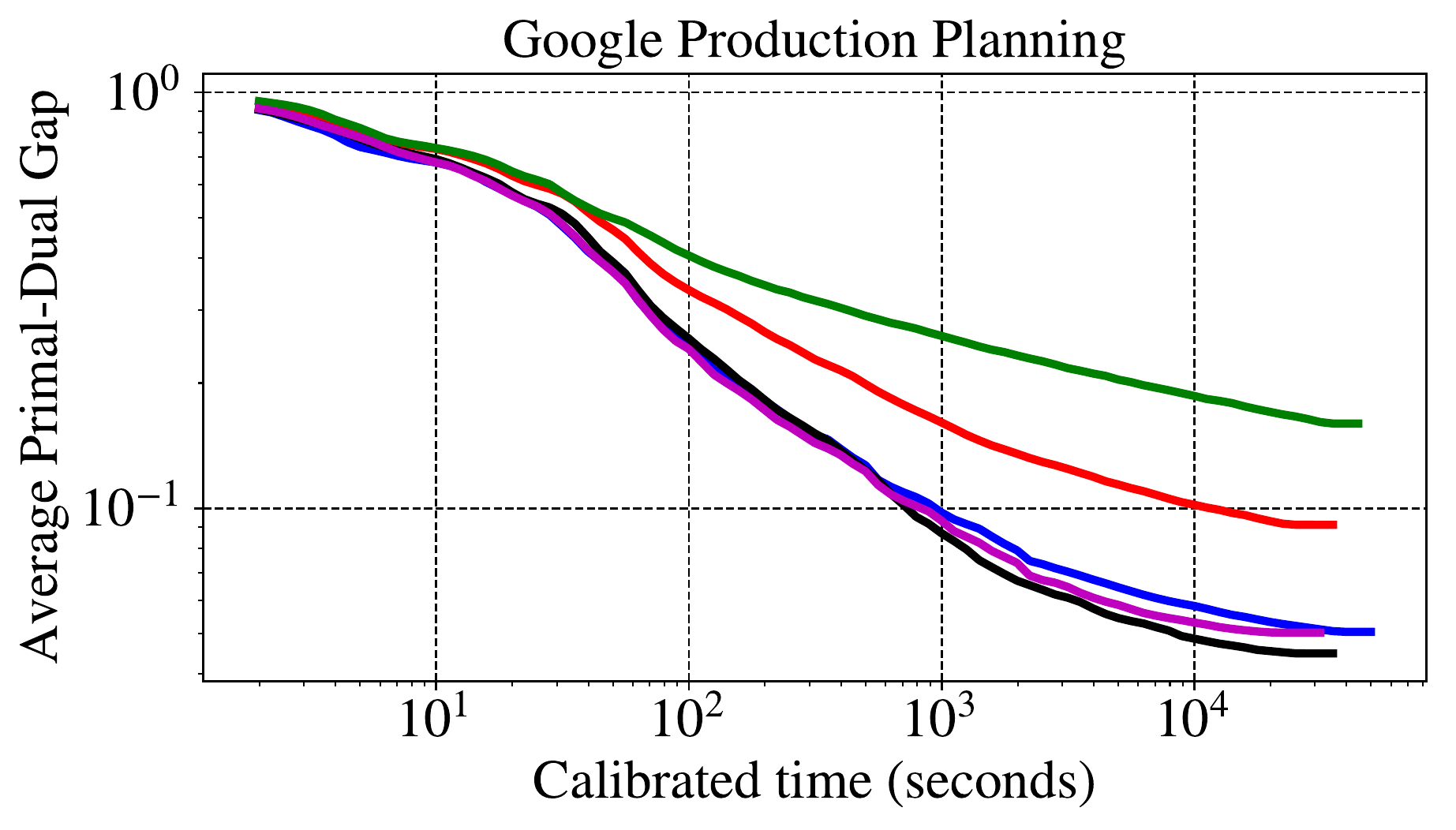} \\
        \includegraphics[width=0.48\textwidth]{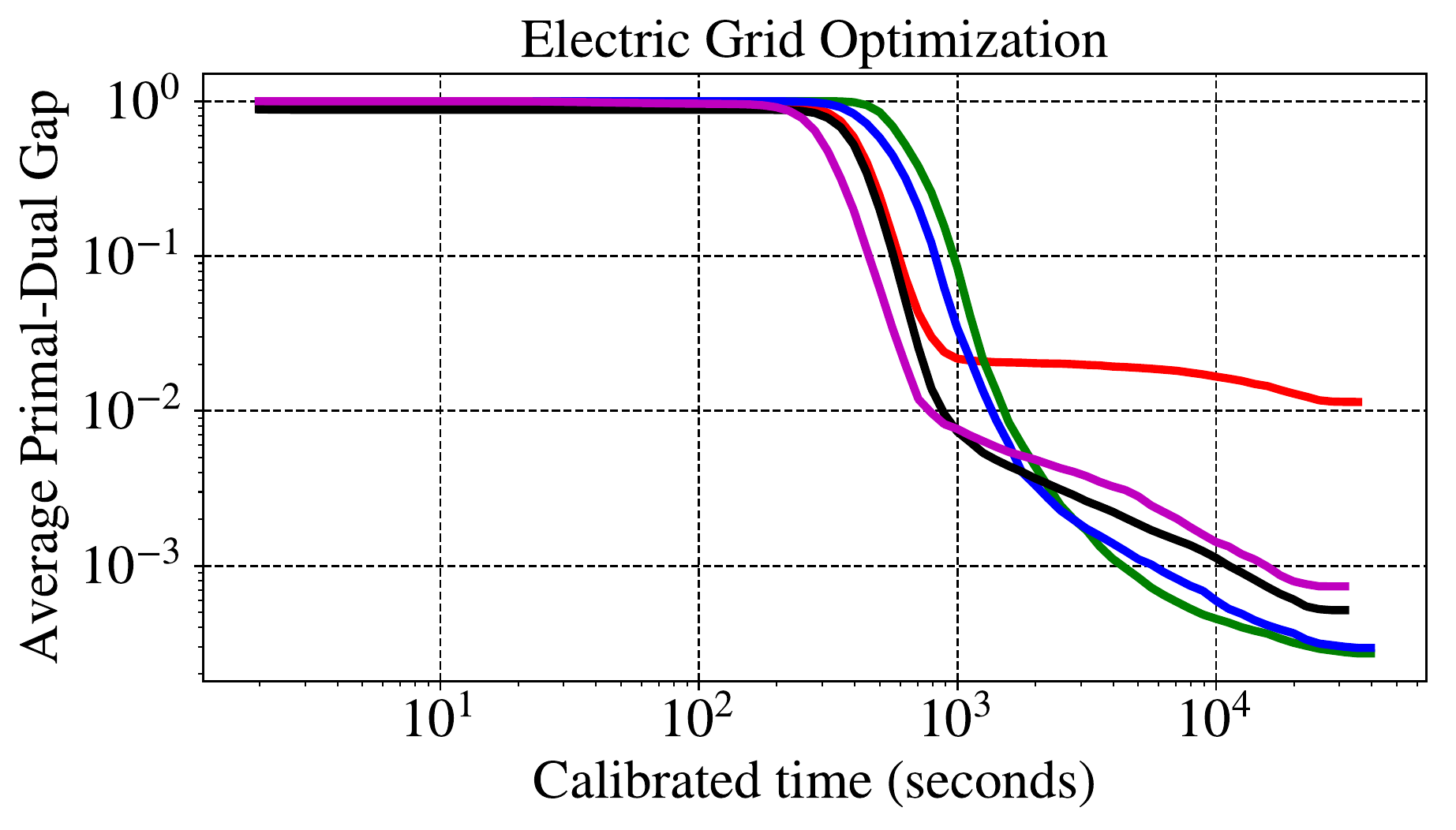} &
        \includegraphics[width=0.48\textwidth]{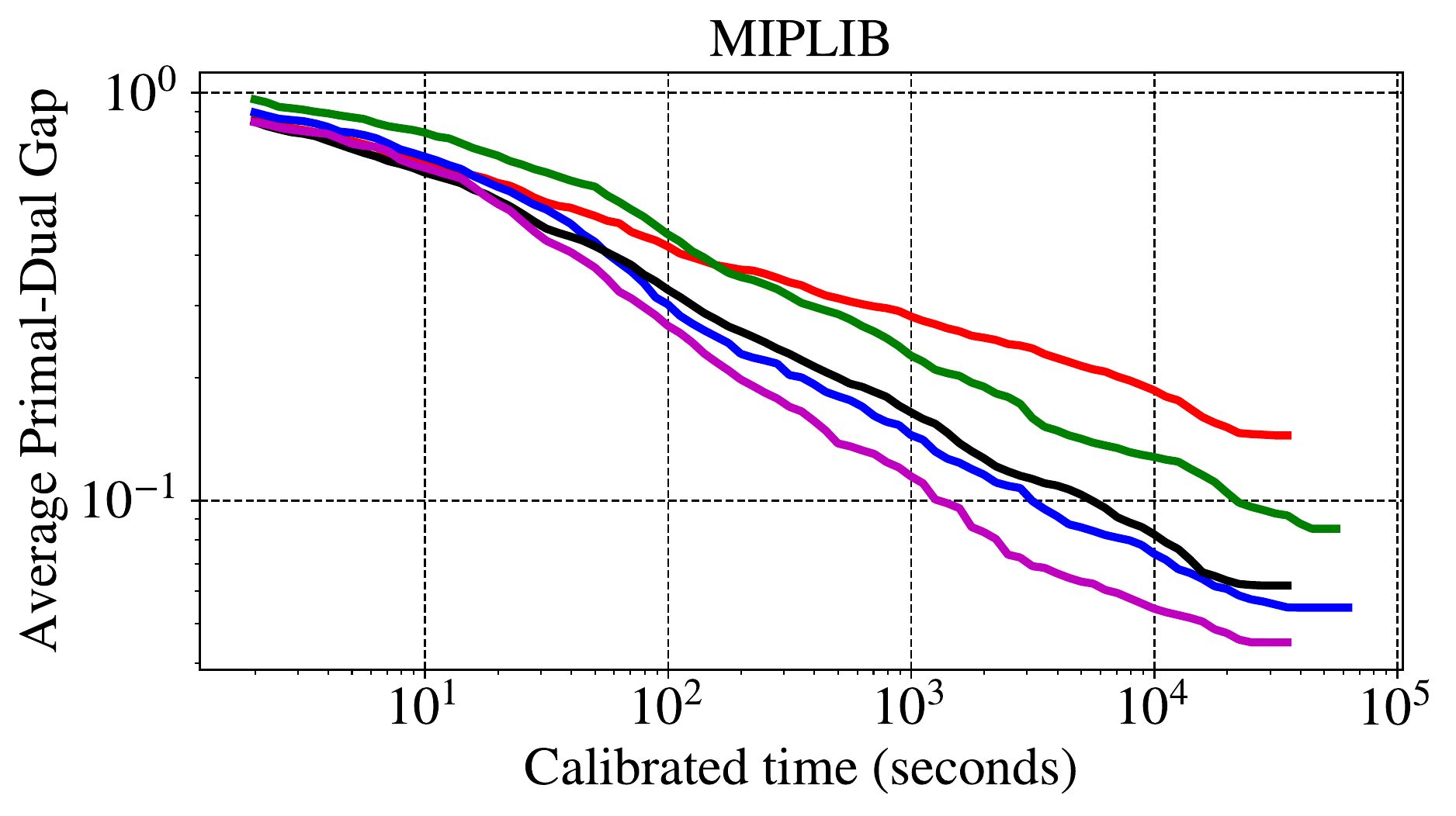} \\

    \end{tabular}
    \caption{ Primal-dual gap achieved by various combinations of primal heuristics and variable selection policies as a function of running time on the benchmark datasets. The gap is computed for each MIP in the test set of each dataset and then aggregated using the average.}
    \label{fig:rel_gap_vs_calibrated_time}
\end{figure}

Figure~\ref{fig:rel_gap_vs_calibrated_time} shows the average primal-dual gap curves as a function of running time. A Neural Solver significantly outperforms Tuned SCIP on four of the datasets. (Note the log scale of the y-axis.) The average primal-dual gap is more than five orders of magnitude better than Tuned SCIP on Neural Network Verification at large time limits. On Google Production Packing, Neural Branching + Neural Diving (Sequential) achieves low gaps much more quickly, reaching a $0.1$ gap in more than $5\times$ less time, but Tuned SCIP eventually catches up. On Electric Grid Optimization, it achieves a $2\times$ lower gap at higher running times. On MIPLIB, Tuned SCIP + Neural Diving (Sequential) achieves a $1.5\times$ better gap at higher running times.

\begin{figure}[t]
    \begin{tabular}{cc}

        \includegraphics[width=0.48\textwidth]{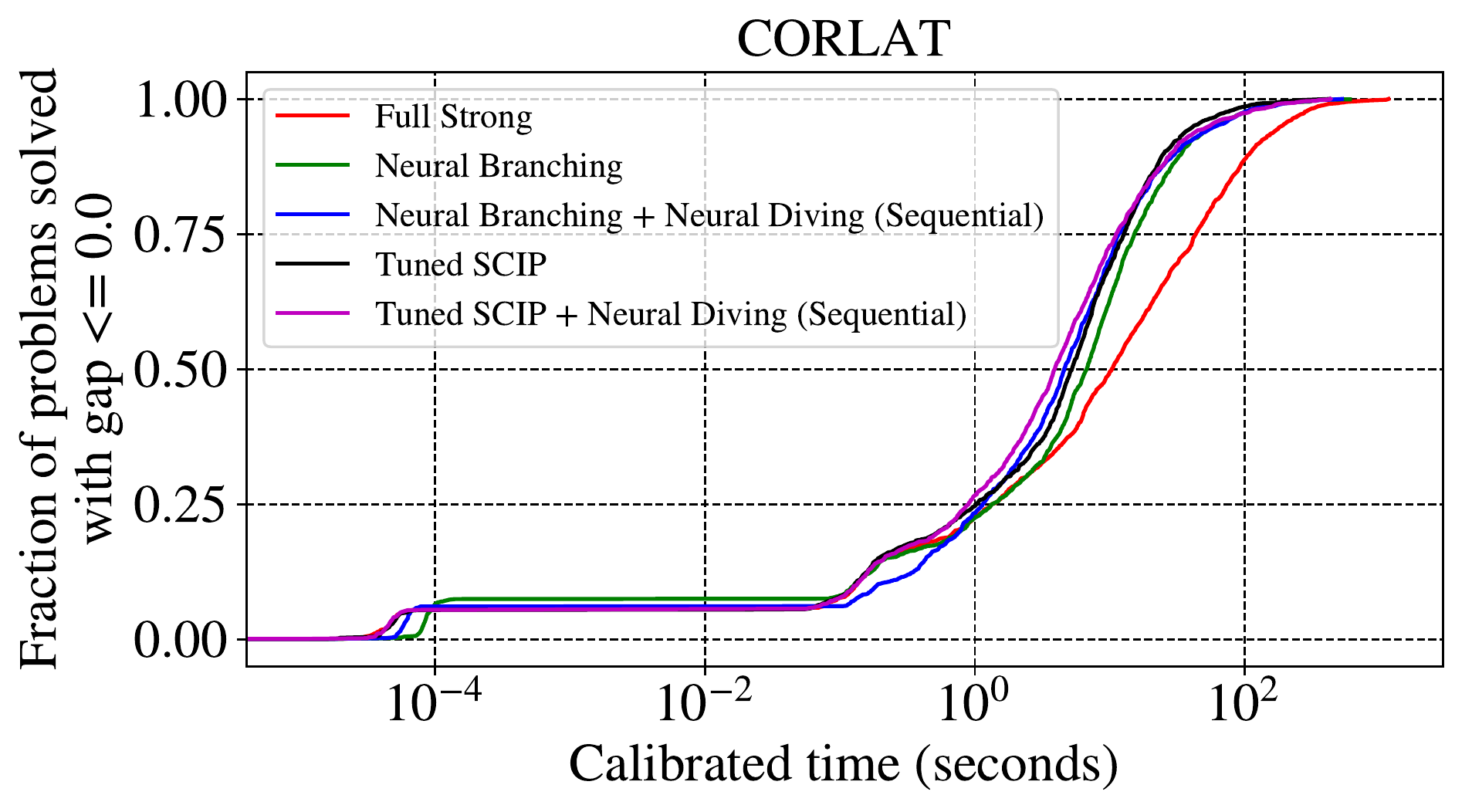} &
        \includegraphics[width=0.48\textwidth]{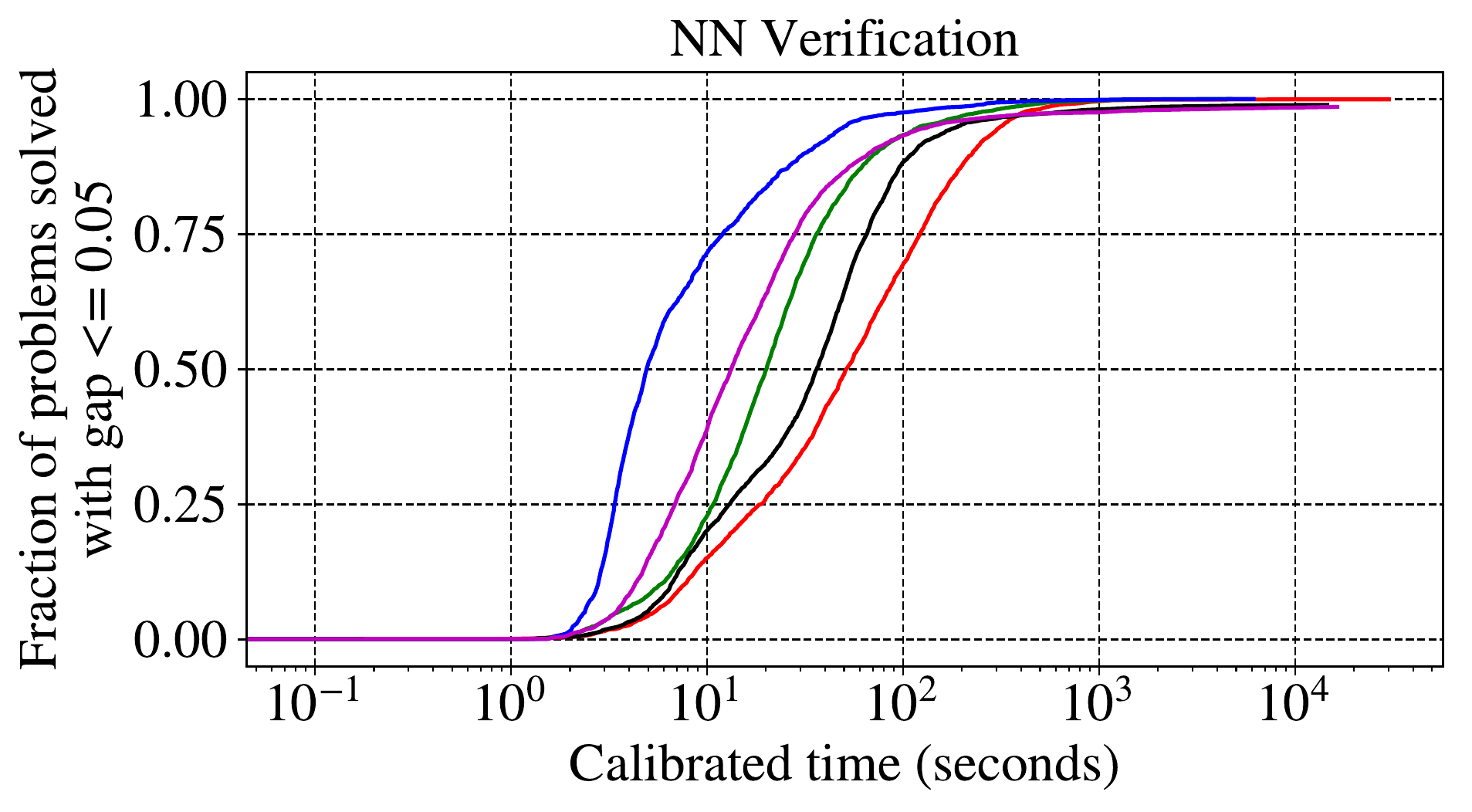} \\
        \includegraphics[width=0.48\textwidth]{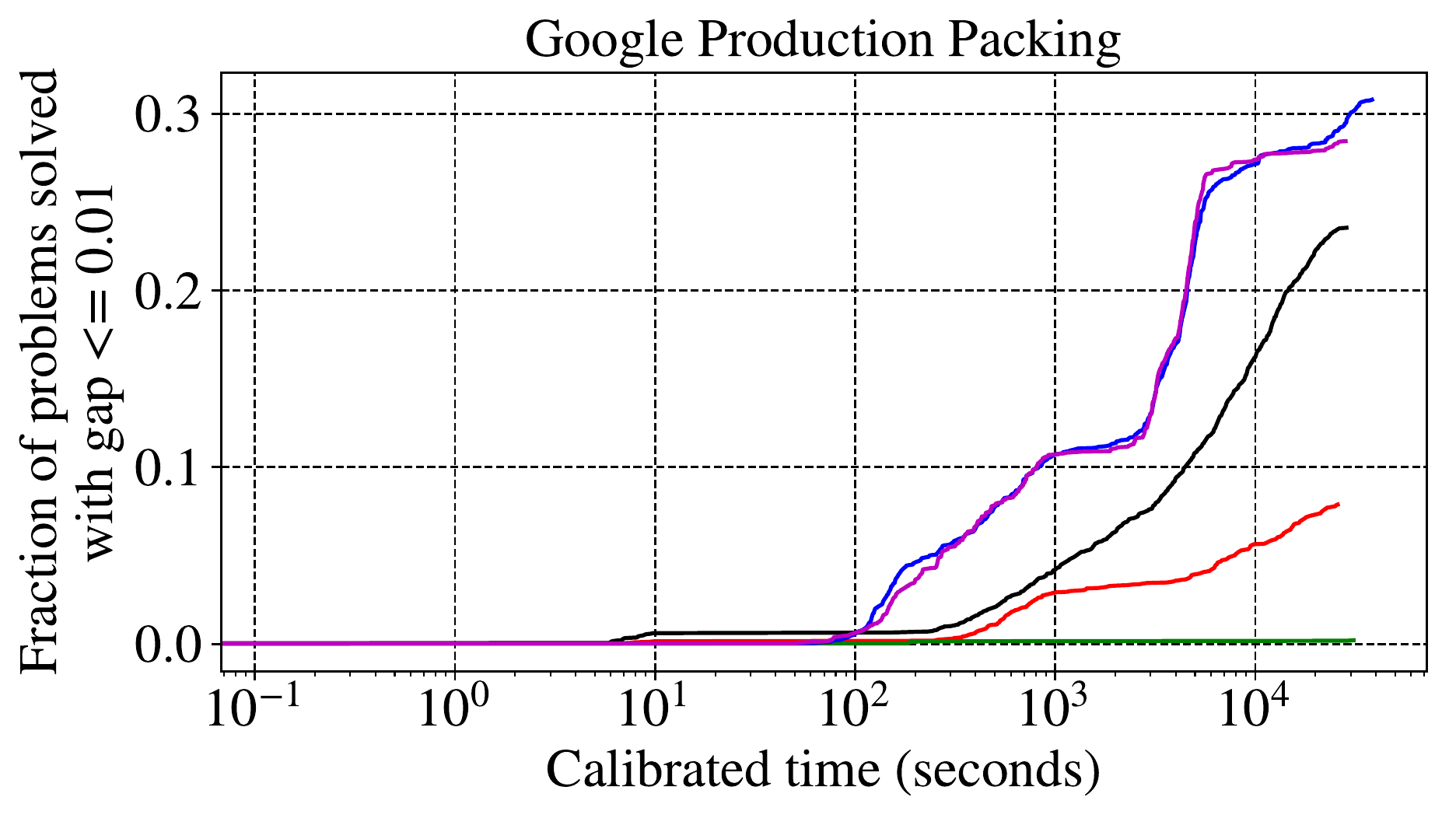} &
        \includegraphics[width=0.48\textwidth]{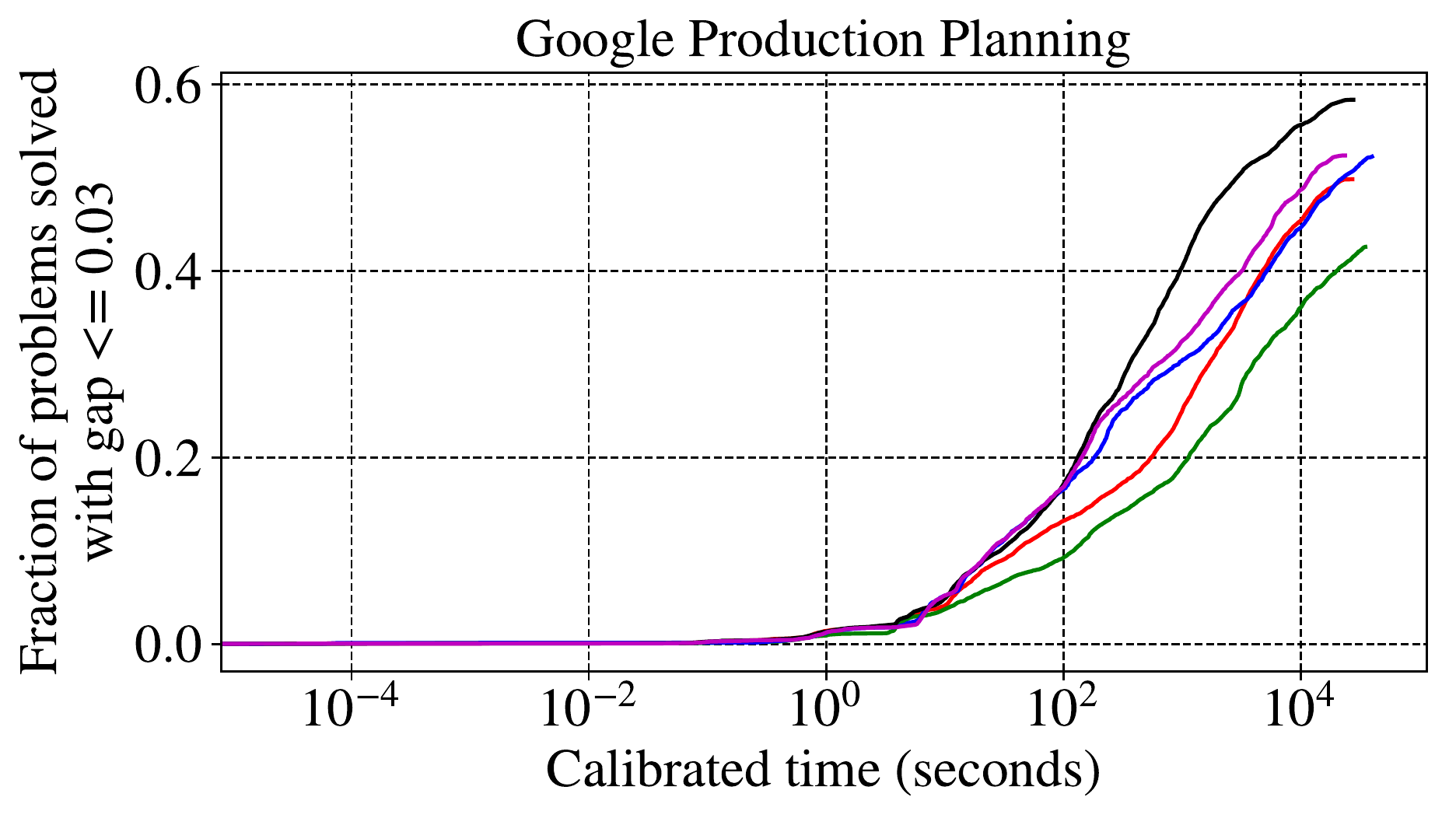} \\
        \includegraphics[width=0.48\textwidth]{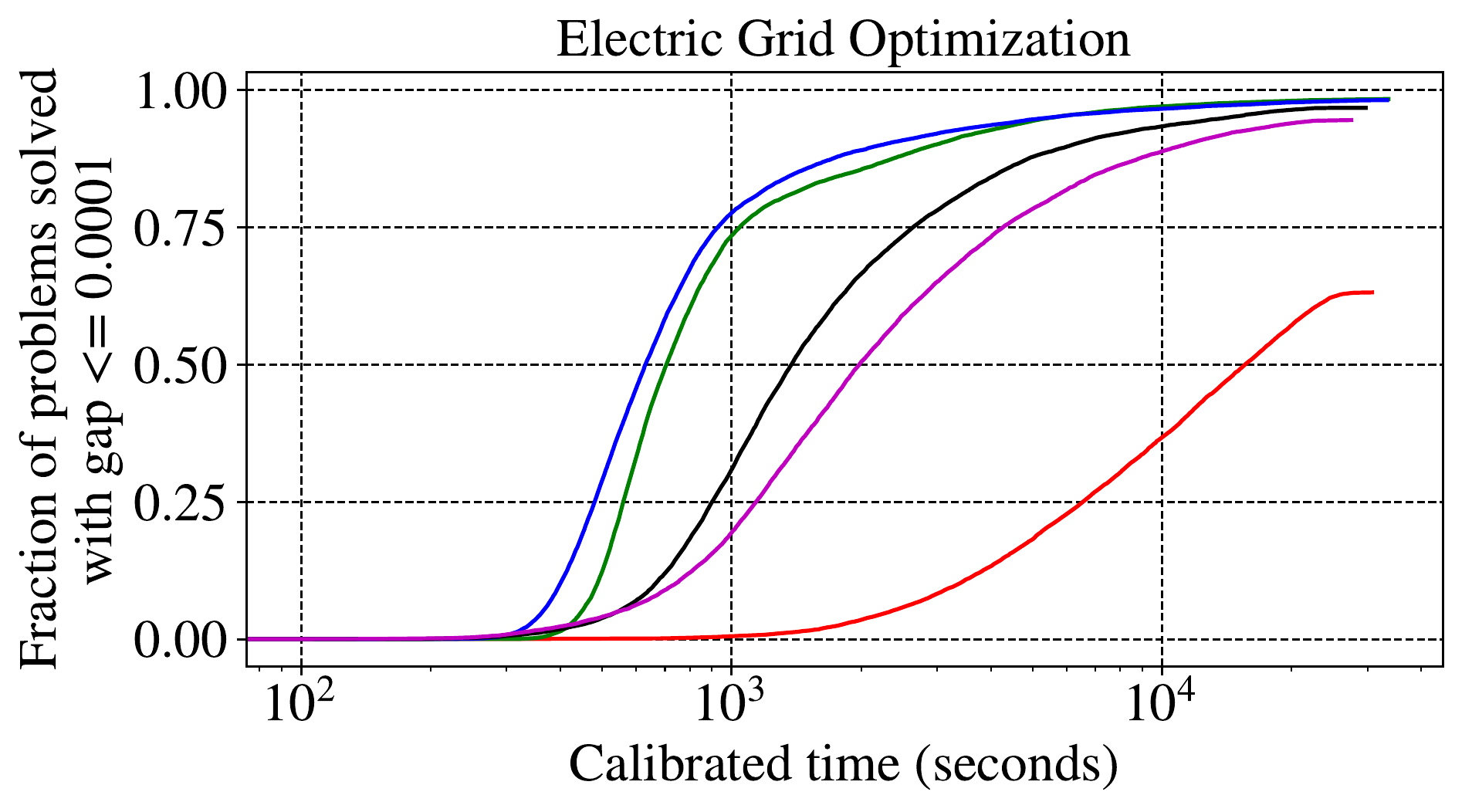} &
        \includegraphics[width=0.48\textwidth]{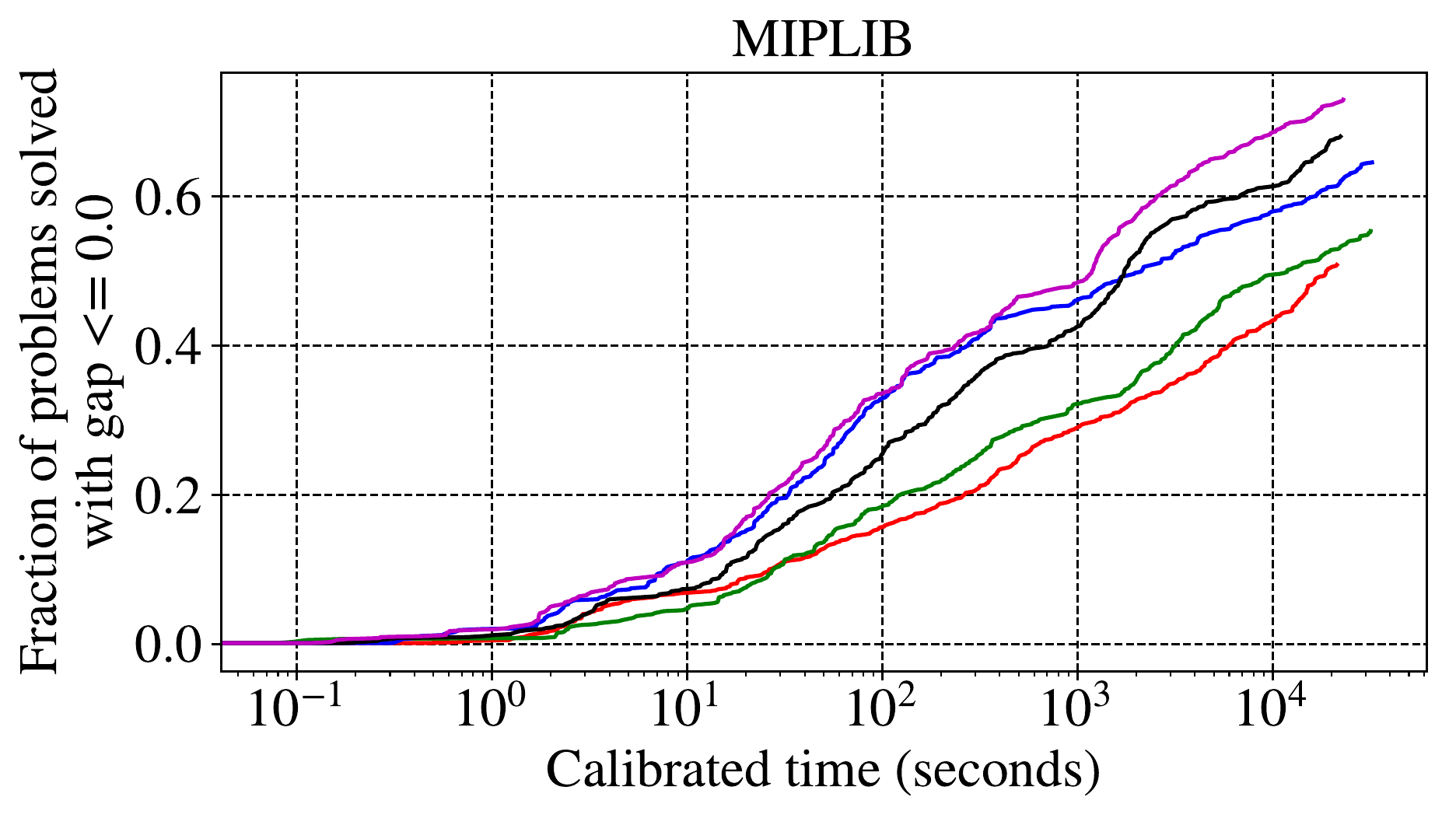} \\
    \end{tabular}
    \caption{Survival plots with respect to calibrated running time for various combinations of primal heuristics and variable selection policies as a function of running time on the benchmark datasets.}
    \label{fig:survival_plots_calibrated_time}
\end{figure}

Figure~\ref{fig:survival_plots_calibrated_time} shows the survival plots as a function of running time. They further confirm the observations from figure~\ref{fig:rel_gap_vs_calibrated_time}, with a Neural Solver solving a higher fraction of test set problems at a given time limit than Tuned SCIP on four of the datasets. Together, the results show that learning can be used to significantly improve a strong solver such as SCIP. The improvements are seen across different problem scales and applications. Surprisingly, learning improves performance even on MIPLIB, a heterogeneous dataset combining many applications that a priori would not be considered similar enough to have sufficient shared structure for learning to exploit. That the improvements are observed across diverse settings confirms the broad usefulness of learning.

A further comparison with respect to the \emph{Penalized Average Runtime} metric can be found in Appendix section~\ref{subsec:par10}. The overall conclusions remain the same.

%% file: related_work.tex
\section{Related Work}
\label{sec:related_work}

\noindent\textbf{Learning and MIP Primal Heuristics:} Learning has previously been used to switch among an ensemble of existing primal heuristics during a branch-and-bound run. \cite{khalil2017primal} learn a binary classifier offline to predict whether applying a primal heuristic at a branch-and-bound node will succeed in finding an incumbent. \cite{Hendel2018alns} formulate a multi-armed bandit approach to learn a switching policy online. Neural Diving constructs a new primal heuristic specifically for a given application. This can be particularly powerful if none of the existing heuristics are well-suited for the application.

Several works consider learning neighborhood search policies \citep{addanki2020nlns, hottung2019neurallns, song2020generallns}. Given a feasible point, a learned policy is used to modify the assignment for a subset of variables such that the resulting feasible assignment has a lower objective. In these works the learned policy is used to either select the variables to be modified, or to assign new values to an already selected subset of variables. They require a feasible assignment as input, while Neural Diving does not. 

Learning to predict variable values in a MIP and combining the predictions with a solver has been studied before. \cite{ding2019accelerating} trains a model that takes a MIP as input and predicts values for a subset of its binary variables. This is then used to define an additional constraint that any solution be within a pre-specified Hamming distance from the predicted values, which can significantly reduce the search space and speed up a MIP solver. Training labels are generated by iteratively improving a feasible assignment and selecting only the subset of binary variables whose values are stable across iterations as prediction targets. Results are presented on several synthetic datasets. \cite{xavier2020mluc} learns to warm start a MIP solver specifically for the electric grid \emph{unit commitment} problem. They use a $k$-Nearest Neighbor binary classifier to predict values for a subset of binary variables in the MIP. Assuming a fixed grid topology and a fixed set of power plants allows them to use a MIP representation with constant dimensionality and a pre-defined distance metric to identify neighbors. Our work differs significantly from both of these: 1) the generative modelling formulation provides a principled approach to capturing the uncertainty in variable values and exploiting it to generate several sub-problems to find high-quality feasible assignments, 2) it handles general MIPs (varying sizes, non-binary integers), and 3) we present results on diverse real-world applications.

\noindent\textbf{Learning and MIP Branching Heuristics:} Several works have used imitation learning to train a branching policy to imitate (a variant of) Full Strong Branching. \cite{Khalil2016LearningToBranch} learns to imitate Full Strong Branching for a single instance, using data collected as it is being solved. They apply a learning-to-rank formulation and handcrafted features. Similarly, \cite{Alvarez2017MLApproxStrongBranching} learns a regression model to predict the scores computed by Full Strong Branching for candidate variables, again using handcrafted features. Neither of these works report current state-of-the-art results for the learned policies. Neural Branching builds on the work by \cite{gasse2019exact}, which trains a Graph Convolutional Network to imitate decisions by Full Strong Branching. The main differences are 1) our use of the ADMM-based expert to scale up Full Strong Branching to large instances, and 2) extensive evaluation on real-world applications and significant performance improvements over SCIP.

More recently, \cite{zarpellon2020parameterizing} show that using a representation of the branch-and-bound tree (rather than that of just a single node) can improve the generalization of a learned branching policy. \cite{gupta2020hybrid} obtain computational improvements similar to \cite{gasse2019exact} in the setting of CPU-restricted machines by using an expensive graph network only at the root node of the branch-and-bound tree and a cheaper MLP elsewhere.

\noindent\textbf{Other Learning Approaches for MIP Solvers:} Learning has been used to automatically tune parameters for MIP solvers. \emph{Algorithm Configuration} \citep{ansotegui2009gga, hutter2009paramils, hutter2011smac, ansotegui2015ggapp} aims to find parameter values that improve aggregate performance on a dataset of instances. For example, \cite{hutter2011smac} proposes Sequential Model-based Algorithm Configuration (SMAC), which iterates over three main steps: 1) evaluate the performance of selected solver parameter values on a set of instances, 2) learn a model to predict the performance of given parameter values from the evaluations done so far, and 3) use the model learned so far to select new candidate parameter values to evaluate. Results on CPLEX show that SMAC finds parameter values that are better than the default and the values found by CPLEX's tuning tool. \emph{Algorithm Selection} \citep{kotthoff2016AlgoSelectSurvey} aims to improve solver parameter values for a specific instance. For example, \cite{hutter2014epm} shows that it is possible to accurately predict the running time of MIP solvers as a function of the solver parameters and the instance being solved. Such a predictor can be used to improve the solver parameters for a specific instance. Note that learning techniques for Algorithm Configuration and Algorithm Selection are complementary to and can be combined with the learning of application-specific heuristics.

\noindent\textbf{Learning for Combinatorial Optimization:} Our work is an instance of the broader topic of learning to solve combinatorial optimization problems. Some of the earliest works in this area are \cite{zhang1995RLforJSS, moll1999routing, boyan97stage}. More recently, deep learning has been applied to the Travelling Salesman Problem \citep{vinyals2015pointernets, bello2016neural}, Vehicle Routing \citep{kool2018attention, nazari2018RLforVRP}, Boolean Satisfiability \citep{selsam2019NeuroSAT}, and general graph-structured combinatorial optimization problems \citep{khalil2017learning, li2018combinatorial}. A survey of the topic is available by \cite{bengio2018survey}.

%% file: conclusions.tex
\section{Conclusion}
This work has demonstrated the long-held promise of machine learning to significantly improve MIP solver performance on both large-scale real-world application datasets and MIPLIB. We believe even bigger improvements are possible with further advances in models and algorithms. 

Some promising future directions are:
\begin{itemize}
    \item Learning to cut: Better selection and generation of cuts using learning is another potential source of performance improvements. \cite{tang2020icml} provide evidence for the usefulness of this direction.
    \item Warm-starting models: The strong performance of learned models on MIPLIB suggests that it is possible to learn heuristics that work well across diverse MIPs. This can be used to overcome the `cold-start' problem in applications where the amount of training data available early on in an application's life cycle may be too small to train good models. We can start by using models trained on heterogeneous datasets and use them as a bridge to more specialized models as more data is collected for the application.
    \item Reinforcement learning: Performance achieved using distillation or behavioral cloning is capped by the best expert available, while reinforcement learning (RL) can potentially exceed it. Efficient exploration, long range credit assignment, and computational scalability of learning are key challenges in applying RL to large-scale MIPs. Addressing them can lead to bigger performance improvements.
\end{itemize}

%% file: acknowledgments.tex
\section{Acknowledgments}
The authors would like to thank 
\begin{itemize}
    \item Alek Andreev, Matko Bosnjak, Elena Buchatskaya, Andrew Cowie, Vali Irimia, Vladimir Macko, Sona Mokra, Rahul Palamuttam, Alvaro Sanchez, and Kevin Waugh for their help in implementing various neural network models,
    \item David Applegate for advice on Neural Diving and providing the Google Production Packing dataset,
    \item Dave Helstroom for help with the Google Production Planning dataset,
    \item Jakob Bauer, Michael King, Andrea Michi, Richard Tanburn, and Sims Witherspoon for providing the Electric Grid Optimization dataset,
    \item Dj Dvijotham and Sven Gowal for providing the Neural Network Verification dataset, and
    \item Theophane Weber for proofreading the paper and suggesting improvements,
    \item Miles Lubin for providing help with PySCIPOpt and references for Algorithm Configuration and Algorithm Selection, and
    \item Ross Anderson and Juan Pablo Vielma for suggesting improvements to the paper.
\end{itemize}

%% file: appendix.tex
\section{Appendix}
In the appendix we provide supplementary information which can be grouped into four strands:
\begin{enumerate}
    \item An autoregressive model approach for Neural Diving, which produces promising initial results, but requires too much further investigation to warrant its inclusion in the main text.
    \item Further data shedding light onto the previously described Neural Brancher work from a slightly
different perspective, namely analysis results with respect to the number of nodes in the branch-and-bound tree,
and also PAR-$k$ metrics.
    \item Additional details on the ADMM batch LP solver.
    \item Extended description of the datasets we used throughout the paper and of the concept of calibrated time we 
used to mitigate the problems arising from trying to measure solving times on large clusters of heterogeneous machines.
\end{enumerate} 

\subsection{Autoregressive Models for Neural Diving}
\label{appendix:autoregressive}
Fairly simple conditionally-independent generative models may provide enough capacity for certain problems, especially when optimal solution is simple to predict.
However, when the training set does not contain a dense region of high-quality solutions, such a model might ultimately fail to adequately express uncertainty and thus either learn to predict an overly narrow sub-region of the solution space or output spread probability mass over many infeasible solutions.

One well-proven way of overcoming limitations of conditionally-independent models is to model the desired distribution in auto-regressive way.
Such a model would predict variable values sequentially based on some pre-defined ordering.
We tried a number of simple extensions of the basic conditionally-independent model such as LSTM~\cite{hochreiter1997lstm} on top of the pre-computed variable node embedding, but found this kind of auto-regressive dependencies not expressive enough to provide a clear advantage over the conditionally-independent model.

Another direction we explored was to only provide auto-regressive dependencies through incrementally solving the underlying LP problem as we sample and assign variable values.
Thus, we first compute the variable node embeddings $\{ v_d \}_{d=1}^D$ (see section~\ref{subsec:architecture} for details). 
Then defining $\tilde{x}(x_{<d})$ to be an LP-solution of $M$ with the first $d-1$ variables assigned to particular values, we predict the value of $x_d$ in the following way:
\begin{align}
    p_{\theta}(x | M) &= \prod_{d=1}^D p(x_d | x_{d-1}, \ldots, x_1, M) = \prod_{d=1}^D \text{Bernoulli}(x_d | \mu_d),  \\
    \mu_d &= \text{MLP}(v_d, \tilde{x}_d(x_{<d}); \theta).
\end{align}

In practice, we do not re-compute the LP-solution after assigning each of the variables and only do so after every $K$ ($K=100$ in the experiment below) steps by warm-starting the LP solver with the last state of the previous solution process.

\begin{figure}
    \centering
    \includegraphics[width=0.49\textwidth]{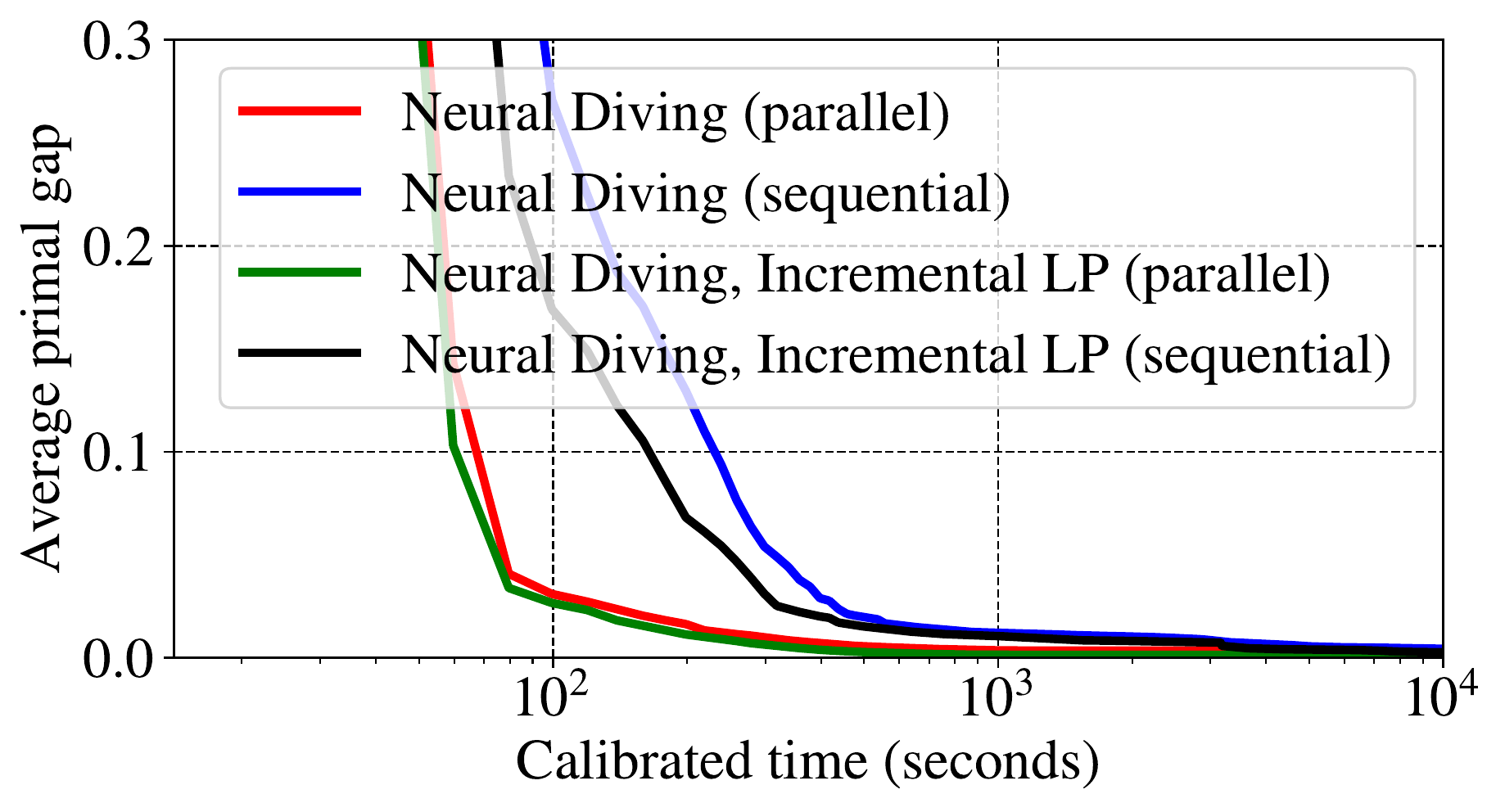}
    \includegraphics[width=0.49\textwidth]{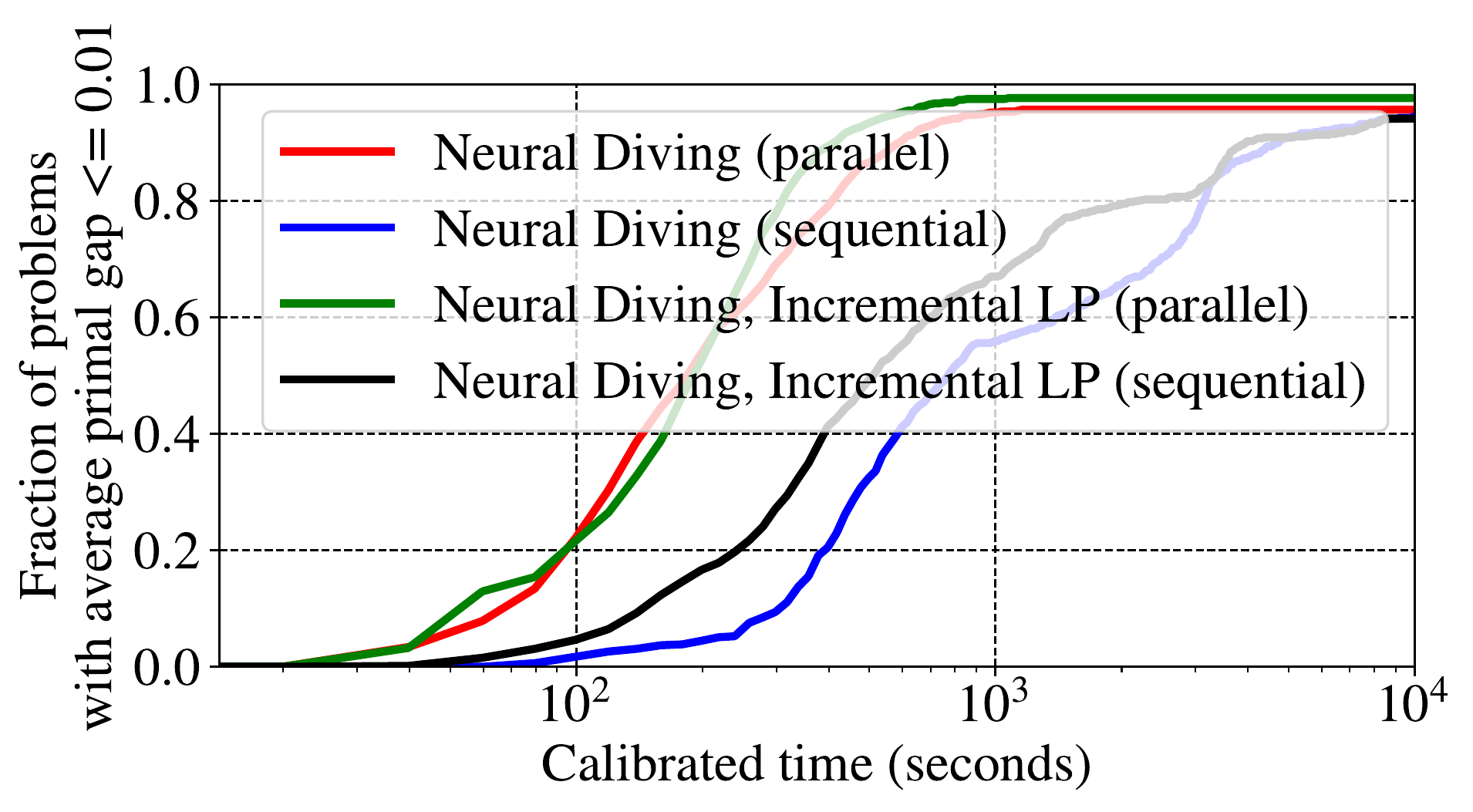}
    \caption{Comparison between the standard conditionally-independent GNN and its autoregressive version with incremental LP solver on Google Production Packing dataset.}
    \label{fig:lp_gnn}
\end{figure}

Figure~\ref{fig:lp_gnn} contains a comparison between the conditionally-independent and the autoregressive models on Google Production Packing dataset.
One can see, that information encoded in incremenetal LP solutions does bring an improvement in terms of the average primal gap and the number of solutions solved to optimality. 
However, most of the gains brought by higher-quality assignments are outweighed by increased sampling time since sampling can no longer be done in parallel for each of the variables.
These preliminary results indicate that there is potential in further architectural advances and leveraging more advanced features, which we leave for future work.

\subsubsection{Ordering of variables}
In many domains, there is often a reasonably natural way to order variables in an autoregressive model. In the case of audio and other temporal sequences, time defines an order. In the case of images, raster scan ordering of pixels is commonly used. 

In the case of integer programs, we need a way to order the binary variables. This ordering should be invariant to row and column permutations of the integer program since such changes do not affect the integer program. Intuitively we want a "canonical" ordering that places variables that are "similar" across instances at similar positions in the ordering so that autoregressive models can learn effectively. Unfortunately, there is no universal canonical ordering of variables that works well for all MIPs. Instead, we list several possible ordering methods as an additional hyperparameter, and pick the one that works best during hyperparameter sweeps.
The orderings that we consider are:
\begin{enumerate}
    \item \textbf{Input order}: We can proceed through the variables in the order as they are given in the MIP definition. Often, the ordering in which the variables are written down in the input already places "similar" variables close to each other.
    \item \textbf{Coefficient order}: We can sort variables by the coefficient by which they contribute to the objective. In some use cases, the objective coefficient of a variable encodes a semantic meaning of importance or priority, and we can hope to use this to e.g. process the most important variables first.  
    \item \textbf{Fractionality order}: We can order variables by their fractionality in the LP solutions. The fractionality of a variable is also often used to determine variable order in diving heuristics in the literature (see, e.g., section 3.1 of \cite{berthold2006thesis}).
    Empirically we found this ordering to perform the best with autoregressive models.
\end{enumerate}

\subsection{Imitation Accuracy of the Learned Branching Policy}
\begin{figure}
\centering
  \begin{subfigure}[t]{0.45\textwidth}
    \centering
    \includegraphics[width=\linewidth]{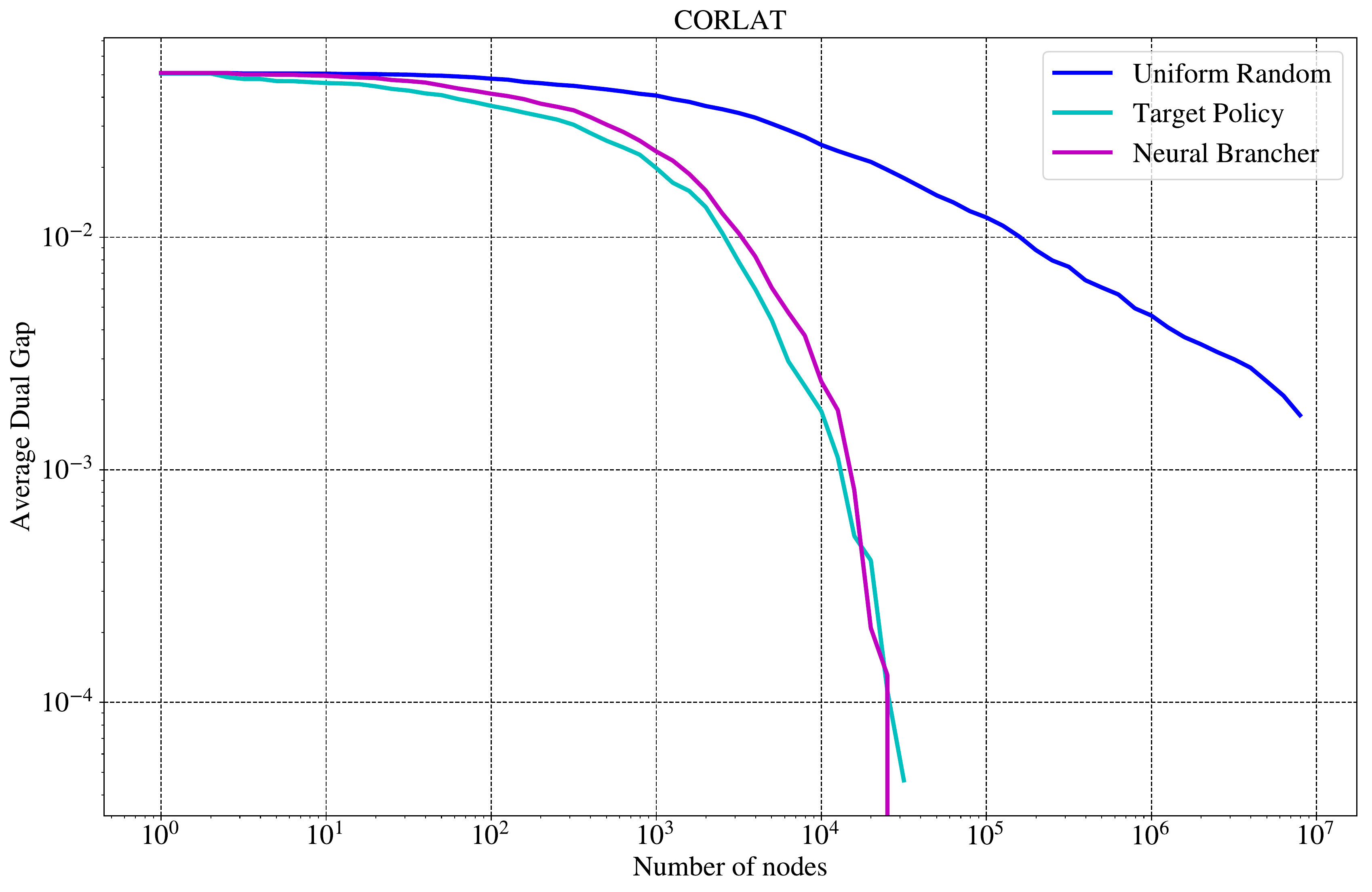}
    \caption{CORLAT Training Set}
  \end{subfigure}
  \begin{subfigure}[t]{0.45\textwidth}
    \centering
    \includegraphics[width=\linewidth]{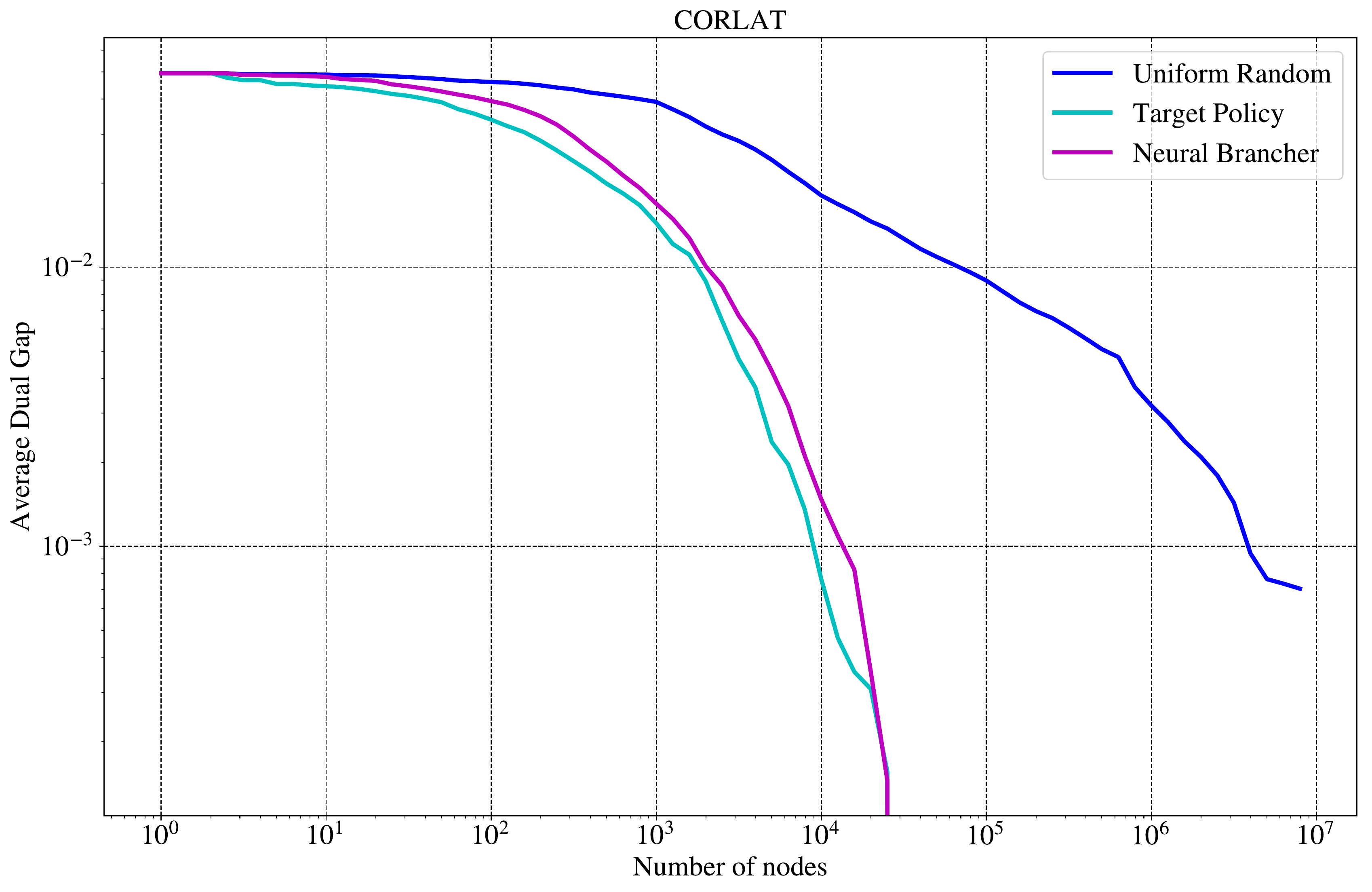}
    \caption{CORLAT Test Set}
  \end{subfigure}\\
  \begin{subfigure}[t]{0.45\textwidth}
    \centering
    \includegraphics[width=\linewidth]{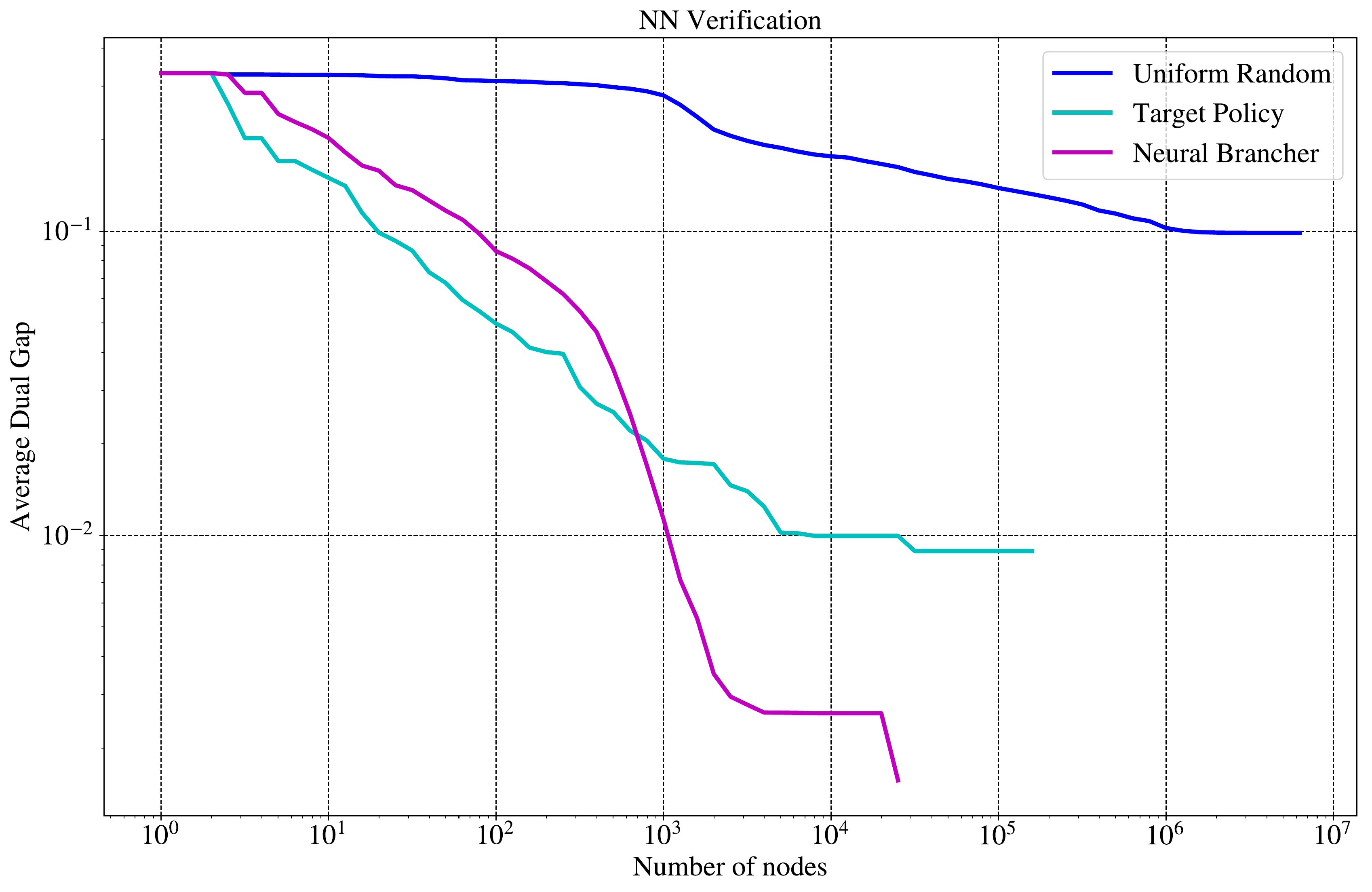}
    \caption{Neural Network Verification Training Set}
  \end{subfigure}
  \begin{subfigure}[t]{0.45\textwidth}
    \centering
    \includegraphics[width=\linewidth]{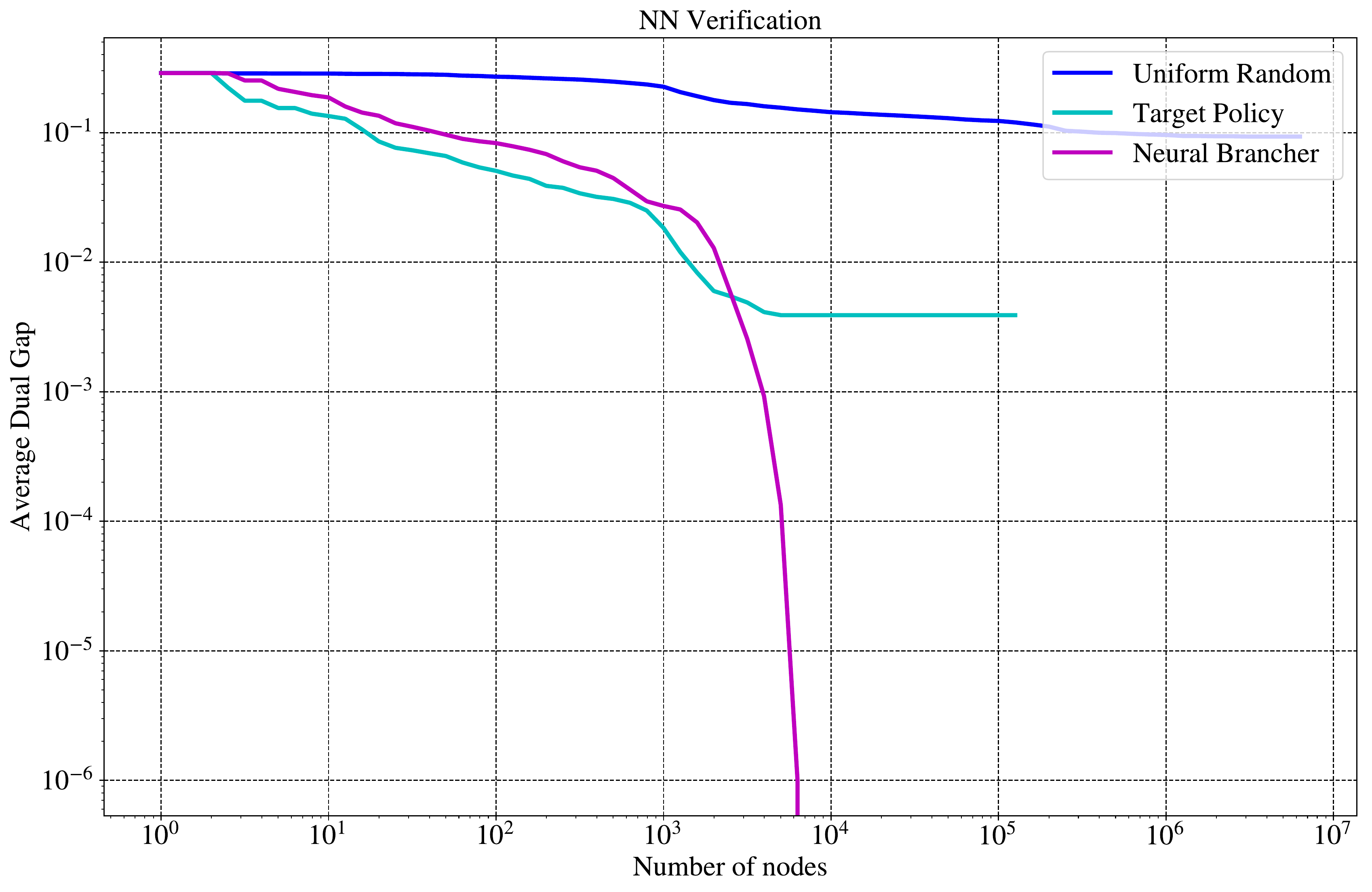}
    \caption{Neural Network Verification Test Set}
  \end{subfigure}
\caption{Comparison of the average dual gap achieved by the target (expert) policy and the learned policy on training and test sets for CORLAT and Neural Network Verification. We also include the uniform random policy as a baseline.}
\label{fig:neural_Branching_target_policy_comparison}
\end{figure}

How well does the learned branching policy approximate the expert policy that it is trained to imitate? Figure~\ref{fig:neural_Branching_target_policy_comparison} compares the average dual gap achieved by the expert policy and the learned policy on both the training and test sets for CORLAT and Neural Network Verification. We also include a policy that takes actions uniformly randomly as a baseline. The results show that imitation learning succeeds in accurately approximating the expert policy both on training and test sets. The learned policy is much closer to the expert than it is to the uniform random policy. The similarity in performance between training and test sets show that the learned policy generalizes well to unseen instances.

\subsection{Evaluation Results with Respect to Number of Nodes}
\begin{figure}
    \centering
    \begin{tabular}{cc}
        \includegraphics[width=0.48\textwidth]{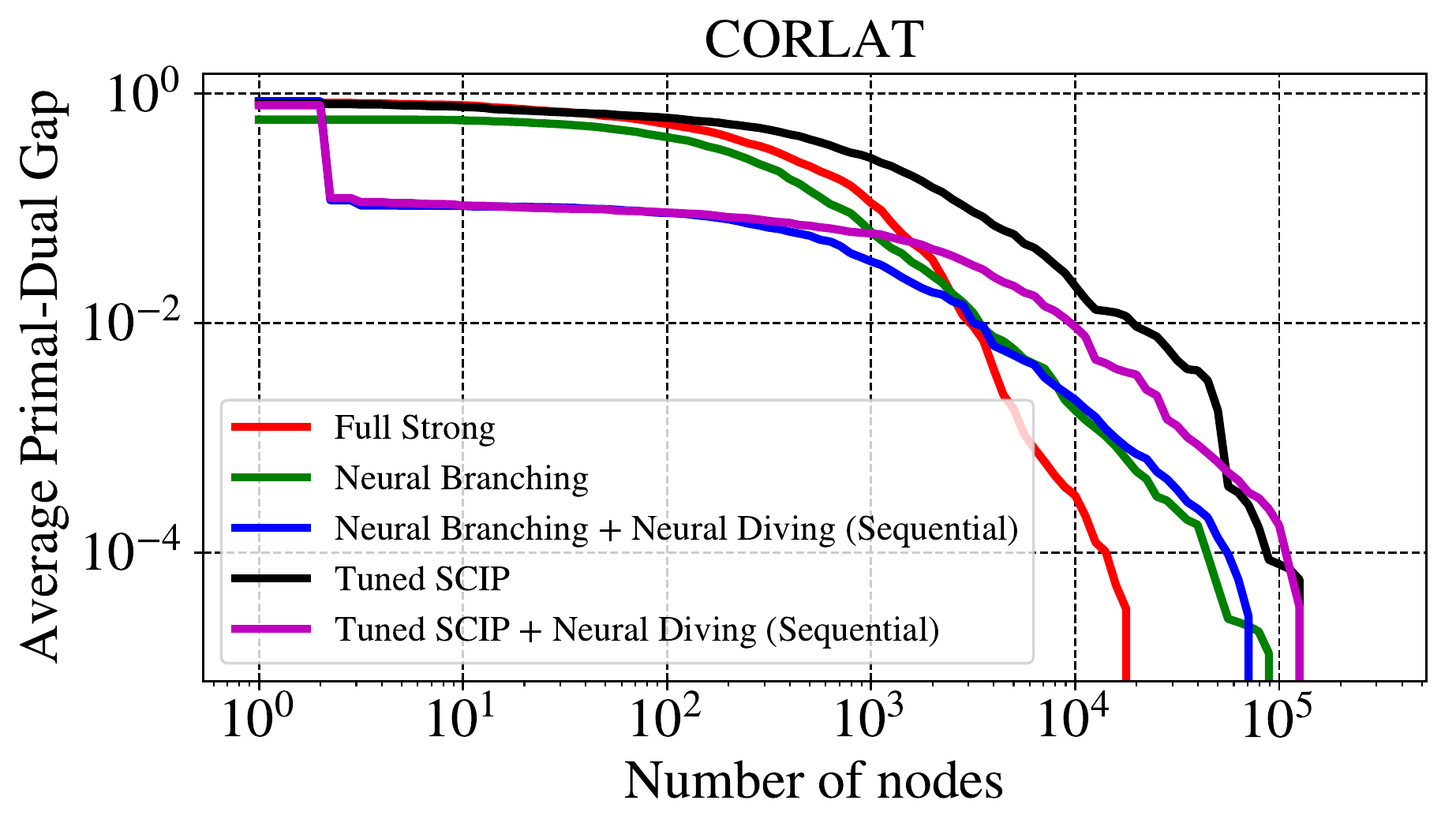} &
        \includegraphics[width=0.48\textwidth]{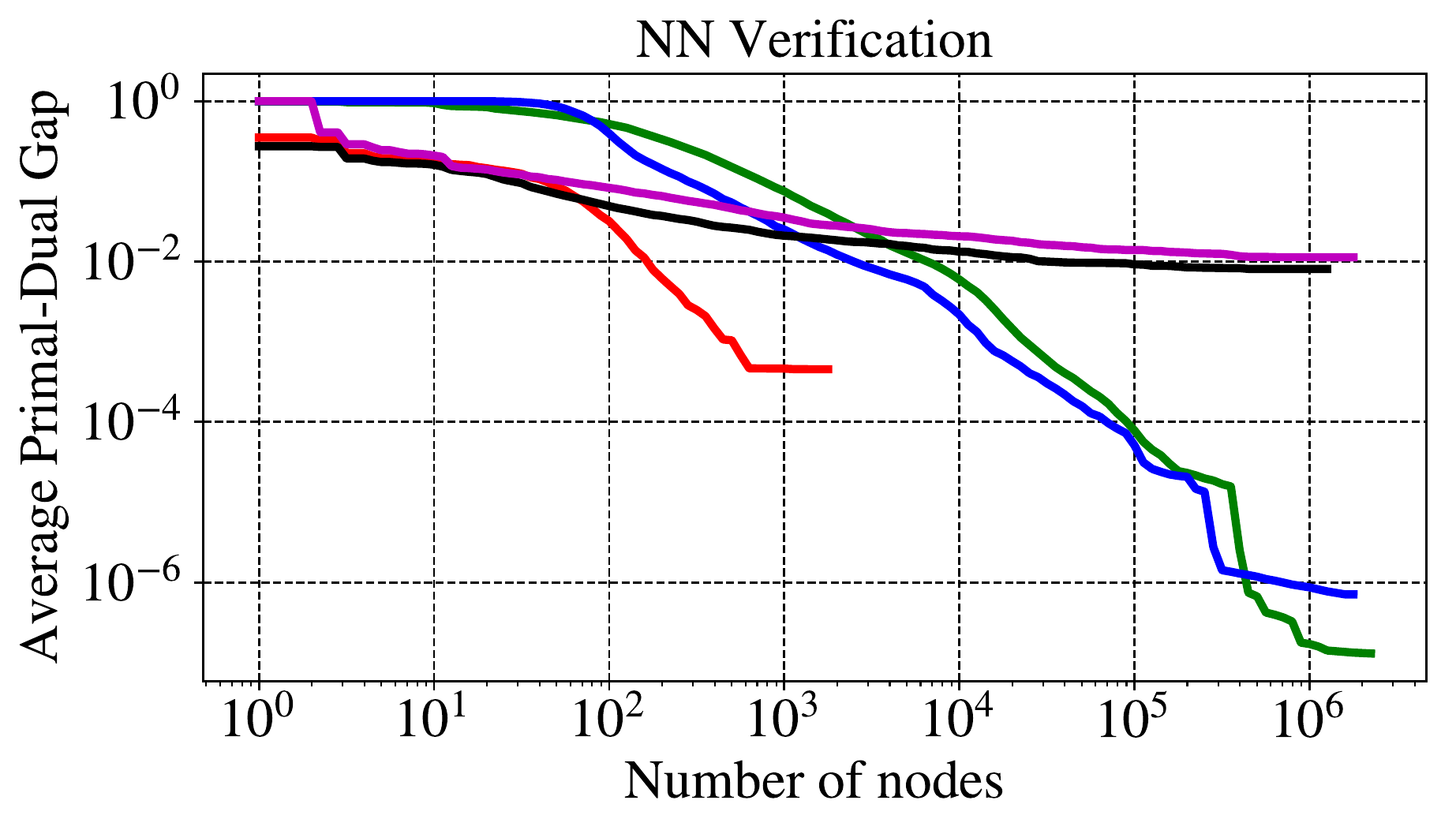} \\
        \includegraphics[width=0.48\textwidth]{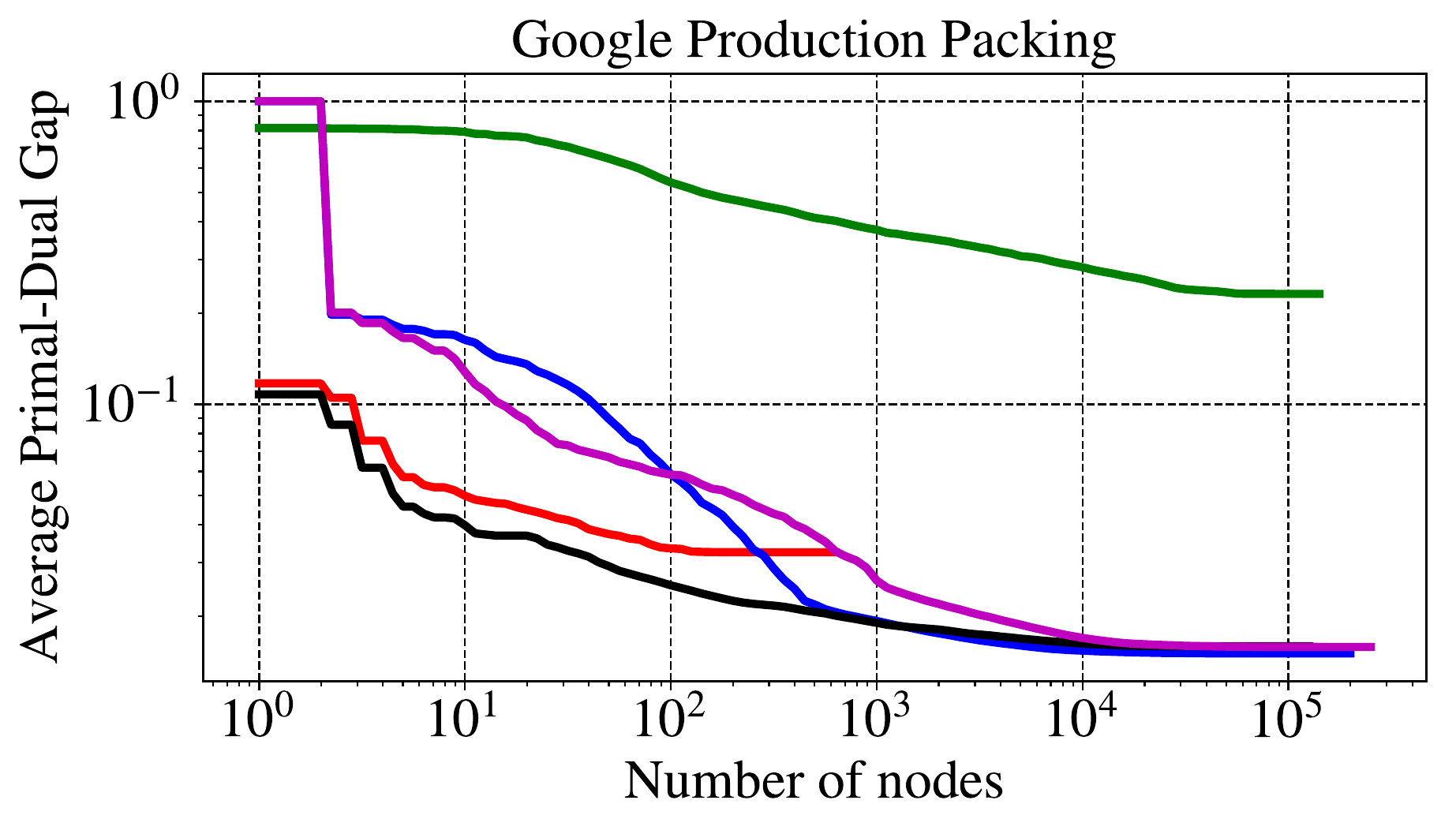} &
        \includegraphics[width=0.48\textwidth]{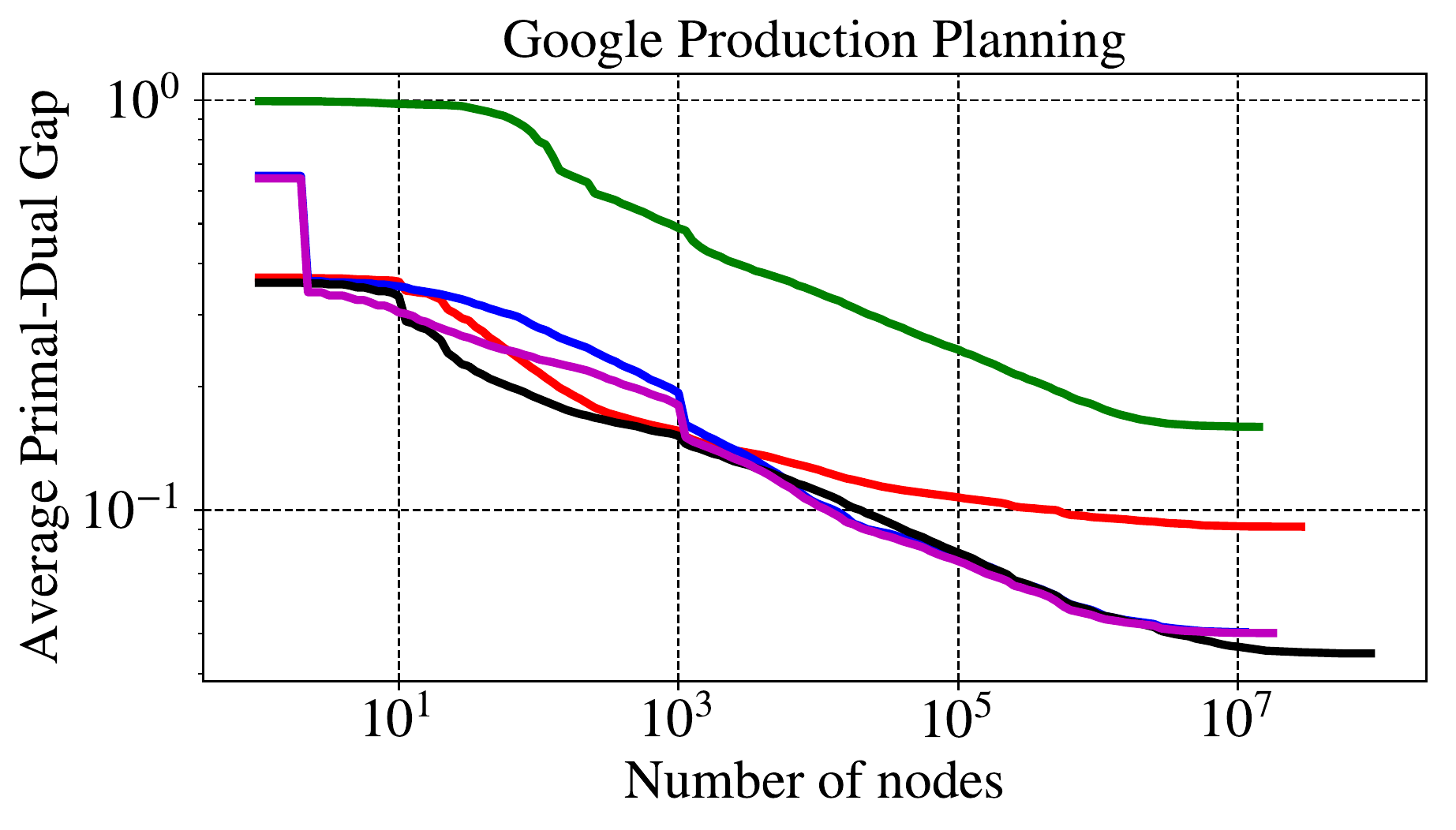} \\
        \includegraphics[width=0.48\textwidth]{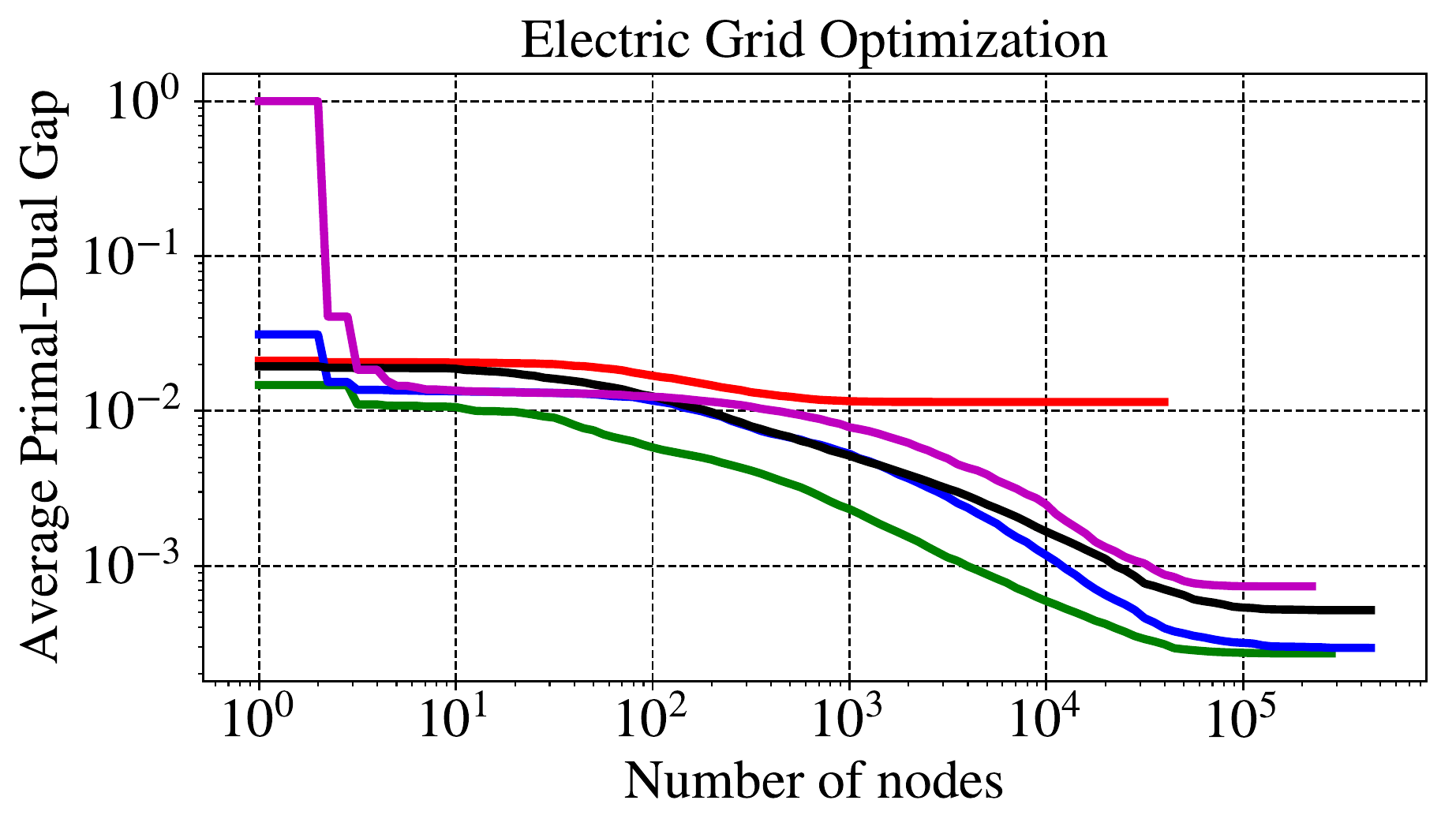} &
        \includegraphics[width=0.48\textwidth]{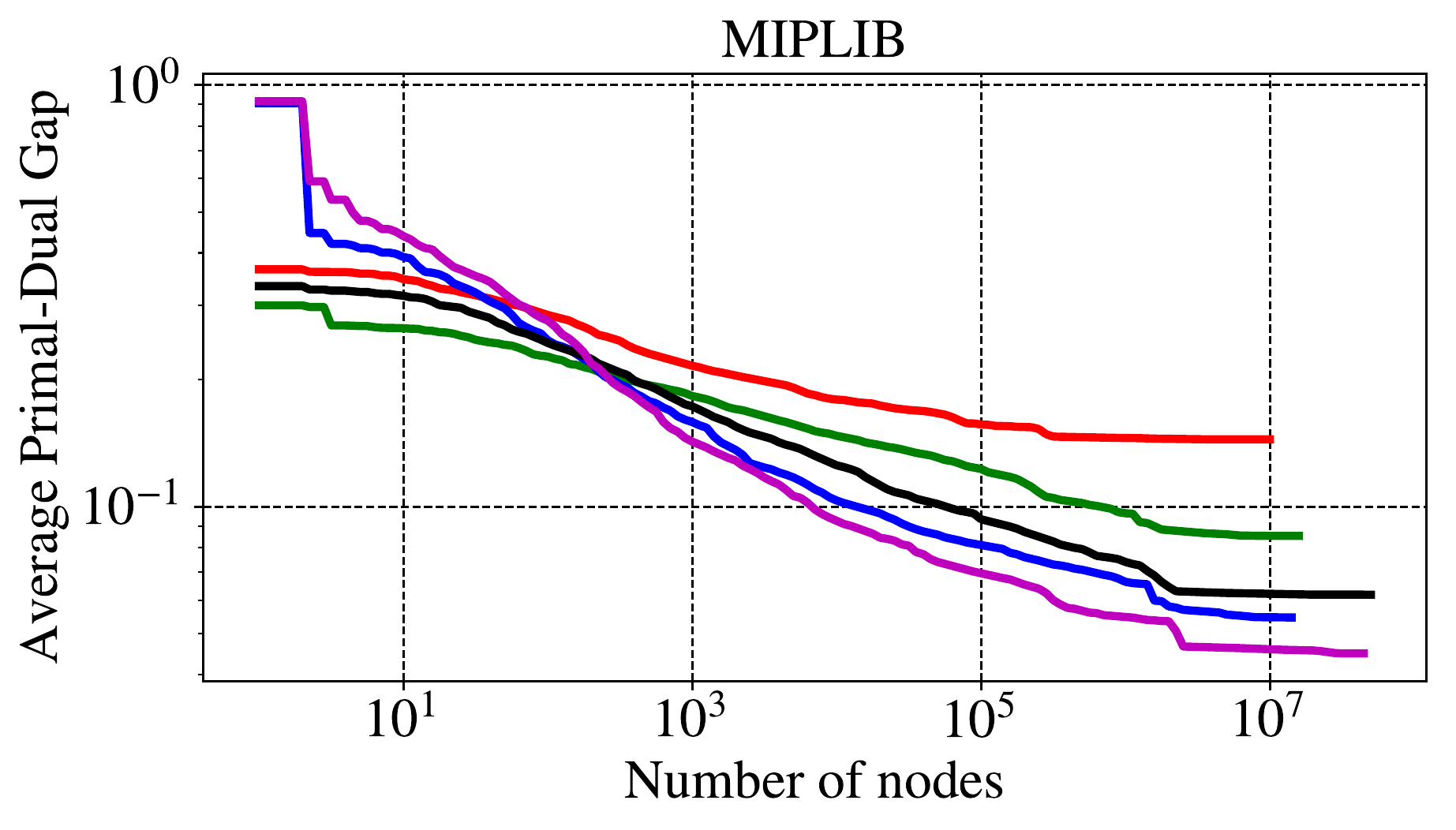} \\
    \end{tabular}
    \caption{Average primal-dual gap achieved on the test set by various solvers as a function of the number of nodes.}
    \label{fig:primal_dual_gap_vs_num_nodes}
\end{figure}

Sections~\ref{subsec:deepbrancher_results} and~\ref{sec:joint_eval} present evaluation results for Neural Branching and the combined Neural Branching + Neural Diving solver with respect to calibrated running time. Here we present the corresponding results with respect to number of nodes. They are shown in Figure \ref{fig:primal_dual_gap_vs_num_nodes}.
While running time is often the practically relevant resource metric, evaluating with respect to number of nodes can provide additional insights by removing the effect of different computational costs and hardware for the solvers being compared.

\begin{figure}
    \centering
    \begin{tabular}{cc}
        \includegraphics[width=0.48\textwidth]{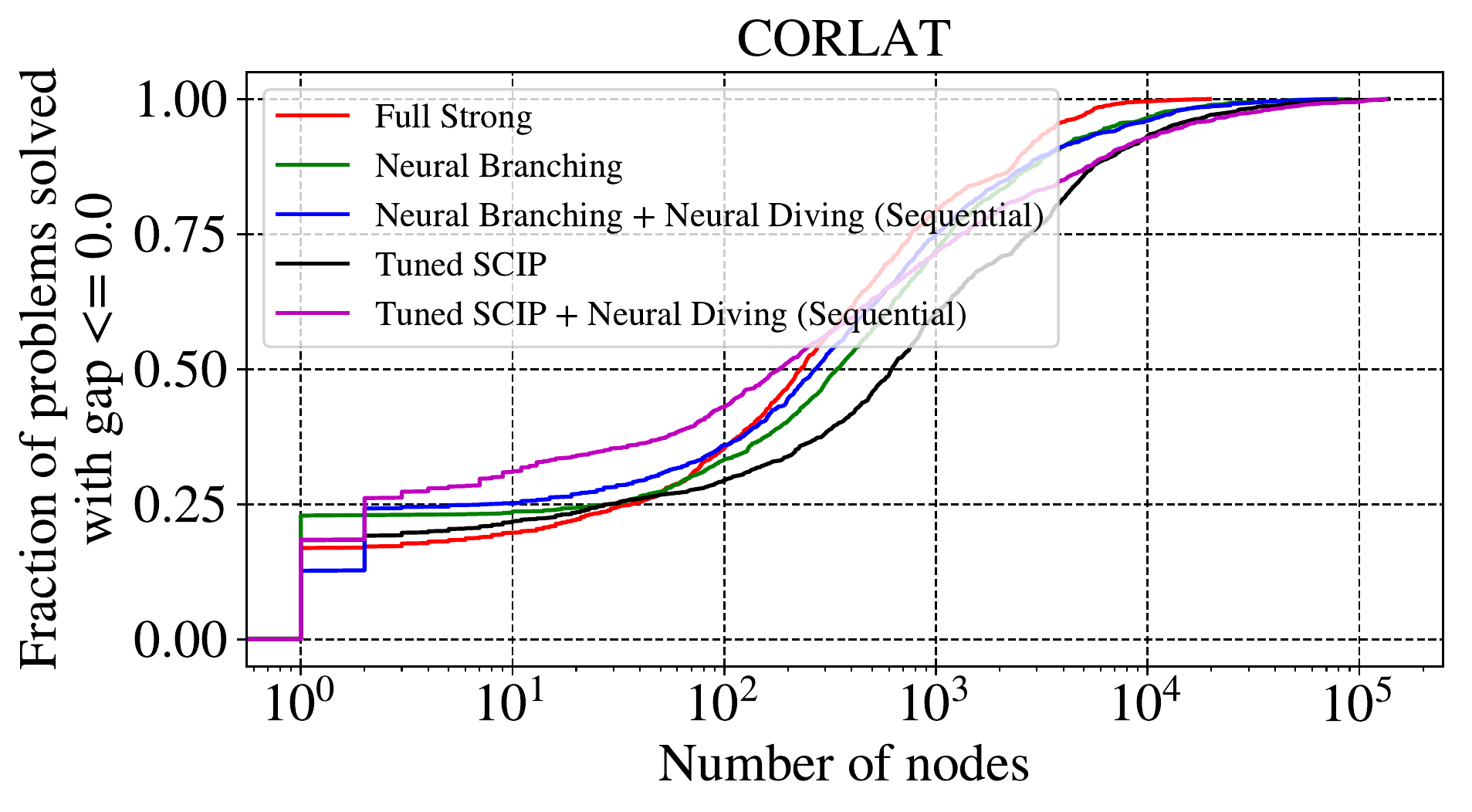} &
        \includegraphics[width=0.48\textwidth]{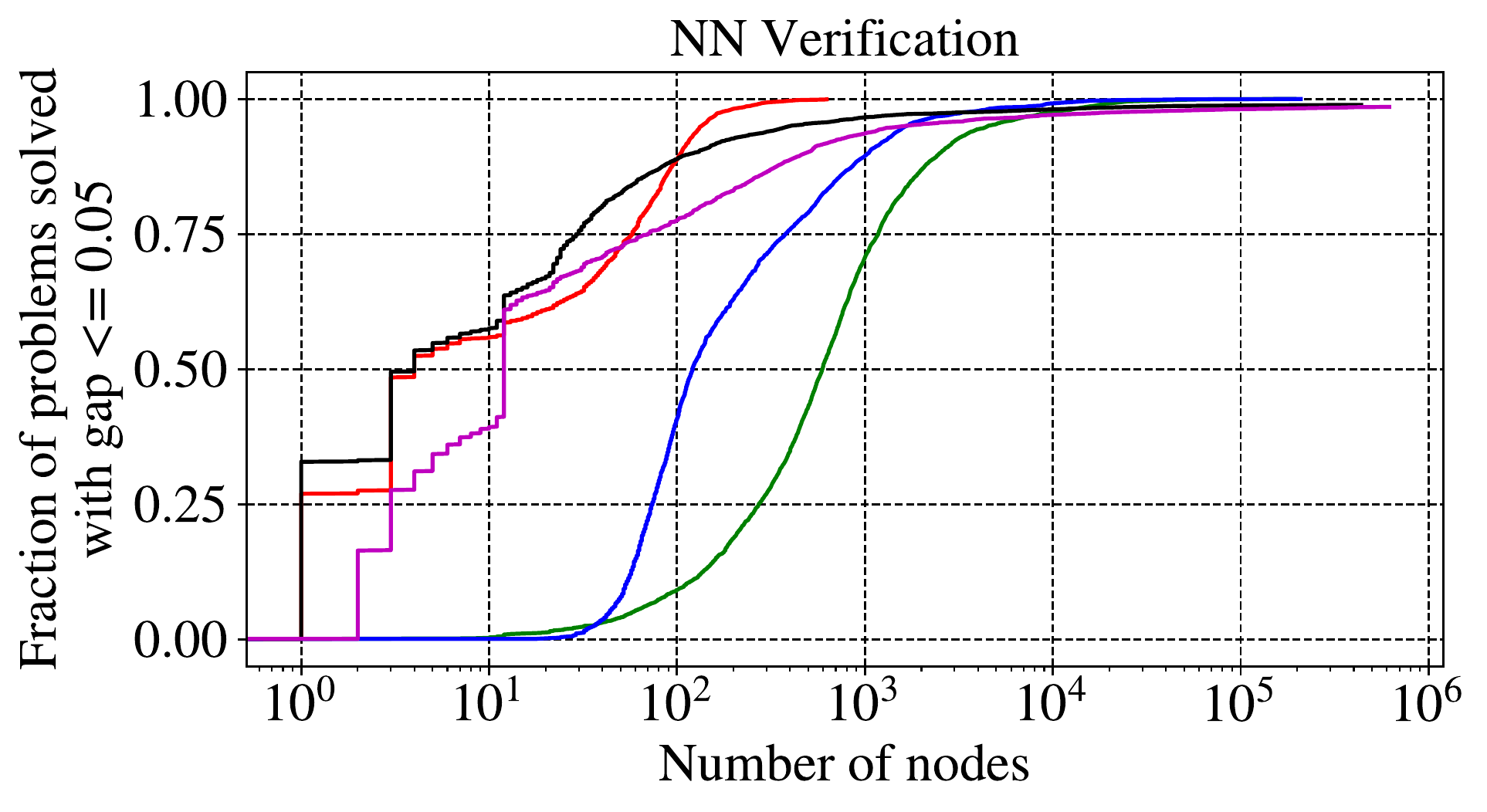} \\
        \includegraphics[width=0.48\textwidth]{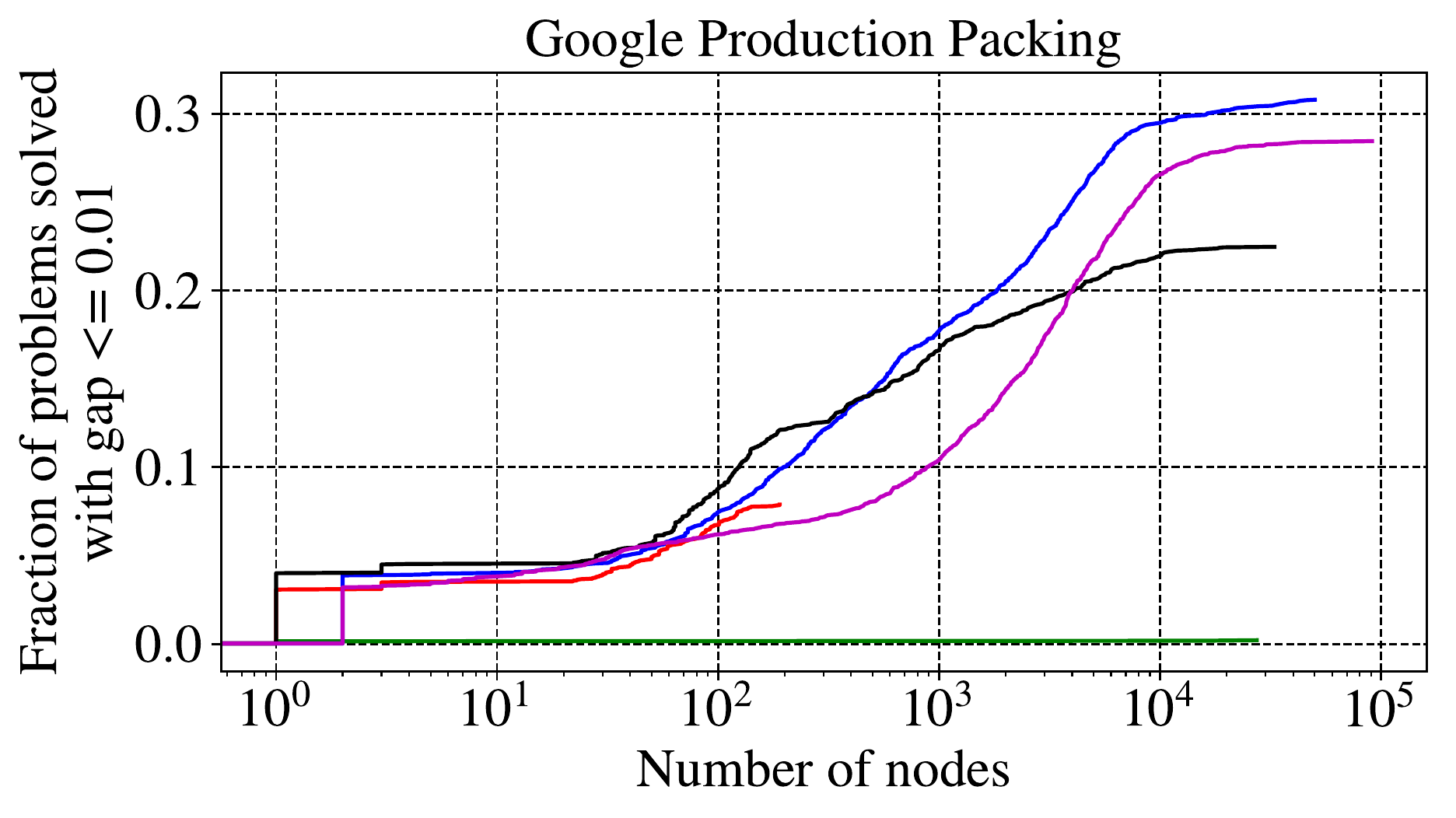} &
        \includegraphics[width=0.48\textwidth]{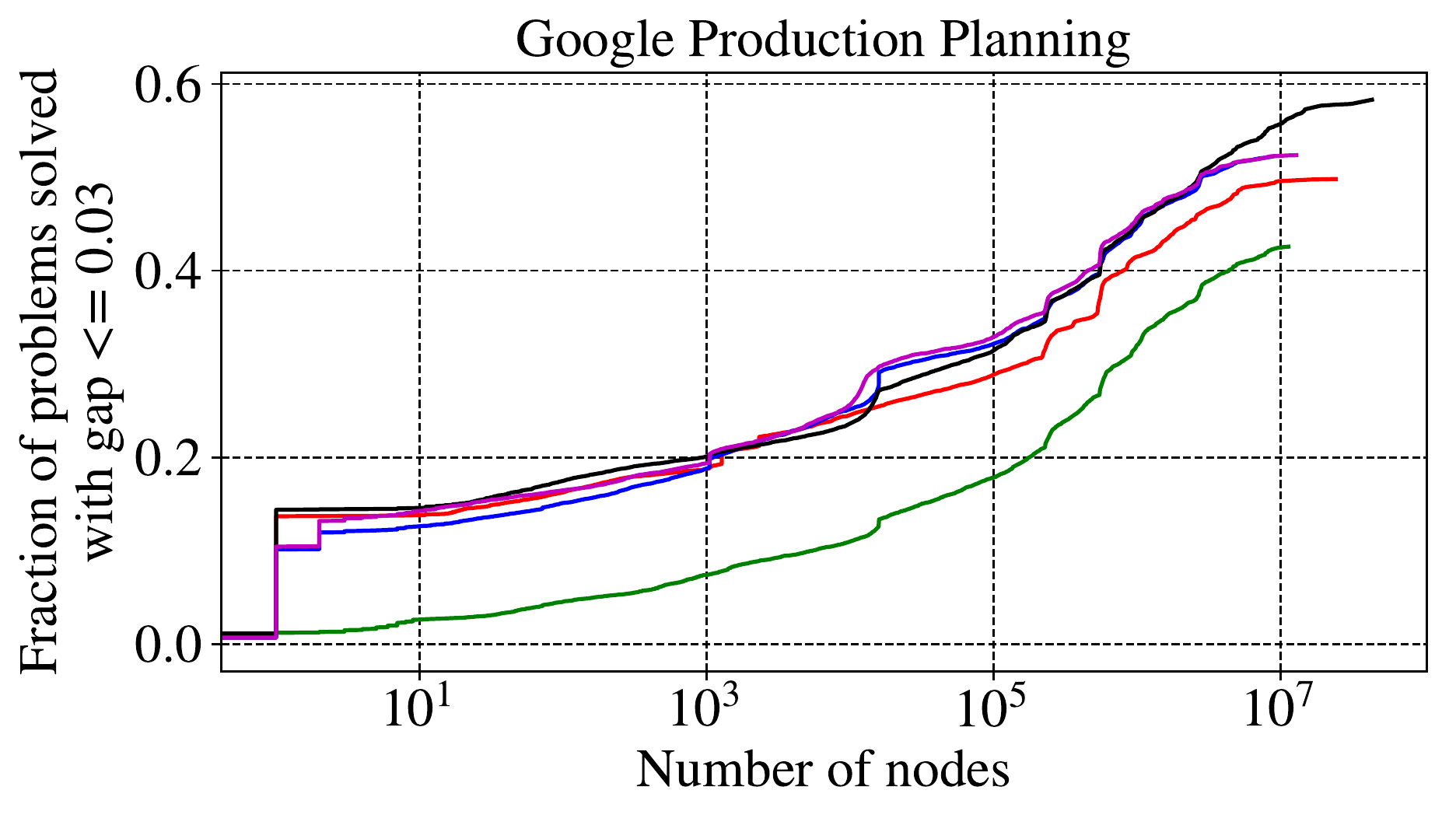} \\
        \includegraphics[width=0.48\textwidth]{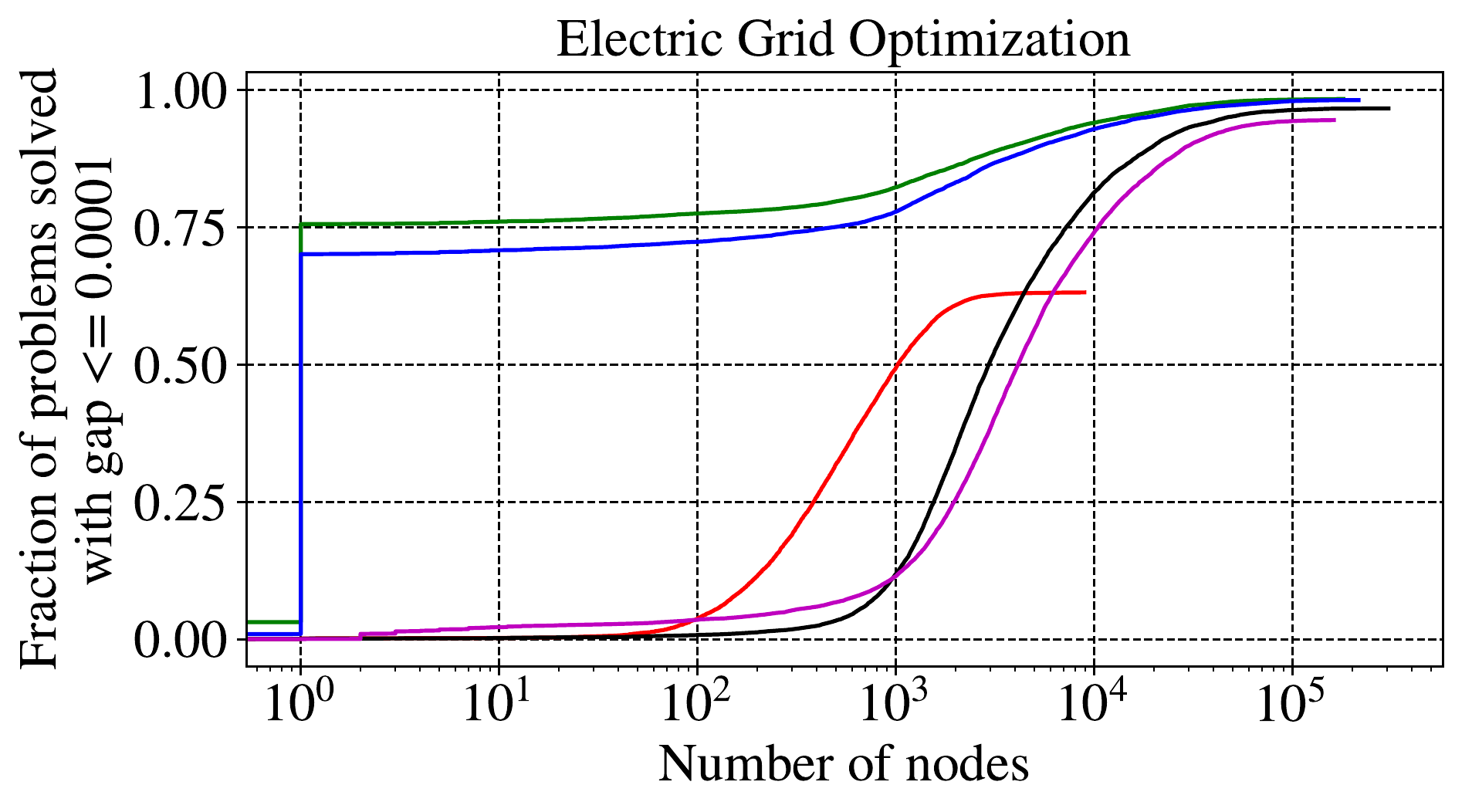} &
        \includegraphics[width=0.48\textwidth]{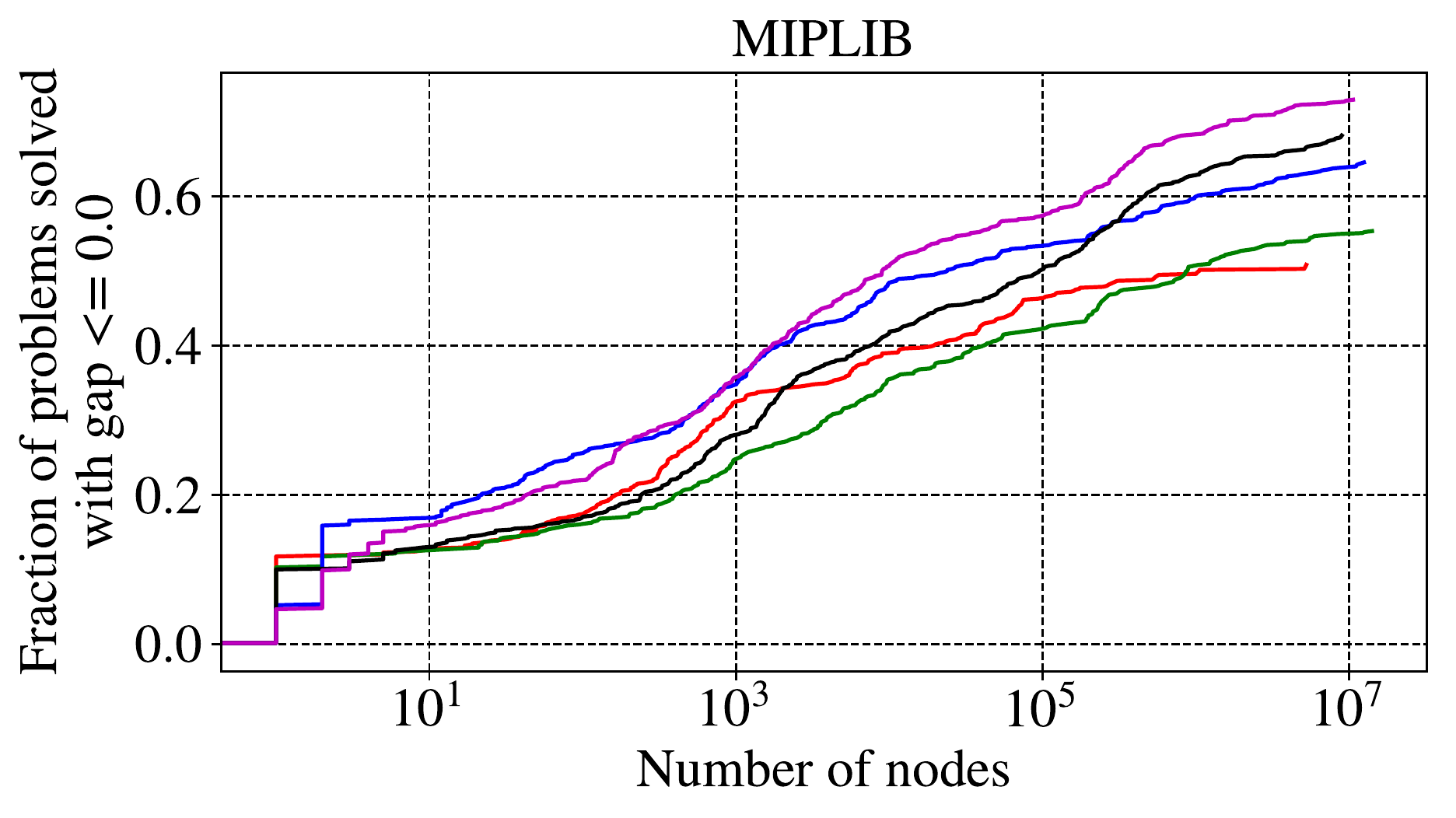} \\
    \end{tabular}
    \caption{Survival plots with respect to number of nodes for various combinations of primal heuristics and variable selection policies as a function of running time on the benchmark datasets.}
    \label{fig:survival_plots_num_nodes}
\end{figure}

\subsection{PAR-10 Results}
\label{subsec:par10}

In Table~\ref{table:par10} we present Penalized Average Running time results with a penalty factor of 10 (PAR-10). PAR-10 averages over a set of instances the time taken to achieve the target optimality gap, with timed out instances being assigned a running time equal to the time limit multiplied by the penalty factor. In our setting we use a time limit of 10000. As can seen in the Table~\ref{table:par10} the winner for four of the five datasets with the largest MIPs is a Neural Solver, supporting our claim that learning is beneficial.

\begin{table}[h!]
\begin{center}
\begin{tabular}{ | m{3cm} | m{1.7cm} | m{2cm} | m{2cm} | m{2cm} | m{2cm} | m{2cm} |} 
 \hline
 Solver & CORLAT & NN Verification & Google Production Packing & Electric Grid Optimization & Google Production Planning & MIPLIB\\
 \hline\hline
 Full Strong & 21.33 & 120.86 & 92389.16 & 39982.2 & 50899.31 & 50097.96\\
 \hline
 Tuned SCIP & \textbf{6.10} & 1128.74 & 77338.42 & 4398.12 & \textbf{42148.14} & 32818.75  \\
 \hline
 Neural Branching &76.81 & 21.72 & 99793.69 & \textbf{2311.02} & 58392.5 & 45637.33  \\
 \hline
 Tuned SCIP + Neural Diving (Sequential) & 6.93 & 1499.49 & 72032.52 & 7001.28 & 48225.51 & \textbf{27762.1} \\
 \hline
 Neural Branching + Neural Diving (Sequential) & 7.35 & \textbf{10.20}& \textbf{70026.23} & 2450.8 & 48819.25 & 36383.16  \\
 \hline\hline
 The PAR-10 of the \textbf{Winner} divided by Tuned SCIP's & 1 & 0.009 & 0.9055  & 0.5255  & 1 & 0.8459 \\
 \hline
\end{tabular}
\end{center}
\caption{PAR-10 results with respect to calibrated time. Lower is better. In \textbf{bold} is the best one for each dataset.}
\label{table:par10}
\end{table}

\subsection{ADMM Batch LP Solver}
\label{subsec:appendix_admm_batch_lp_solver}
In this section we give more details about the ADMM batch LP solver we use to generate the
expert training data that the neural Branching is trained on.
Consider the following composite convex optimization problem
\[
\begin{array}{ll}
\mbox{minimize} & f(x) + g(z)\\
\mbox{subject to} & x = z
\end{array}
\]
over variables $x \in \reals^n$, $z \in \reals^n$ and where $f : \reals^n
\rightarrow \reals \cup \{ \infty \}$, $g : \reals^n
\rightarrow \reals \cup \{ \infty \}$ are closed, proper, convex functions.
The \emph{Alternating direction method of multipliers} (ADMM) \citep{boyd2011distributed, lions1979splitting, eckstein1992douglas} applied to this
problem is the following procedure for any $\rho > 0$:
\[
\begin{array}{ll}
  x^{k+1} &= \argmin_x \Big(f(x) + (\rho/2) \|x - z^k - \lambda^k \|_2^2\Big)\\
  z^{k+1} &= \argmin_z \Big(g(z) + (\rho/2) \|x^{k+1} - z - \lambda^k\|_2^2\Big)\\
\lambda^{k+1} &= \lambda^k + x^{k+1} - z^{k+1}.
\end{array}
\]
Under benign conditions this procedure produces a sequence of iterates that converge
to optimal objective values, \ie, $f(x^k) + g(z^k) \rightarrow f(x^\star) + g(z^\star)$, and $\|x^k - z^k\| \rightarrow 0$, where $x^\star$
and $z^\star$ are optimal primal points \citep[\S 3.2]{boyd2011distributed}.
We can write the LP relaxation of MIP \eqref{eqn:mip1} in a form more amenable to ADMM:
\[
\begin{array}{ll}
\mbox{minimize} & c^T x + \Ic_{Ax = y}(x, y) + \Ic_{[b_l, b_u]}(\tilde y)+ \Ic_{[l, u]}(\tilde x) \\
\mbox{subject to} & x = \tilde x, \quad y = \tilde y,
\end{array}
\]
over variables $x, \tilde x \in \reals^n$, $y, \tilde y \in \reals^m$, and
where $\Ic_\Xc$ denotes the convex indicator function of set $\Xc$, \ie,
\[
  \Ic_\Xc(z) = \left\{
    \begin{array}{rl}
      0 & \quad z \in \Xc \\
      \infty & \quad {\rm otherwise.}
    \end{array}
    \right.
\]
When ADMM is applied to this problem it yields the following algorithm with iterates $x, \tilde x \in \reals^n$, $y, \tilde y \in
\reals^m$, $\lambda \in \reals^{m+n}$:
\[
\begin{array}{ll}
(x^{k+1}, y^{k+1}) &=
\argmin_{Ax = y}
\left\{
c^Tx +
\rho/2 \left\| \begin{bmatrix}x \\ y\end{bmatrix} - \begin{bmatrix}\tilde x^k \\
  \tilde y^k\end{bmatrix} - \begin{bmatrix}\lambda_x^k \\ \lambda_y^k\end{bmatrix} \right\|_2^2
\right\}\\
\tilde x^{k+1} &= \Pi_{[l, u]}(x^{k+1} - \lambda_x^k)\\
\tilde y^{k+1} &= \Pi_{[b_l, b_u]}(y^{k+1} - \lambda_y^k)\\
  \begin{bmatrix}\lambda_x^{k+1} \\ \lambda_y^{k+1}\end{bmatrix} &=
    \begin{bmatrix}\lambda_x^k \\ \lambda_y^k\end{bmatrix}+ \begin{bmatrix}x^{k+1} \\ y^{k+1}\end{bmatrix} - \begin{bmatrix}\tilde x^{k+1}\\ \tilde y^{k+1}\end{bmatrix},
\end{array}
\]
where $\Pi_{[l, u]}(z)$ denotes the Euclidean projection of $z$ onto the interval $[l, u]$ elementwise.
It turns out that only first step is computationally challenging since it involves solving
a (potentially large) set of linear equations. Concretely, for some right-hand side $r
\in \reals^{m+n}$ we want to solve
\[
  (x^{k+1}, y^{k+1}) = \argmin_{Ax = y} \left\{c^Tx + \rho/2 \left\| \begin{bmatrix}x \\
    y\end{bmatrix} - \begin{bmatrix}r_x^k \\ r_y^k \end{bmatrix} \right\|_2^2\right\}.
\]
The solution to the which is given by
\[
  \begin{bmatrix}x^{k+1} \\ y^{k+1} \end{bmatrix} =
\begin{bmatrix}
I & \quad A^T \\ A & \quad -I
\end{bmatrix}^{-1}
\begin{bmatrix}
r_x^k - c / \rho \\ r_y^k
\end{bmatrix} +
\begin{bmatrix}
0 \\ r_y^k
\end{bmatrix}.
\]
Note that the matrix to be inverted is the same \emph{for all iterations}.
Not only that, but the matrix does not depend on which variable we are branching on;
it is the same for all LPs we need to solve at any given node.
This means we can batch and solve all the LPs (or however many will fit into the GPU memory) together on the GPU simultaneously.
This yields significant hardware benefits for both direct and
indirect approaches to solving the linear system. Although GPUs are designed to provide very fast dense matrix-matrix routines the sparse matrix operations make less efficient
uses of the GPU stream parallelism. However, by batching
sparse matrix-vector operations into matrix-matrix, then it is easier to keep the GPU
busy and make better use of the parallelism. We shall see that this yields very large
speedups in practice.

\subsubsection{GPU speedups}
\begin{figure}
  \includegraphics[width=0.5\linewidth]{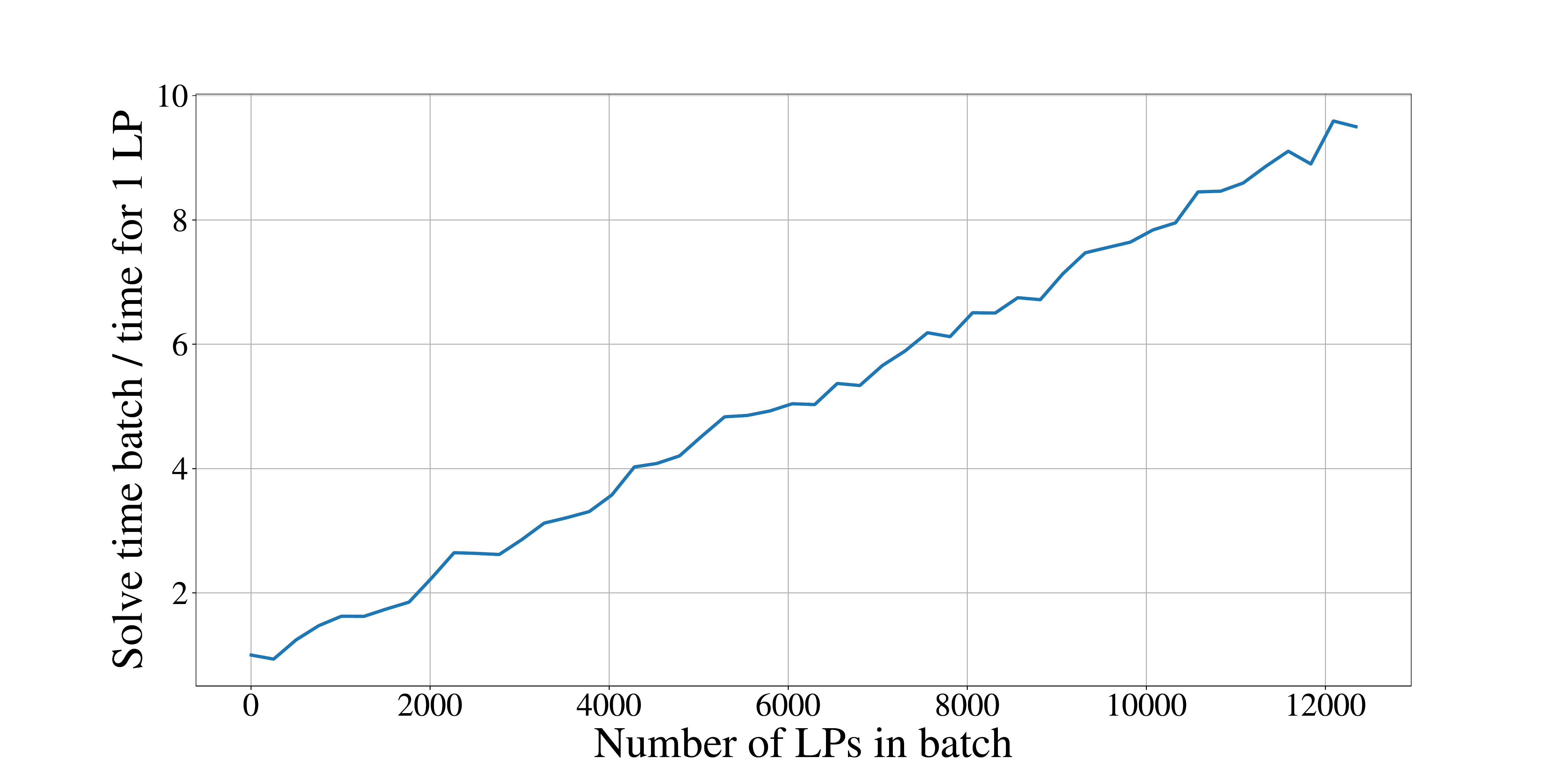}
  \hfill
  \includegraphics[width=0.5\linewidth]{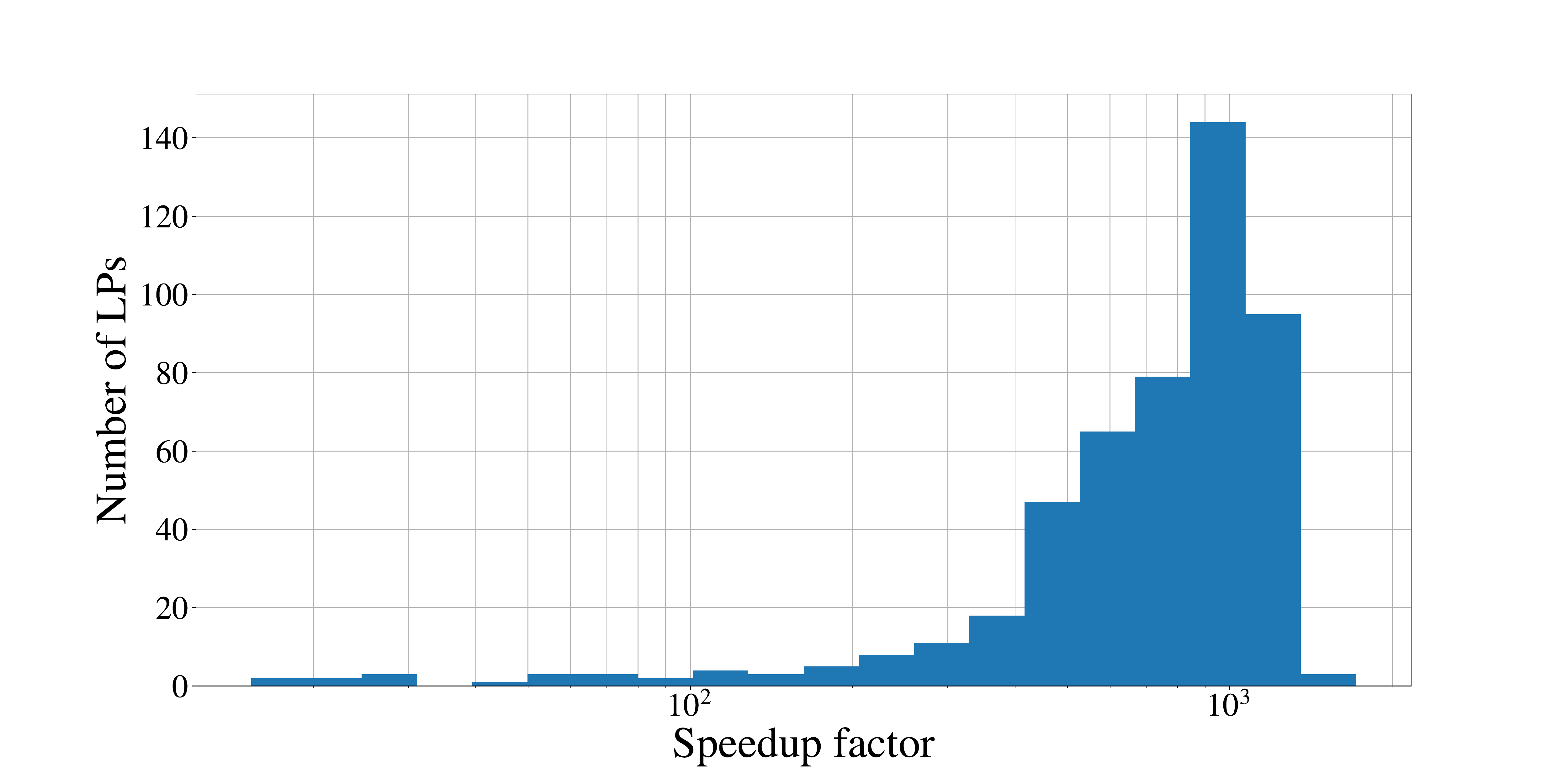}
\caption{Left: Solve time vs number of LPs in a batch for the air05 MIP. Right: Histogram of speedup factors for all MIPs in MIPLIB at root nodes.}
  \label{f-air05}
\end{figure}

Here we present some speedup results for our batching strategy on a single GPU. We
shall consider all possible LPs for all branching candidates at the root nodes of all the MIPs in the MIPLIB dataset.
In all cases we ran the ADMM solver for 100 iterations and used the LP solution
of the root node as a warm-start with which to initialize the ADMM solver.
First, as a concrete example, consider the air05 binary MIP which at the root node (after presolve and any reductions that happens at the root node) has $6.1$k variables, $426$ constraints.
The matrix that defines the problem is $98\%$ sparse with $43.4$k non-zeros. This high sparsity
is not amenable to acceleration on GPU, however, when batched the solver can make much better use
of the GPU parallelism.
In Figure~\ref{f-air05} we show the time required to approximately solve the LPs vs batch-size on GPU, normalized by time for (approximately) solving a single solver using the same solver.
We show that we can compute all 12.2k branch scores in roughly $10$ times the cost of a single solve,
which is approximately a $1200\times$ speedup from batching relative to solving all the LPs
sequentially.
As a point of comparison, for the same set of LPs ECOS \citep{ecos} takes about $4.5$ seconds to solve each one (which is approximately the same as SOPLEX, the built-in LP solver in SCIP, which takes $3$ seconds to solve the LP at the root) so to solve all $12k$ LPs would require more than $15$ hours. This is far too long for use in an industrial MIP solver and significantly larger than the ADMM batch LP solver, which finds the (approximate) solutions for all the LPs in $5.4$ seconds.
Next we show that this pattern holds true for the entire MIPLIB dataset.
In Figure~\ref{f-air05} we show the histogram of this speedup across the root nodes of the entire MIPLIB dataset.
If there are $n$ LPs to be solved at the root node, then we 
define the speedup factor to be
\[
{\rm speedup} = \frac{n \times (\text{time to solve $1$ LP)}}{\text{time to solve $n$ LPs in batch}}.
\]
We remove the effect of the fixed costs (\eg, transferring the data
to the GPU) by subtracting the time taken
to run zero iterations, which is non-zero, from the denominator and numerator (if we didn't remove
this time it would artificially inflate the speedup factors we would observe).
This figure shows that the speedup is significant across the board,
averaging approximately $800\times$, almost three orders of magnitude. Typically
larger MIPs yield larger speedups, which also explains the sudden drop off at around
$1000\times$ speedup, since the largest MIPs eventually become too large
for the GPU memory to handle.

\subsubsection{Accuracy}
\begin{figure}
  \includegraphics[width=0.5\textwidth]{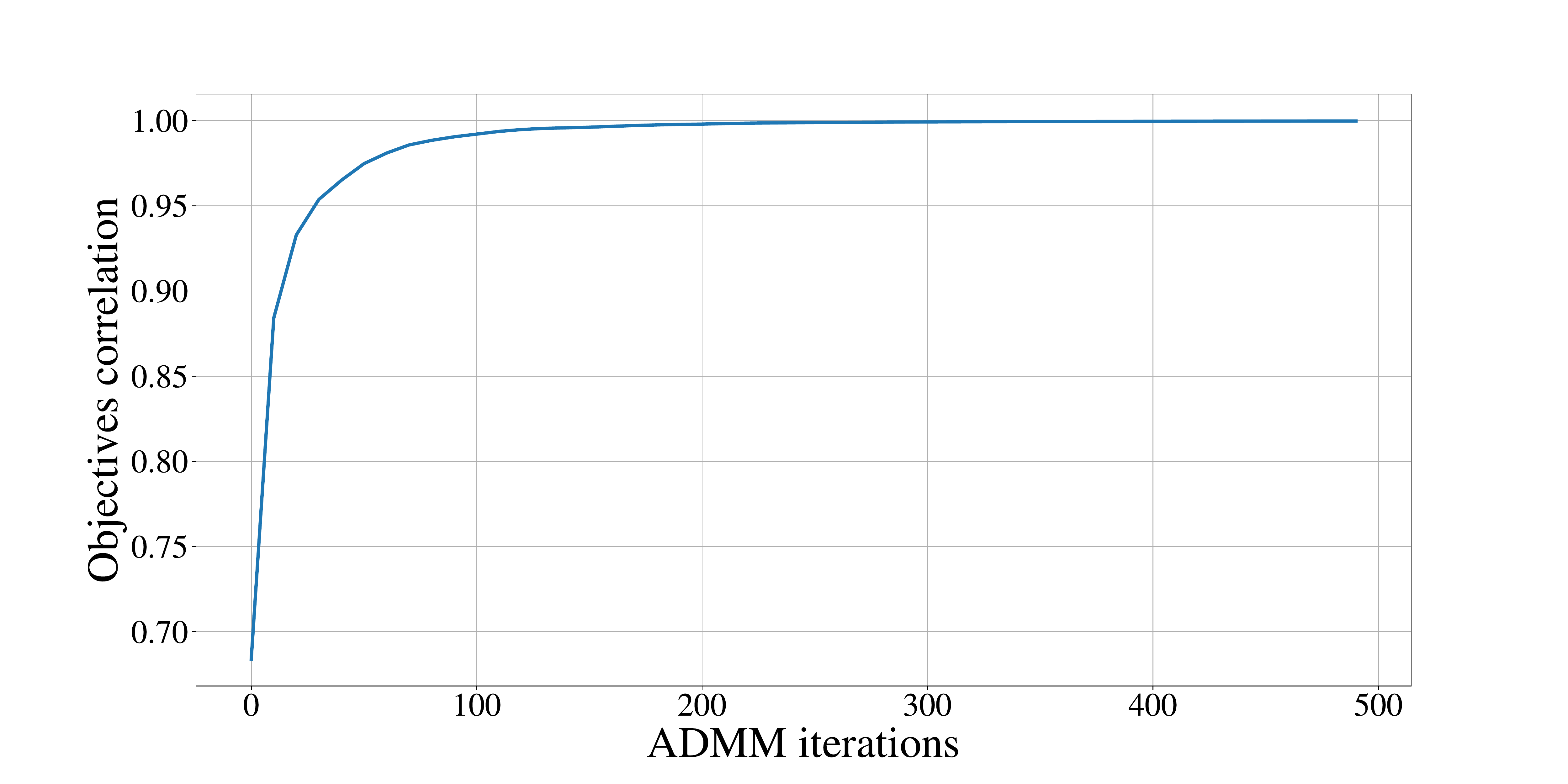}
   \hfill
  \includegraphics[width=0.5\textwidth]{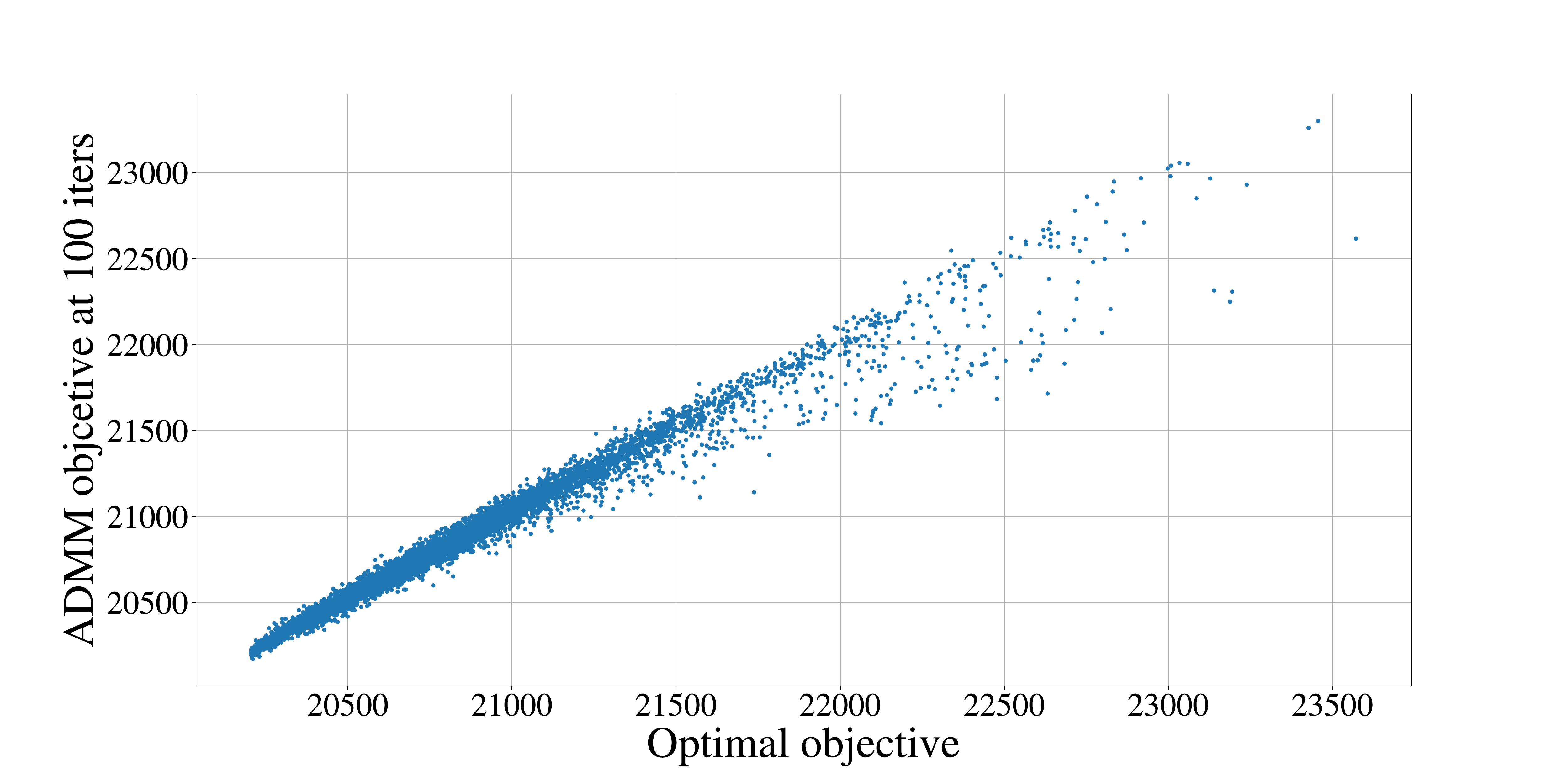}
\caption{Correlation of ADMM objective and true objectives at root of air05 MIP.}
\label{f-accuracy}
\end{figure}

Warm-starting from the solution at the previous node means not many iterations are
required for reasonable accuracy, which is a typical property of ADMM. This justifies
our use of a relatively small number of ADMM iterations to `solve' the problem.
In Figure~\ref{f-accuracy} we plot two quantities at the root node of the air05
MIP, at which we need to solve approximately $12k$ LPs. On the left hand side
we plot the correlation of the $12k$ \emph{true}
objective values with the objective values returned by ADMM after 100 iterations.
The true objective values were found by solving the same LPs using the high-accuracy ECOS interior point
solver \citep{ecos} with parser-solver CVXPY \citep{diamond2016cvxpy}.
As we can see it converges towards one as the number
of iterations increases, but even for a handful of iterations the approximate
objective values are well correlated with the true objective values.
Since the branch and bound decisions are made using the objective values only (see \eqref{eqn:bnbvalue})
this high correlation
indicates that we have a strong signal for the branch and bound decision
with relatively little computational effort, and typically only low accuracy
LP objective values are required for use in FSB \citep{achterberg2005branching}.
On the right we show the
scatterplot of the true objective values and the objective values returned by
ADMM after $100$ iterations, and we can see the high correlation visually.

\subsection{Details of Datasets}
\label{subsec:appendix_datasets}
The target optimality gap used for each dataset in our evaluation as SCIP's stopping criterion for a MIP solve is given in table~\ref{tab:target_opt_gaps}. In the case of Neural Network Verification, the actual criterion used in the application is to stop when the objective value becomes negative, but this is not expressible as a constant target gap across all instances. In order to treat all datasets consistently, we have selected a gap of 0.05.
\begin{table}[]
    \caption{Target optimality gaps for the datasets in our evaluation used as the stopping criterion for SCIP.}
    \centering
    \begin{tabular}{c|c}
        \toprule
        Dataset & Target Optimality Gap \\
        \midrule
        CORLAT &  0\\
        Neural Network Verification & 0.05\\
        Google Production Packing & 0.01 \\
        Google Production Planning & 0.03\\
        Electric Grid Optimization & 0.0001 \\
        MIPLIB & 0\\
        \bottomrule
    \end{tabular}
    \label{tab:target_opt_gaps}
\end{table}

\begin{table}
    \caption{Median of the number of constraints and variables (per type) for the different datasets before / \textit{after} presolving using SCIP 7.0.1.}
    \begin{adjustbox}{center}
    {\footnotesize
    \begin{tabular}{ccccccc}
    \toprule
        Dataset & Constraints & Variables & Binary & Integer (non-binary) & Continuous  \\
        \midrule 
        \emph{CORLAT} & 470/\textit{454} & 466/\textit{458} & 100/\textit{97} & 0/\textit{0} & 366/\textit{361} \\
        \emph{NN Verification} & 6531/\textit{1407} & 7142/\textit{1629} & 171/\textit{170} & 0/\textit{0} & 6972/\textit{1455} \\
        \emph{Google Production Packing} & 36905/\textit{24495} & 10046/\textit{8919} & 3773/\textit{3437} & 0/\textit{0} & 6231/\textit{5421} \\
        \emph{Google Production Planning} & 11910/\textit{478} & 9884/\textit{404} & 833/\textit{119} & 462/\textit{119} & 8337/\textit{136} \\
        \emph{Electric Grid Optimization} & 61851/\textit{45834} & 60720/\textit{58389} & 42240/\textit{42147} & 0/\textit{0} & 18768/\textit{16623} \\
        \emph{MIPLIB} & 7706/\textit{4388} & 11090/\textit{9629} & 4450/\textit{3816} & 0/\textit{0} & 218/\textit{96} \\
        \bottomrule
    \end{tabular}
    }
    \end{adjustbox}
    \label{tab:dataset_stats}
\end{table}

\begin{figure}
    \centering
    \begin{tabular}{cc}
        \includegraphics[width=0.48\textwidth]{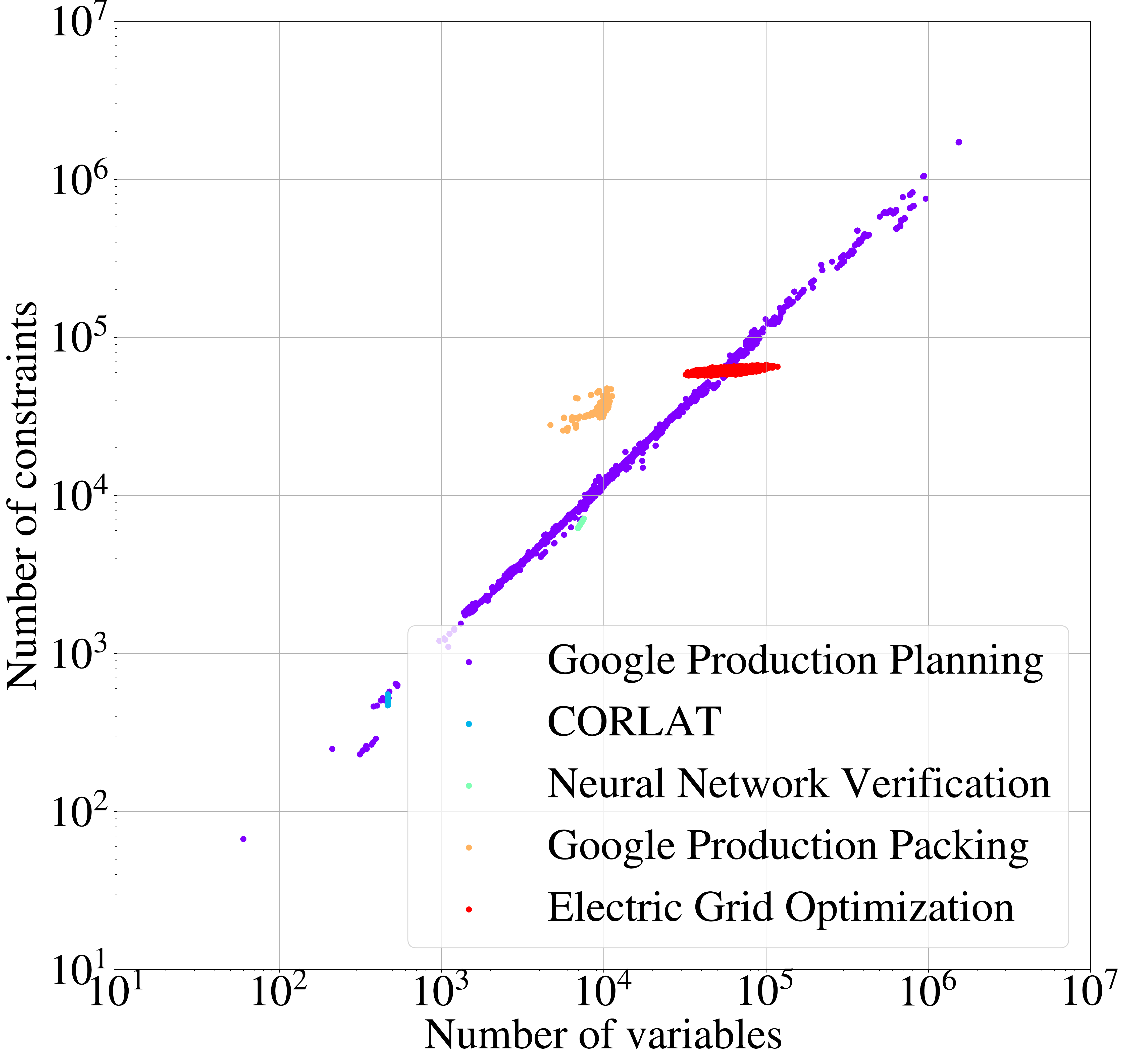}
        \hspace{0.1in}
        \includegraphics[width=0.48\textwidth]{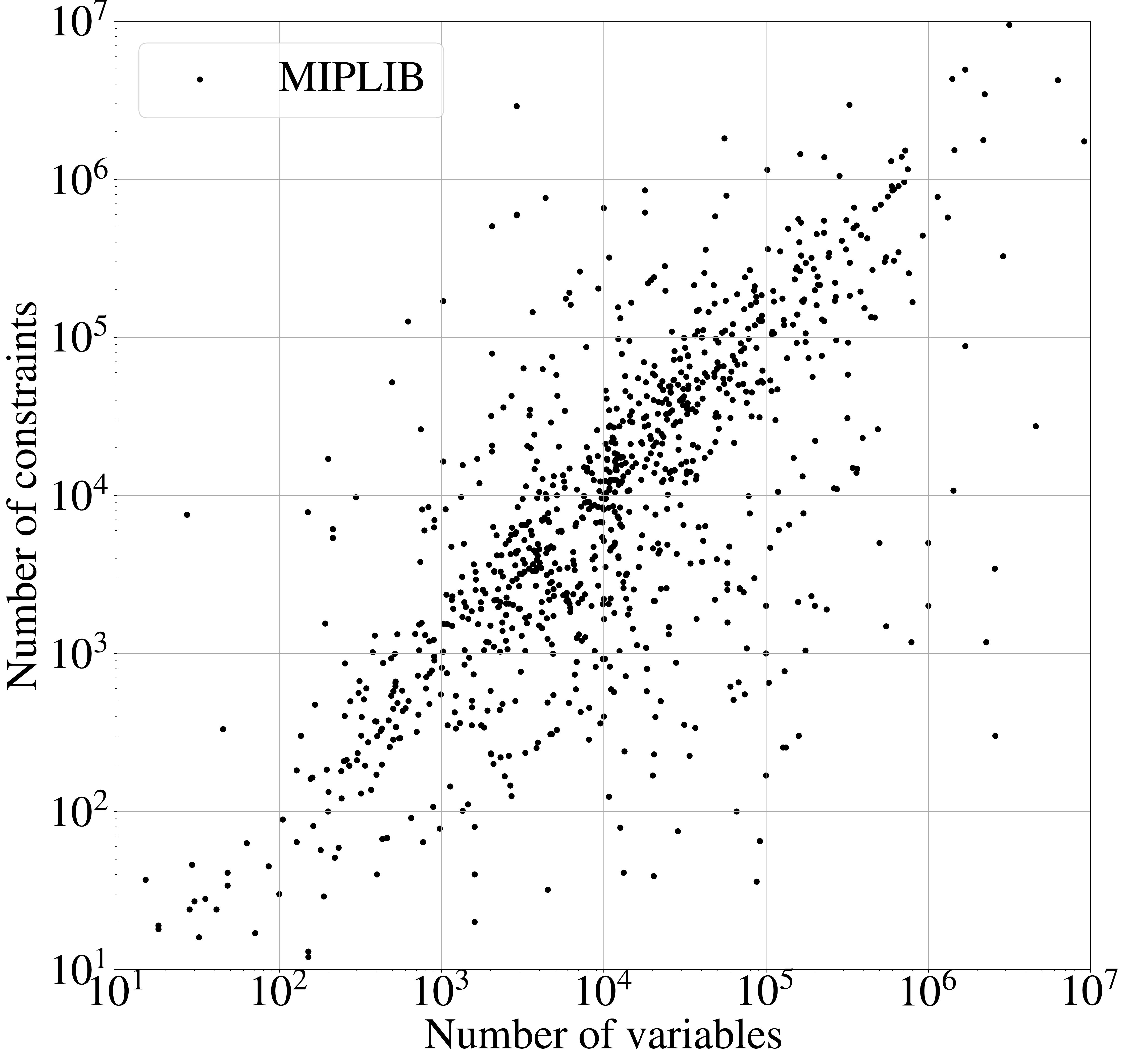}
    \end{tabular}
    \caption{
        {\small Original number of variables and constraints (before presolving) for the datasets used in our evaluation. Application-specific datasets are shown on the left, and MIPLIB is shown on the right.}}
    \label{fig:raw_mip_sizes}
\end{figure}

\begin{figure}
    \centering
    \begin{tabular}{cc}
        \includegraphics[width=0.48\textwidth]{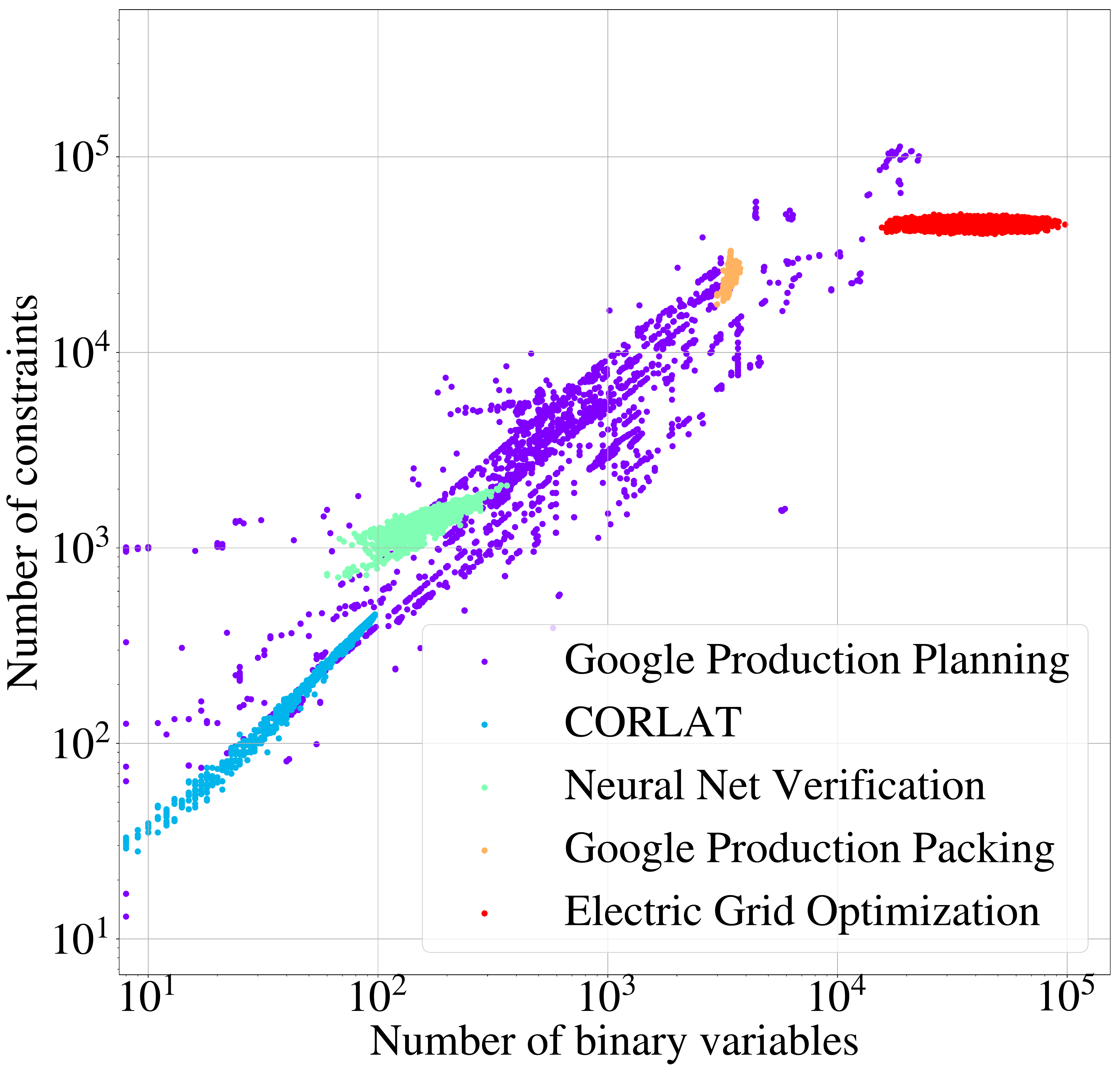}
        \hspace{0.1in}
        \includegraphics[width=0.48\textwidth]{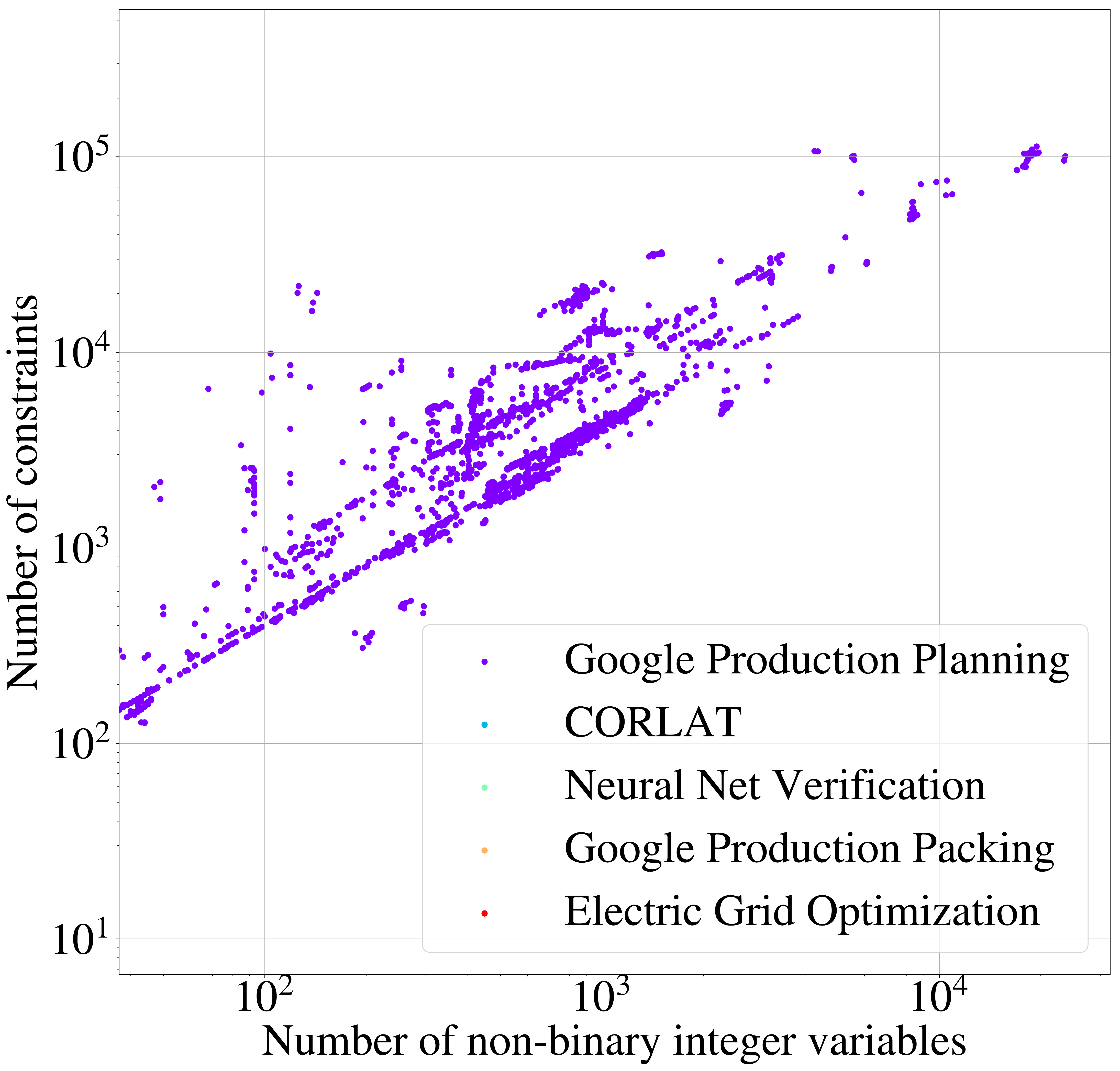}
    \end{tabular}
    \caption{
        {\small Number of binary variables (left) and non-binary integer variables (right) versus number of constraints after presolving using SCIP 7.0.1 for the application-specific datasets used in our evaluation.}
    }
    \label{fig:presolved_bin_nonbinint_mip_sizes}
\end{figure}

Figures~\ref{fig:raw_mip_sizes} and \ref{fig:presolved_bin_nonbinint_mip_sizes} show the MIP sizes for the datasets used in our evaluation with and without presolving using SCIP 7.0.1. Note that, among the application-specific datasets, only Google Production Planning contains non-binary integer variables. It is also the most heterogeneous (along with MIPLIB) in terms of instance sizes. Table~\ref{tab:dataset_stats} summarizes the characteristics of those datasets before and after presolving with the SCIP 7.0.1 solver. The number of constraints and variables ranges different orders of magnitude across the datasets. It is worth noting that during presolving some instances might be deemed infeasible and, hence, dropped from the dataset.

\begin{figure}
    \centering
    \includegraphics[width=\textwidth]{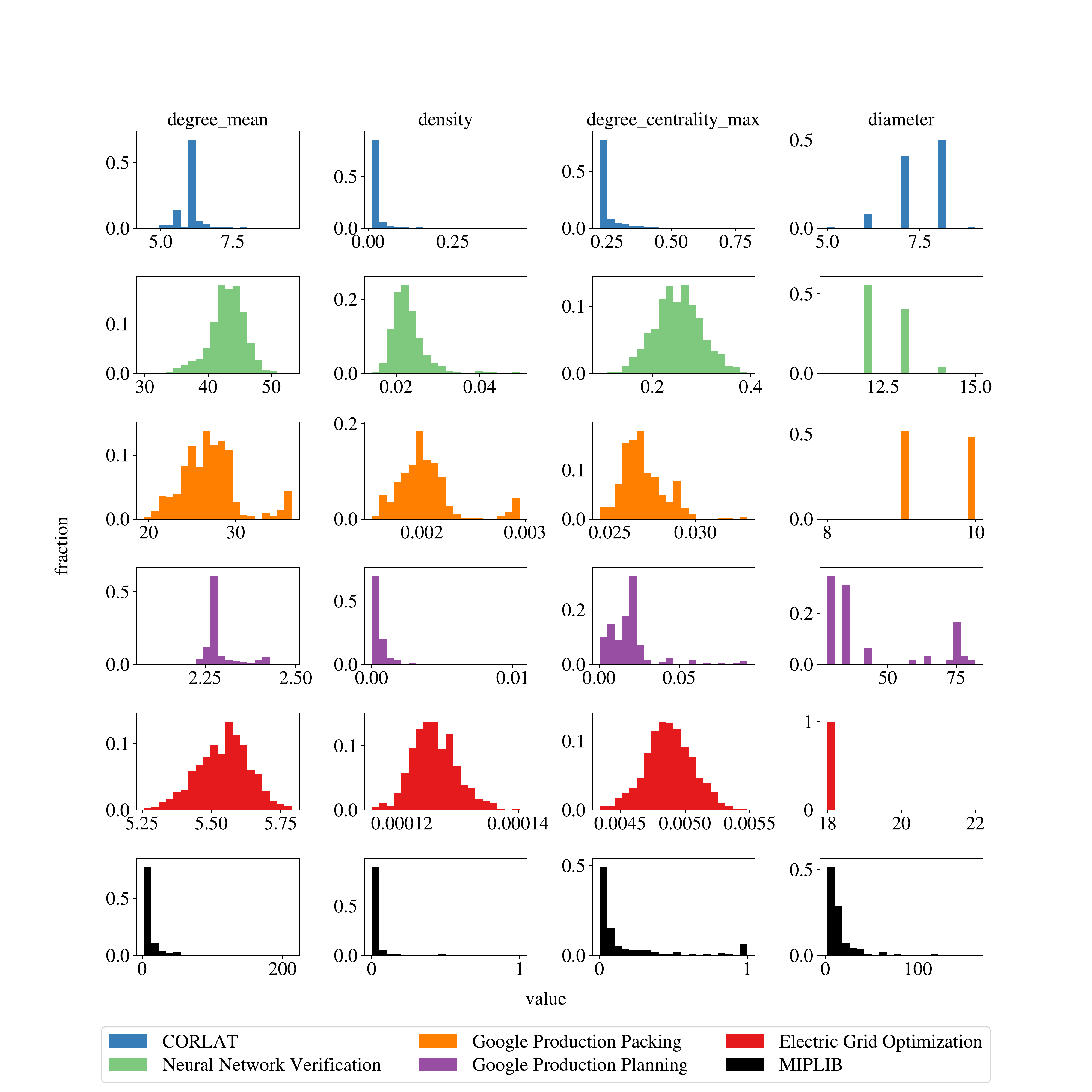}
    \caption{
        {\small Network statistics for the presolved datasets.}
    }
    \label{fig:presolved_network_stats}
\end{figure}

We additionally report some key network statistics for the bipartite graphs that represent the resulting MIPs of the presolved instances. These include: a) the average degree across all nodes, b) the density of the bipartite graph, i.e. the number of present edges over the number of all possible edges between the nodes, c) the maximum degree centrality for the variables nodes, and d) the diameter of the graph, i.e. the longest among all shortest paths between any pair of nodes. The normalized histograms for the graph characteristics are presented in Figure~\ref{fig:presolved_network_stats}.

As mentioned in the main paper, the above datasets were randomly split into training, validation and test sets with sizes 70\%, 15\% and 15\%, respectively. This split is handled differently for MIPLIB as it constitutes an externally defined `Benchmark Set'\footnote{\url{https://miplib.zib.de/tag_benchmark.html}}. This set contains 240 MIPs on which benchmarking results are reported for commercial and open source solvers and, therefore, constitutes our test set. The remaining instances in MIPLIB2017 and MIPLIB2010 correspond to our training and validation sets, respectively. Specifically, we use as the training set those instances from the MIPLIB2017 `Collection Set'\footnote{\url{https://miplib.zib.de/tag_collection.html}} which do not belong to the `Benchmark Set'. Similarly, we use as validation set the instances from MIPLIB2010 which do not belong to MIPLIB2017. This process yields 825 training, 155 validation and 240 test non-overlapping MIPs. The CORLAT dataset is available at \url{https://bitbucket.org/mlindauer/aclib2/src/master/}. The MIPLIB dataset (with the aforementioned tags) is available at \url{https://miplib.zib.de/download.html}. The Neural Network Verification dataset is available at \url{https://github.com/deepmind/deepmind-research/tree/master/neural_mip_solving}.

\subsection{Details of Calibrated Time}
\label{appendix:calibrated_time_results}
As explained in section \ref{sec:evaluation}, we use \emph{calibrated time} to accurately measure running time of evaluation solve tasks on a shared, heterogeneous compute cluster. For each solve task, we periodically solve a small \emph{calibration MIP} on a different thread as the solve task on the same machine. We define the speed of the machine to be
\begin{equation}
    \text{Speed} = \frac{1}{\text{Wall clock time to solve calibration MIP}}.
\end{equation}
For each periodic measurement of calibrated time, we estimate the speed $K$ times and use the average. $K$ is set to be the number of samples needed to estimate mean speed with 95\% confidence, with a minimum of 3 samples and a maximum of 30. The elapsed calibrated time $\Delta t_{\text{calibrated}}$ since the last measurement is
\begin{equation}
\label{eqn:calibrated_time}
    \Delta t_{\text{calibrated}} = \text{Speed} \times \Delta t_{\text{wallclock}},
\end{equation}
where $\Delta t_{\text{wallclock}}$ is the elapsed wallclock time since the last measurement. We use the MIP named \emph{vpm2} from MIPLIB2003 \cite{miplib2003} as the calibration MIP.

Note that the above definition of calibrated time does not have a time unit. Instead it is (in effect) a count of the calibrated MIP solves during the evaluation solve task. To give it a unit of seconds, one can choose a reference machine with respect to which we want to report evaluation solve times, accurately measure the calibration MIP's solve time on it (without other tasks interfering), and multiply the calibrated time in equation \ref{eqn:calibrated_time} by the reference machine's estimated calibration MIP solve time. The resulting quantity has a unit of seconds. It can be interpreted as the time the evaluation solve task would have taken if it ran on the reference machine. We select Intel Xeon 3.50GHz CPU with 32GB RAM as the reference machine. We estimate the mean solving time of the calibration MIP \emph{vpm2} for SCIP 7.0.1 on a single core to be 1.989 seconds (estimated from 1000 samples, collected after a warm up solve which is discarded). All the calibrated time results in the paper are expressed with respect to the reference machine in seconds.

To illustrate the benefit of calibrated time, table \ref{tab:calibrated_time_coeff_of_variation} presents results for four instances from MIPLIB 2017 \citep{gleixner2019miplib} comparing 1) wallclock time measurement on a shared cluster machine, 2) calibrated time measurement on a shared cluster machine, and 3) wallclock time measurement on the reference machine. We solve a MIP to optimality using SCIP 7.0.1 1000 times for each of these three settings. We then compute the \emph{coefficient of variation} as the sample standard deviation divided by the sample mean for the 1000 measurements. Calibrated time reduces the coefficient of variation by about 30$\times$, 8.3$\times$, 7$\times$, and 1.5$\times$, on \emph{air05}, \emph{n5-3}, \emph{swath1}, and \emph{dano3\_3}, respectively, compared to wallclock time measurements, and brings it closer to that of the reference machine.

\begin{table}
\caption{Coefficient of variation for different time measurements on MIPs from MIPLIB 2017.}
\begin{adjustbox}{center}
{\small
    \begin{tabular}{cccc}
    \toprule
        MIP & Shared cluster & Shared cluster & Reference machine \\
            & wallclock time & calibrated time & wallclock time \\
            & coeff. of var. & coeff. of var. & coeff. of var. \\
        \midrule
        \emph{air05} & 0.979 & 0.032 & 0.006 \\
        \emph{dano3\_3} & 0.224 & 0.144 & 0.006 \\
        \emph{n5-3} & 0.318 & 0.038 & 0.007 \\
        \emph{swath1} & 0.623 & 0.087 & 0.010 \\
    \bottomrule
    \end{tabular}
    }
\end{adjustbox}
\label{tab:calibrated_time_coeff_of_variation}
\end{table}